\documentclass[11 pt,draft]{article}
\title{An infinite presentation of the Torelli group}
\author{Andrew Putman}
\usepackage{amsmath}
\usepackage{amssymb}
\usepackage{amsthm}
\usepackage{epsfig}
\usepackage[margin=1.25in]{geometry}
\usepackage{amsfonts}
\usepackage[font=small,format=plain,labelfont=bf,up,textfont=it,up]{caption}
\usepackage{amscd}
\usepackage{pinlabel}
\usepackage{stmaryrd}
\usepackage{type1cm}
\usepackage{mathptmx}  %roman=Times, math=Times where possible

\theoremstyle{plain}
\newtheorem{theorem}{Theorem}[section]
\newtheorem{proposition}[theorem]{Proposition}
\newtheorem{lemma}[theorem]{Lemma}
\newtheorem{corollary}[theorem]{Corollary}

\newtheorem{claim}{Claim}
\newtheorem{step}{Step}

\newcommand\BeginClaims{\setcounter{claim}{0}}
\newcommand\BeginClaimProof{\begin{proof}[{Proof of Claim \theclaim}]}
\newcommand\EndClaimProof{\end{proof}}
\newcommand\BeginSteps{\setcounter{step}{0}}
\newcommand\BeginCases{\setcounter{case}{0}}

\theoremstyle{definition}
\newtheorem{definition}[theorem]{Definition}
\newtheorem{example}[theorem]{Example}

\newtheorem{case}{Case}

\theoremstyle{remark}
\newtheorem*{remark}{Remark}
\newtheorem*{remarks}{Remarks}
\newtheorem*{warning}{Warning}

\newtheorem*{alternatederivation}{Alternate Derivation}

\DeclareMathOperator{\Ker}{ker}

\DeclareMathOperator{\Mod}{Mod}
\newcommand\Torelli{\text{${\mathcal I}$}}
\DeclareMathOperator{\Sp}{Sp}

\newcommand\Curves{\text{$\mathcal{C}$}}
\newcommand\Sphere[1]{\text{$\mathbf{S}^{#1}$}}
\newcommand\Ball[1]{\text{$\mathbf{B}^{#1}$}}

\newcommand\Z{\text{$\mathbb{Z}$}}

\DeclareMathOperator{\HH}{H}

\DeclareMathOperator{\Max}{max}

\DeclareMathOperator{\Aut}{Aut}

\DeclareMathOperator{\Interior}{Int}
\newcommand\Span[1]{\text{$\langle \text{#1} \rangle$}}
\newcommand\CaptionSpace{\hspace{0.2in}}
\DeclareMathOperator{\Dim}{dim}

\newcommand\Figure[3]{
\begin{figure}[t]
\centering
\centerline{\psfig{file=#2,scale=60}}
\caption{#3}
\label{#1}
\end{figure}}

\newcommand\BigFreeProd{\mathop{\text{\Huge{$\ast$}}}}
\newcommand\ModCurves{\text{$\mathcal{MC}$}}
\DeclareMathOperator{\nosep}{nosep}
\newcommand\DComm[1]{\text{$\llbracket #1 \rrbracket$}}
\newcommand\PtPsh[1]{\text{$\overline{\text{Push}}(#1)$}}
\newcommand\LPtPsh[1]{\text{$\text{Push}(#1)$}}

\DeclareMathOperator{\Star}{star}
\DeclareMathOperator{\Link}{link}
\newcommand\CNoSep{\mathcal{C}^{\text{nosep}}}
\newcommand\Lines{\text{$\mathcal{L}$}}

\DeclareMathOperator{\Rank}{rk}

\newcommand\GeomI{\text{$i_{\text{geom}}$}}
\newcommand\AlgI{\text{$i_{\text{alg}}$}}

\begin{document}

\maketitle

\begin{abstract}
In this paper, we construct an infinite presentation of the Torelli subgroup of the mapping
class group of a surface whose generators consist of the set of
all ``separating twists'', all ``bounding pair maps'', and
all ``commutators of simply intersecting pairs'' and whose relations all come from
a short list of topological configurations of these generators on the surface.  Aside from a few obvious ones, all of
these relations come from a set of embeddings of groups derived from surface groups into
the Torelli group.  In the process of analyzing these embeddings, we derive a novel
presentation for the fundamental group of a closed surface whose generating set is the
set of {\em all} simple closed curves.  
\end{abstract}

\section{Introduction}
\label{section:introduction}

Let $\Sigma_{g}$ be a closed genus $g$ surface and $\Mod_{g}$ be the {\em mapping class
group} of $\Sigma_g$, that is, the group of homotopy classes of orientation-preserving diffeomorphisms
of $\Sigma_{g}$.  The action of $\Mod_{g}$ on $\HH_1(\Sigma_g;\Z)$ preserves the algebraic intersection form, so it
induces a representation $\Mod_{g} \rightarrow \Sp_{2g}(\Z)$.  The kernel 
$\Torelli_{g}$ of this representation is known as the {\em Torelli group}.  It
plays an important role in both low-dimensional topology and algebraic geometry.  See
\cite{JohnsonSurvey} for a survey of $\Torelli_{g}$, especially the remarkable work of Dennis Johnson. 

Despite the Torelli group's importance, little is known about its combinatorial group theory.  Generators
for $\Torelli_{g}$ were first found by Birman and Powell \cite{BirmanSiegel, PowellTorelli} (see below).  Later, 
Johnson \cite{JohnsonFinite} constructed a finite generating set for $\Torelli_{g}$ for
$g \geq 3$, while McCullough and Miller \cite{McCulloughMiller} proved that $\Torelli_{2}$ is not finitely
generated.  The investigation of the genus $2$ case was completed by Mess \cite{MessTorelli}, who proved that $\Torelli_{2}$ 
is an infinitely generated free group, though no explicit free generating set is known.  However, the basic question
of whether $\Torelli_{g}$ is ever finitely presented for $g \geq 3$ remains open.  

In this paper, we construct an infinite presentation for $\Torelli_{g}$ whose
generators and relations have simple topological interpretations.  This is not the first presentation
of the Torelli group in the literature -- another appears in a paper of Morita and Penner \cite{MoritaPenner}.
However, while their generators and relations have nice interpretations in terms of a certain triangulation of
Teichm\"{u}ller space, they are topologically and group-theoretically extremely complicated.  Indeed, their
generating set contains infinitely many copies of {\em every} element of the Torelli group.  Our methods and
perspective are very different from theirs. 

\Figure{figure:generators}{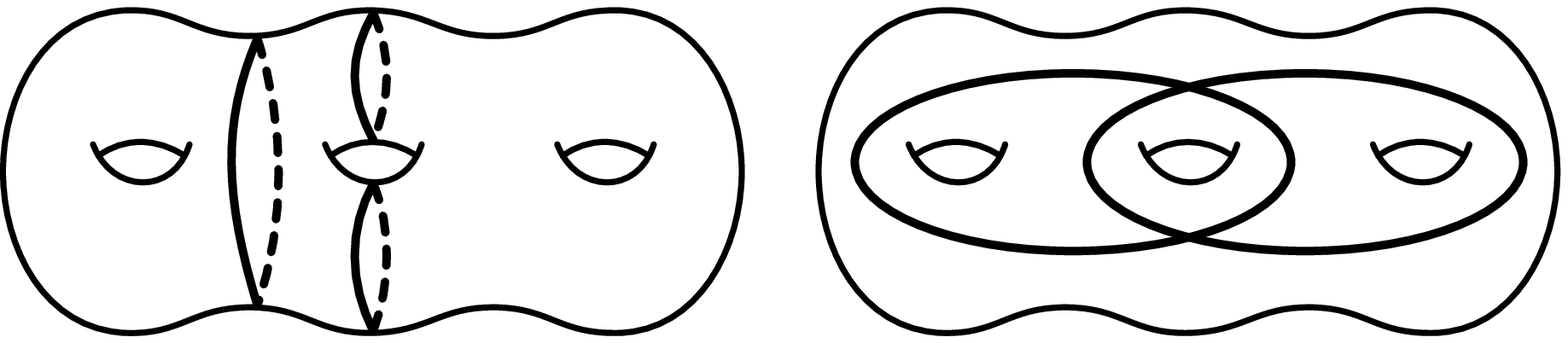}{a. A separating curve $x_1$ and a bounding pair $\{x_2,x_3\}$ \CaptionSpace b. A simply
intersecting pair $\{x_4,x_5\}$}

\paragraph{Generators.}
Letting $T_{\gamma}$ be the 
right Dehn twist about a simple closed curve $\gamma$, 
our generators are all mapping classes of the following types.
\begin{enumerate}
\item Let $\gamma$ be a simple closed curve that separates the surface (for instance,
      the curve $x_1$ in Figure \ref{figure:generators}.a).  Then it is not hard
      to see that $T_\gamma \in \Torelli_{g}$.  These are known as {\em separating twists}.
\item Let $\{\gamma_1,\gamma_2\}$ be a pair of non-isotopic disjoint homologous
      curves (for instance, the pair of curves $\{x_2,x_3\}$ from Figure \ref{figure:generators}.a).  Then
      $T_{\gamma_1}$ and $T_{\gamma_2}$ map to the same element of $\Sp_{2g}(\Z)$, so
      $T_{\gamma_1} T_{\gamma_2}^{-1} \in \Torelli_{g}$.  These are known as {\em bounding
      pair maps}.  We will denote them by $T_{\gamma_1,\gamma_2}$.
\item Let $\{\gamma_1,\gamma_2\}$ be a pair of curves whose algebraic
      intersection number is $0$.  Then the images of $T_{\gamma_1}$ and $T_{\gamma_2}$
      in $\Sp_{2g}(\Z)$ commute, so
      $[T_{\gamma_1},T_{\gamma_2}] \in \Torelli_{g}$.  We will make use of
      such commutators for simple closed curves $\gamma_1$ and $\gamma_2$ whose
      geometric intersection number is $2$ (for instance, the pair of curves
      $\{x_4,x_5\}$ from Figure \ref{figure:generators}.b).  We will call these {\em commutators of simply intersecting
      pairs} and denote them by $C_{\gamma_1,\gamma_2}$.
\end{enumerate}

\begin{remarks}
\mbox{}
\begin{itemize}
\item The fact that $\Torelli_{g}$ is generated by separating twists
      and bounding pair maps follows from work of Birman and
      Powell (\cite{BirmanSiegel,PowellTorelli}; see also \cite{PutmanCutPaste} for a different proof, as well as
      generalizations)
\item {\bf Warning :} Traditionally, the curves in a bounding pair are required to
      be nonseparating; however, to simplify our statements we allow them to be separating.
\item Commutators of simply intersecting pairs are not needed to generate $\Torelli_{g}$,
      but their presence greatly simplifies our relations.  We remark
      that the expression of a mapping class as a commutator of a simply intersecting pair
      is not unique; see Example \ref{example:nonunique} for an example.
\end{itemize}
\end{remarks}

\paragraph{Relations.}
Our relations are as follows; a more detailed description follows.
\begin{enumerate}
\item The {\em formal relations} \eqref{F.1}-\eqref{F.8}.  An example
      is $T_{\gamma_1,\gamma_2} = T_{\gamma_2,\gamma_1}^{-1}$.  
\item Two families of relations (the {\em lantern relations} and the {\em crossed
      lantern relations}) that arise from easy identities among various ways of ``dragging subsurfaces around''.
\item Two families of relations (the {\em Witt--Hall relations} and the {\em commutator shuffle
      relations}) that arise from easy identities among various ways of ``dragging bases of handles around''.  
\end{enumerate}

\Figure{figure:handledrag}{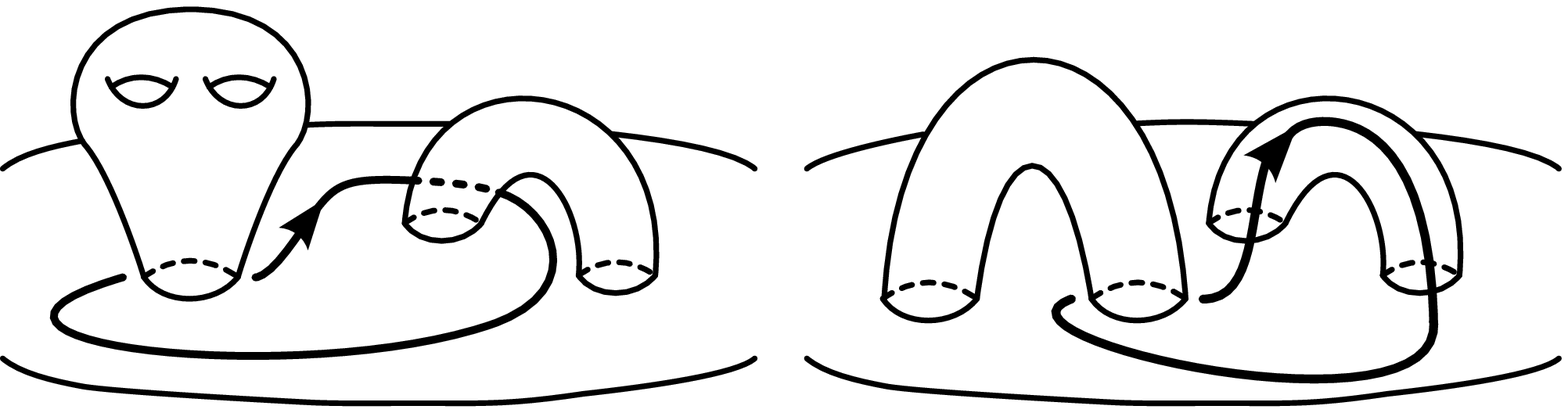}{a. Dragging a copy of $\Sigma_{h_2,1}$ around a
curve $\gamma$. \CaptionSpace b. Dragging the end of a handle around a curve $\gamma$}

\paragraph{Formal relations.}
These relations are formal in the sense that they are either
immediate consequences of the standard expressions of our generators as products of Dehn twists
or are consequences of the conjugation relation $f T_x f^{-1} = T_{f(x)}$, where $x$ is a simple closed
curve and $f$ is a mapping class.  The first three are immediate, and
are true for any curves $x_1$, $x_2$, and $x_3$ so that the expressions make sense.
\begin{align*}
T_{x_1,x_2}             &= T_{x_2, x_1}^{-1},            \tag{F.1}\label{F.1} \\
C_{x_1,x_2} &= C_{x_2, x_1}^{-1},                        \tag{F.2}\label{F.2} \\
T_{x_1,x_2} T_{x_2,x_3} &= T_{x_1,x_3}.                  \tag{F.3}\label{F.3} \\
\end{align*}
Next, if $\{x_1,x_2\}$ is a bounding pair so that both $x_1$ and $x_2$ are separating curves, we need
\begin{equation*}
T_{x_1,x_2} = T_{x_1} T_{x_2}^{-1}. \tag{F.4}\label{F.4}\\
\end{equation*}
If $\{x_1,x_2\}$ is a bounding pair and $\{x_3,x_2\}$ is a simply intersecting pair so that $x_1$ and
$x_3$ are disjoint, we need
\begin{equation*}
T_{x_1,T_{x_3}^{-1}(x_2)} = C_{x_3,x_2} T_{x_1, x_2}. \tag{F.5}\label{F.5}\\
\end{equation*}
Finally, we will also need the following conjugation relations.  In them,
$A$ is any generator and $x$, $x_1$,
and $x_2$ are any curves so that the expressions make sense.
\begin{align*}
A T_{x} A^{-1}       &= T_{A(x)},          \tag{F.6}\label{F.6} \\
A T_{x_1,x_2} A^{-1} &= T_{A(x_1),A(x_2)}, \tag{F.7}\label{F.7} \\
A C_{x_1,x_2} A^{-1} &= C_{A(x_1),A(x_2)}. \tag{F.8}\label{F.8} \\
\end{align*}

\paragraph{Lantern and crossed lantern relations.} 
Letting $\Sigma_{h,n}$ denote a genus $h$ surface
with $n$ boundary components, consider a subsurface $S$ of $\Sigma_g$ that is homeomorphic
to $\Sigma_{h_1,1}$ for some $h_1 < g-1$.  The closure $S'$ of the complement of $S$ is
then homeomorphic to $\Sigma_{h_2,1}$ with $h_1 + h_2 = g$ and $h_2 > 1$.  Informally, we
can obtain elements of $\Mod_{g}$ by ``dragging'' $S$  around a curve $\gamma$ in $S'$ (see Figure \ref{figure:handledrag}.a).
Using results of Birman \cite{BirmanExactSequence} and Johnson \cite{JohnsonFinite}, we will formalize this
and show that it yields an injection $i:\pi_1(U\Sigma_{h_2}) \rightarrow \Mod_{g}$, where $U\Sigma_{h_2}$ is the unit
tangent bundle of $\Sigma_{h_2}$ (see \S \ref{section:subsurfacedragprelim} for the details; we need the unit 
tangent bundle because $S$ may ``rotate'' as it is being dragged).  Moreover, $i(\pi_1(U\Sigma_{h_2})) \subset \Torelli_{g}$.
If $b = \partial S$, then $i$ of the loop around the fiber (with an appropriate choice of orientation) is $T_b$.

We can thus find relations in $\Torelli_{g}$ from relations in $\pi_1(U\Sigma_{h_2})$.  It will
be easier to describe these relations in terms of the group $\pi_1(\Sigma_{h_2})$.  Let
$\rho : \pi_1(U\Sigma_{h_2}) \rightarrow \pi_1(\Sigma_{h_2})$ be the projection.  We thus have
an exact sequence
$$1 \longrightarrow \Z \longrightarrow \pi_1(U\Sigma_{h_2}) \stackrel{\rho}{\longrightarrow} \pi_1(\Sigma_{h_2}) \longrightarrow 1.$$
Since $h_2 > 1$, this exact sequence does not split.  However, in \S \ref{section:birmanexactsequence} we
will give a procedure which takes any nontrivial $\gamma \in \pi_1(\Sigma_{h_2})$ that can be represented by
a simple closed curve and produces a well-defined $\tilde{\gamma} \in \pi_1(U\Sigma_{h_2})$ so that 
$\rho(\tilde{\gamma}) = \gamma$.
Define $\LPtPsh{\gamma} = i(\tilde{\gamma}) \in \Torelli_{g}$.  We will prove that $\LPtPsh{\gamma}$ is
a bounding pair map.

Let $\gamma_1,\ldots,\gamma_n \in (\pi_1(\Sigma_{h_2}) \setminus \{1\})$ be elements all of which
can be represented by simple closed curves and which satisfy $\gamma_1 \cdots \gamma_n = 1$.
If $\tilde{\gamma}_i$ is the aforementioned lift of $\gamma_i$ to $\pi_1(U\Sigma_{h_2})$ for $1 \leq i \leq n$, then
$\tilde{\gamma}_1 \cdots \tilde{\gamma}_n$ is equal to some power of the loop around the fiber.  We conclude that
for some $k \in \Z$ we have the following relation in $\Torelli_{g}$ :
$$\LPtPsh{\gamma_n} \cdots \LPtPsh{\gamma_1} = T_b^k.$$
The order of the product on the left hand
side is reversed because fundamental group elements are composed via concatenation order while
mapping classes are composed via functional order.

We thus need to find all relations between simple closed curves in $\pi_1(\Sigma_{h_2})$.  This
is provided by the following theorem.

\begin{theorem}
\label{theorem:surfacegroup}
Let $\Gamma$ be the abstract group whose generating set consists
of the symbols
$$\{s_\gamma \text{ $|$ $\gamma \in (\pi_1(\Sigma_g) \setminus \{1\})$ is represented by a simple closed curve}\}$$
and whose relations are $s_{\gamma} s_{\gamma^{-1}}=1$ for all simple closed curves $\gamma$,
\begin{equation*}
s_x s_y s_z = 1 \tag{$\overline{\text{L}}$}\label{OL} \\
\end{equation*}
for all curves $x$, $y$, and $z$ arranged like the curves in Figure \ref{figure:surfacerelations}.a, and
\begin{equation*}
s_x s_y = s_z \tag{$\overline{\text{CL}}$}\label{OCL}\\
\end{equation*}
for all curves $x$, $y$, and $z$ arranged like the curves in Figure \ref{figure:surfacerelations}.b.  Then
the natural map $\Gamma \rightarrow \pi_1(\Sigma_g)$ is an isomorphism.
\end{theorem}

\noindent
We will see that via the above procedure the relation \eqref{OL} lifts to the well-known {\em lantern relation} \eqref{L}
$$T_{\tilde{z}_1,\tilde{z}_2} T_{\tilde{y}_1,\tilde{y}_2} T_{\tilde{x}_1,\tilde{x}_2} = T_b$$
depicted in Figure \ref{figure:surfacerelations}.c, while the relation \eqref{OCL} lifts
to the relation
$$T_{\tilde{y}_1,\tilde{y}_2} T_{\tilde{x}_1,\tilde{x}_2} = T_{\tilde{z}_1,\tilde{z}_2}$$
depicted in Figure \ref{figure:surfacerelations}.d.  We will call this the {\em crossed
lantern relation} \eqref{CL}.

\Figure{figure:surfacerelations}{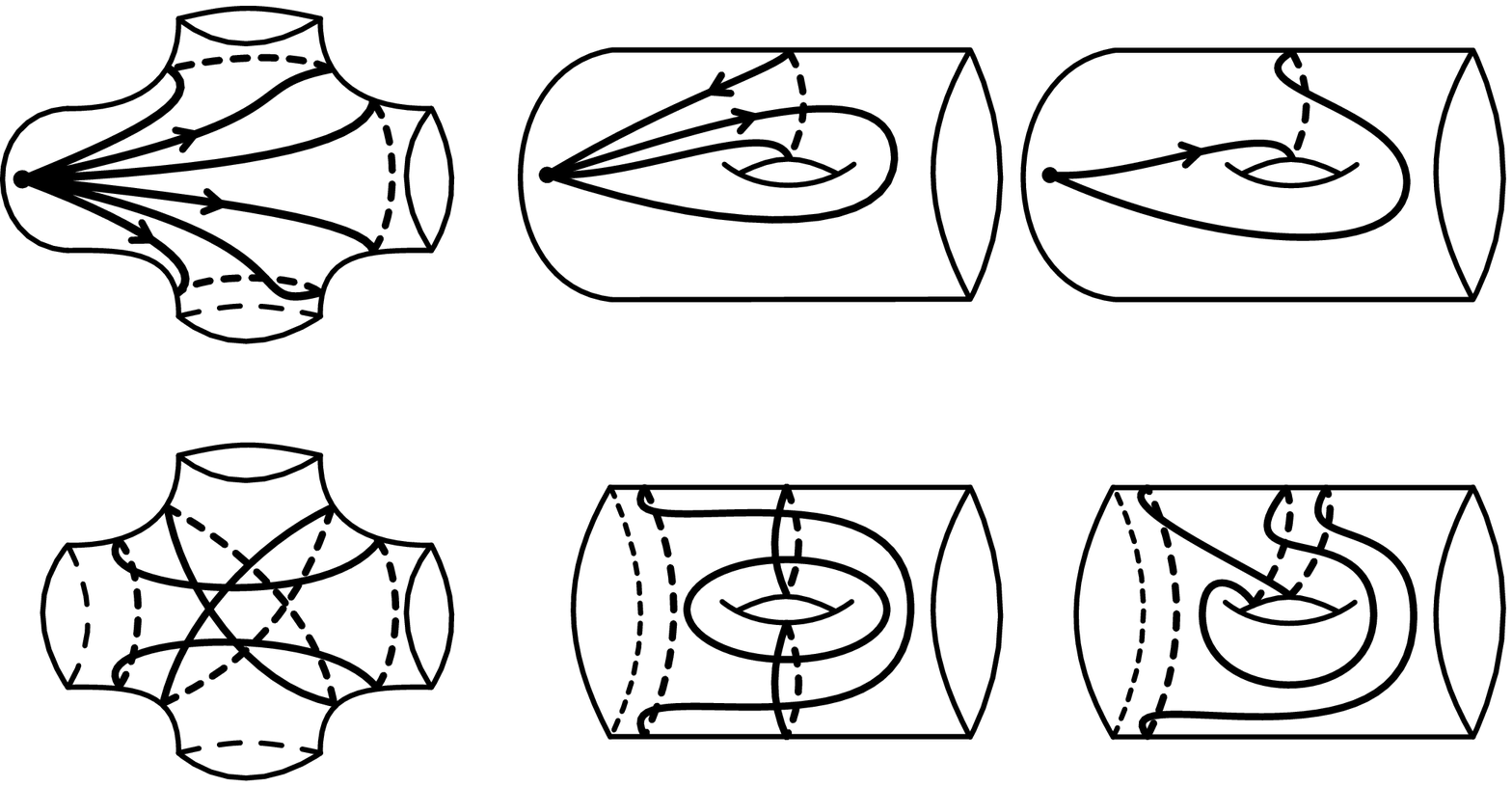}{a. Relation $\overline{L}$ \CaptionSpace b. Relation $\overline{CL}$ \CaptionSpace
c. The lantern relation $T_{\tilde{z}_1,\tilde{z}_2} T_{\tilde{y}_1,\tilde{y}_2} T_{\tilde{x}_1,\tilde{x}_2} = T_b$ 
\CaptionSpace d. The crossed lantern relation
$T_{\tilde{y}_1,\tilde{y}_2} T_{\tilde{x}_1,\tilde{x}_2} = T_{\tilde{z}_1,\tilde{z}_2}$}

\paragraph{Witt--Hall and commutator shuffle relations.}
Let $H$ be a handle on $\Sigma_g$; i.e.\ an
embedded annulus that does not separate the surface.  The closure of the 
complement of $H$ is homeomorphic to $\Sigma_{g-1,2}$.  In
a manner similar to the previous case, dragging
one end of $H$ around curves $\gamma$ on $\Sigma_{g-1,2}$ (see Figure \ref{figure:handledrag}.b)
yields an injection $j : \pi_1(U\Sigma_{g-1,1}) \rightarrow \Mod_{g}$.  

However, in this case we do {\em not} have $j(\pi_1(U\Sigma_{g-1,1})) \subset \Torelli_{g}$.  Using previous
results of the author (see \S \ref{section:handledragprelim}), we will show there is an isomorphism
$j^{-1}(\Torelli_{g}) \cong [\pi_1(\Sigma_{g-1,1}),\pi_1(\Sigma_{g-1,1})]$.  We thus
have an induced map $j' : [\pi_1(\Sigma_{g-1,1}),\pi_1(\Sigma_{g-1,1})] \rightarrow \Torelli_{g}$.  Throughout
the paper, we will say that two curves $x$ and $y$ in the fundamental group of a surface are
{\em completely distinct} if $x \neq y$ and $x \neq y^{-1}$.  We
we will then show that if $x,y \in \pi_1(\Sigma_{g-1,1})$ are completely distinct nontrivial elements that 
can be represented by simple closed curves that only
intersect at the basepoint, then $j'([x,y])$ has a simple expression in terms of our generators.
It follows that we can use commutator identities
between appropriate simple closed curves to obtain relations in $\Torelli_{g}$.  In what follows, we will
frequently use the observation that if $x,y \in \pi_1(\Sigma_{g-1,1})$ can be represented by simple closed
curves that only intersect at the basepoint and $z \in \pi_1(\Sigma_{g-1,1})$ is arbitrary, then $x^z$ and
$y^z$ can also be represented by simple closed curves that only intersect at the basepoint (here $x^z$ and
$y^z$ denote $z^{-1} x z$ and $z^{-1} y z$).

For the Witt--Hall relations, let $g_1,g_2,g_3 \in (\pi_1(\Sigma_{g-1,1}) \setminus \{1\})$ be elements so that for each of the
sets $\{g_1,g_2,g_3\}, \{g_1 g_2, g_3\} \subset \pi_1(\Sigma_{g-1,1})$,
the elements of the set can be represented by completely distinct
simple closed curves that only intersect at the basepoint.  Via the above procedure, we will 
use the Witt--Hall commutator identity
$$[g_1 g_2, g_3] = [g_1,g_3]^{g_2} [g_2,g_3]$$
to derive a family of relations \eqref{WH} which we will call the {\em Witt--Hall relations}.  

For the commutator shuffle relations, let $g_1,g_2,g_3 \in (\pi_1(\Sigma_{g-1,1}) \setminus \{1\})$ be completely distinct
elements which can
be realized by simple closed curves that only intersect at the basepoint.  Via the above procedure, we will 
we will use the easily-verified commutator identity
$$[g_1,g_2]^{g_3} = [g_3, g_1] [g_3, g_2]^{g_1} [g_1,g_2] [g_1,g_3]^{g_2} [g_2,g_3]$$
to obtain a family of relations \eqref{CS} that we will call
the {\em commutator shuffles}.  This final commutator identity may be viewed as a variant
of the classical Jacobi identity. 

\begin{remark}
For each Witt--Hall and commutator shuffle relation, the above procedure gives
a relation that is supported on a subsurface of $\Sigma_g$.  This subsurface
may be embedded in the surface in many different ways, and we will need all
relations come from such embeddings.  See the beginning of \S \ref{section:witthall}
for a precise description of this.
\end{remark}

\paragraph{Main theorem.}
We can now state our Main Theorem.

\begin{theorem}
\label{theorem:maintheoreminformal}
For $g \geq 2$, the group $\Torelli_{g}$ has a presentation
whose generators are the set of all separating twists, all bounding pair maps, and
all commutators of simply intersecting pairs and whose relations are the formal 
relations \eqref{F.1}-\eqref{F.8}, the lantern relations \eqref{L}, the crossed lantern relations
\eqref{CL}, the Witt--Hall relations \eqref{WH}, and the commutator shuffle relations \eqref{CS}.
\end{theorem}
\noindent
We also prove a similar statement for surfaces with boundary (see \S \ref{section:precisestatement}).

The proof of Theorem \ref{theorem:maintheoreminformal} is by induction on $g$.  The base case 
$g=2$ is derived from the theorem of Mess \cite{MessTorelli} mentioned above that says 
that $\Torelli_{2}$ is an infinitely generated free group.  For the inductive
step, the key is to show that $\Torelli_{g}$ has a presentation most of whose
relations ``live'' in the subgroups of $\Torelli_{g}$ stabilizing simple closed
curves (these subgroups are supported on ``simpler'' subsurfaces).  

The proof of this, like many constructions of group presentations, relies on the study
of a natural simplicial complex upon which the group acts.  We will use a suitable
modification of the nonseparating complex of curves, whose definition is as follows.

\begin{definition}
The {\em complex of curves} on $\Sigma_{g,n}$, denoted
$\Curves_{g,n}$, is the simplicial complex whose $(k-1)$-simplices
are sets $\{\gamma_1,\ldots,\gamma_k\}$ of distinct nontrivial isotopy classes of simple closed curves on $\Sigma_{g,n}$
that can be realized disjointly. 
The {\em nonseparating complex of curves} on $\Sigma_{g,n}$, denoted $\CNoSep_{g,n}$, is the subcomplex
of $\Curves_{g,n}$ whose
$(k-1)$-simplices are sets $\{\gamma_1,\ldots,\gamma_k\}$ of isotopy classes
that can be realized so that $\Sigma_{g,n} \setminus (\gamma_1 \cup \cdots \cup \gamma_k)$ is connected.
\end{definition}

\noindent
The complex of curves was introduced by Harvey \cite{HarveyComplex}, while the nonseparating complex
of curves was introduced by Harer \cite{HarerStability}.  We will usually omit the $n$ on $\Curves_{g,n}$
and $\CNoSep_{g,n}$ when it equals $0$.

Now, there are several standard methods for writing down a presentation from a group
action in terms of the stabilizers (see, e.g., the work of K. Brown \cite{BrownPresentation}).
However, we are unable to use these methods here,
as they all require an explicit fundamental domain
for the action, which seems quite difficult to pin down in our situation.
We instead use a theorem of the author (\cite{PutmanPresentation}; see Theorem \ref{theorem:presentation} below)
that allows us to derive presentations from group actions without identifying
a fundamental domain.  

The hypotheses of this theorem require that the quotient
of the simplicial complex by the group be $2$-connected.  Unfortunately, $\CNoSep_{g} / \Torelli_{g}$
is only $(g-2)$-connected (see Lemma \ref{lemma:modcurveslines} and Proposition \ref{proposition:main}), and
hence $\CNoSep_g$ does not work for the case $g = 3$.  Our solution is to attach additional
cells to $\CNoSep_{g}$ to increase the connectivity of its quotient by $\Torelli_{g}$.  The 
complex we make use of is as follows.  Denote by $\GeomI(\gamma_1,\gamma_2)$ the {\em geometric intersection number} 
of two simple closed curves $\gamma_1$ and $\gamma_2$, i.e.\ the minimum over all curves
$\gamma_1'$ and $\gamma_2'$ with $\gamma_i'$ isotopic to $\gamma_i$ for $1 \leq i \leq 2$ of the number
of points of $\gamma_1' \cap \gamma_2'$.

\Figure{figure:modifiedcomplex}{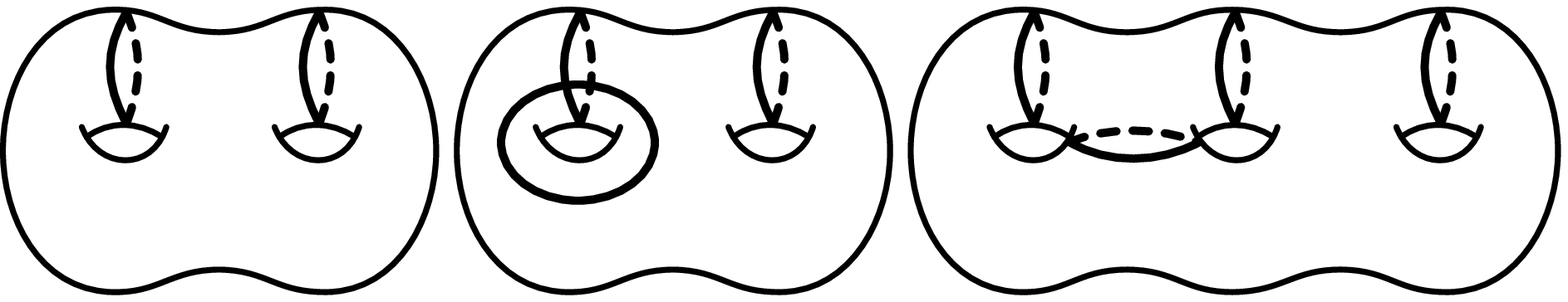}{a,b,c. Examples of the three kinds of simplices in $\ModCurves_g$}

\begin{definition}
The complex $\ModCurves_g$ is the simplicial complex whose $(k-1)$-simplices are sets
$\{\gamma_1,\ldots,\gamma_k\}$ of isotopy classes of simple closed nonseparating curves
on $\Sigma_{g}$ satisfying one of the following three conditions (for some ordering of the $\gamma_i$).
\begin{itemize}
\item The $\gamma_i$ are disjoint and $\gamma_1 \cup \cdots \cup \gamma_k$ does not separate
$\Sigma_{g}$ (see Figure \ref{figure:modifiedcomplex}.a).  
\item The $\gamma_i$ satisfy
$$\GeomI(\gamma_i,\gamma_j) = \begin{cases}
        1 & \text{if $(i,j)=(1,2)$} \\
        0 & \text{otherwise}
\end{cases}$$
and $\gamma_1 \cup \cdots \cup \gamma_k$ does not separate $\Sigma_{g}$ (see Figure \ref{figure:modifiedcomplex}.b).
\item The $\gamma_i$ are disjoint, $\gamma_1 \cup \gamma_2 \cup \gamma_3$ cuts off a copy of $\Sigma_{0,3}$ 
from $\Sigma_{g}$, and $\{\gamma_1,\ldots,\gamma_k\} \setminus \{\gamma_1\}$ is a standard simplex (see 
Figure \ref{figure:modifiedcomplex}.c).  
\end{itemize}
\end{definition}

\noindent
Our main result about $\ModCurves_g$ (Proposition \ref{proposition:curvestorelliconnected} below)
says that $\ModCurves_g / \Torelli_{g}$ is $(g-1)$-connected.  In particular, it is $2$-connected for $g=3$.

\paragraph{History and comments.}
Three additional results concerning presentations of the Torelli group should be mentioned.  First,
Krsti\'{c} and McCool \cite{KrsticMcCoolIAn} have proven that the analogue of the Torelli group
in $\Aut(F_n)$ is not finitely presentable for $n=3$.  Second, using algebreo-geometric methods,
Hain \cite{HainMalcev} has computed a finite presentation for the Malcev Lie algebra of
$\Torelli_{g}$ for $g \geq 6$.  Finally, in addition to their infinite presentation
of the Torelli group, Morita and Penner \cite{MoritaPenner} used
Johnson's finite generating set for the Torelli group to give a finite presentation of the fundamental {\em groupoid}
of a certain cell decomposition of the quotient of Teichm\"{u}ller space by the Torelli group.

As far as relations in the Torelli group go, Johnson's paper \cite{JohnsonFinite} contains
a veritable zoo of relations, most of which are derived from clever combinations of lantern
relations in the mapping class group.  An excellent discussion of these relations, plus some
generalizations of them, can be found in Brendle's unpublished thesis \cite{BrendleThesis}.  The
rest of our relations seem to be new, though it is unclear which of them can be derived from
Johnson's relations.

We finally wish to draw attention to a paper of Gervais \cite{GervaisPresentation} that
constructs an infinite presentation for the whole mapping class group using the set of all Dehn
twists as generators.  Gervais's presentation was later simplified by Luo \cite{LuoPresentation}.

\paragraph{Outline.}
We begin in \S \ref{section:preliminaries} with a review of the Birman
exact sequence together with some basic group theory.  Next, in \S \ref{section:relations} 
we derive the nonformal relations in our presentation.  The proof of Theorem 
\ref{theorem:maintheoreminformal} is in \S \ref{section:maintheorem}.  This proof
depends on two propositions that are proven in \S \ref{section:gammaexactsequences}
and \S \ref{section:curvestorelliconnected}.  

\paragraph{Conventions and notation.}
All homology groups will have $\Z$ coefficients.
Throughout this paper, we will systematically confuse simple closed curves with their homotopy
classes.  Hence (based/unbased) curves are said to be simple closed curves if they are (based/unbased)
homotopic to simple closed curves, etc.  If $\gamma_1$ and $\gamma_2$ are two simple closed curves,
then $\GeomI(\gamma_1,\gamma_2)$ will denote the {\em geometric intersection number} of $\gamma_1$
and $\gamma_2$; i.e.\ the minimum over all curves
$\gamma_1'$ and $\gamma_2'$ with $\gamma_i'$ isotopic to $\gamma_i$ for $1 \leq i \leq 2$ of the number
of points of $\gamma_1' \cap \gamma_2'$.  If $\gamma_1$ and $\gamma_2$ are either
oriented simple closed curves or elements of $\HH_1(\Sigma_g)$, then $\AlgI(\gamma_1,\gamma_2)$ will
denote the algebraic intersection number of $\gamma_1$ and $\gamma_2$.  Finally, we will say
that $x,y \in \pi_1(\Sigma_{g,n})$ are {\em completely distinct} if $x \neq y$ and $x \neq y^{-1}$.

For surfaces with boundary, the group $\Mod_{g,n}$ is defined to be the group of homotopy classes of orientation-preserving
homeomorphisms of $\Sigma_{g,n}$ that fix the boundary pointwise (the homotopies also must
fix the boundary).  Like in the closed surface case, the group $\Torelli_{g,1}$ is defined
to be the subgroup of $\Mod_{g,1}$ consisting of mapping classes that act trivially on $\HH_1(\Sigma_{g,1})$.  For surfaces
with more than 1 boundary component, there is more than one useful definition for the Torelli group (see
\cite{PutmanCutPaste} for a discussion).  We discuss one special definition in \S \ref{section:handledragprelim}.
As far group-theoretic conventions go, we define
$[g_1,g_2] = g_1^{-1} g_2^{-1} g_1 g_2$ and $g_1^{g_2} = g_2^{-1} g_1 g_2$.  Finally, we wish
to draw the reader's attention to the warning at the end of \S \ref{section:birmanexactsequence};
it is the source of several somewhat counterintuitive formulas.

\paragraph{Acknowledgements.}
I wish to thank my advisor Benson Farb for his enthusiasm and encouragement
and for commenting extensively on previous incarnations
of this paper.  I also wish to thank Joan Birman, Matt Day, Martin Kassabov, Justin Malestein, and
Ben Wieland for their comments on this project.  I particularly wish to thank an anonymous referee
for a very careful reading and many useful suggestions.
Finally, I wish to thank the Department
of Mathematics of the Georgia Institute of Technology for their hospitality
during the time in which parts of this paper were conceived.

\section{Preliminaries}
\label{section:preliminaries}

\subsection{The Birman exact sequence}
\label{section:birmanexactsequence}
In this section, we review the exact sequences of Birman and
Johnson \cite{BirmanExactSequence,BirmanBook,JohnsonFinite} that describe the effect on the mapping class group of gluing a
disc to a boundary component; these will be the basis for our inductive arguments.  We will need the following definition.

\begin{definition}
Consider a surface $\Sigma_{g,n}$.  Let $\ast \in \Sigma_{g,n}$ be a point.  We define $\Mod_{g,n}^{\ast}$, the 
{\em mapping class
group relative to $\ast$}, to be the group of orientation-preserving homeomorphisms of $\Sigma_{g,n}$ that fix $\ast$ and
the boundary pointwise modulo isotopies fixing $\ast$ and the boundary pointwise. 
\end{definition}

\noindent 
Let $b$ be a boundary component of $\Sigma_{g,n}$.  There is a natural embedding $\Sigma_{g,n} \hookrightarrow \Sigma_{g,n-1}$
induced by gluing a disc to $b$.  Let $\ast \in \Sigma_{g,n-1}$ be a point in the interior of the new disc.  Clearly we can
factor the induced map $\Mod_{g,n} \rightarrow \Mod_{g,n-1}$ into a composition
$$\Mod_{g,n} \longrightarrow \Mod_{g,n-1}^{\ast} \longrightarrow \Mod_{g,n-1}.$$
Now let $U\Sigma_{g,n-1}$ be the unit tangent bundle of $\Sigma_{g,n-1}$ and $\tilde{\ast}$
be any lift of $\ast$ to $U\Sigma_{g,n-1}$.  The combined work of Birman \cite{BirmanExactSequence} and 
Johnson \cite{JohnsonFinite} shows that (except for the degenerate cases where $(g,n)$ equals $(0,1)$, $(0,2)$, or $(1,1)$) 
all of our groups fit into the following commutative diagram with exact rows and columns.

\begin{center}
\begin{tabular}{c@{\hspace{0.05 in}}c@{\hspace{0.05 in}}c@{\hspace{0.05 in}}c@{\hspace{0.05 in}}c@{\hspace{0.05 in}}c@{\hspace{0.05 in}}c@{\hspace{0.05 in}}c@{\hspace{0.05 in}}c}
    &               & $1$                        &               & $1$                      &               &                        &               &     \\
    &               & $\downarrow$               &               & $\downarrow$             &               &                        &               &     \\
    &               & $\Z$                       & $=$           & $\Z$                     &               &                        &               &     \\
    &               & $\downarrow$               &               & $\downarrow$             &               &                        &               &     \\
$1$ & $\longrightarrow$ & $\pi_1(U\Sigma_{g,n-1},\tilde{\ast})$ & $\longrightarrow$ & $\Mod_{g,n}$     & $\longrightarrow$ & $\Mod_{g,n-1}$ & $\longrightarrow$ & $1$ \\
    &               & $\downarrow$               &               & $\downarrow$             &               & $\parallel$            &               &     \\
$1$ & $\longrightarrow$ & $\pi_1(\Sigma_{g,n-1},\ast)$  & $\longrightarrow$ & $\Mod_{g,n-1}^{\ast}$ & $\longrightarrow$ & $\Mod_{g,n-1}$ & $\longrightarrow$ & $1$ \\
    &               & $\downarrow$               &               & $\downarrow$             &               &                        &               &     \\
    &               & $1$                        &               & $1$                      &               &                        &               &
\end{tabular} 
\end{center} 

\Figure{figure:birmanexactsequence}{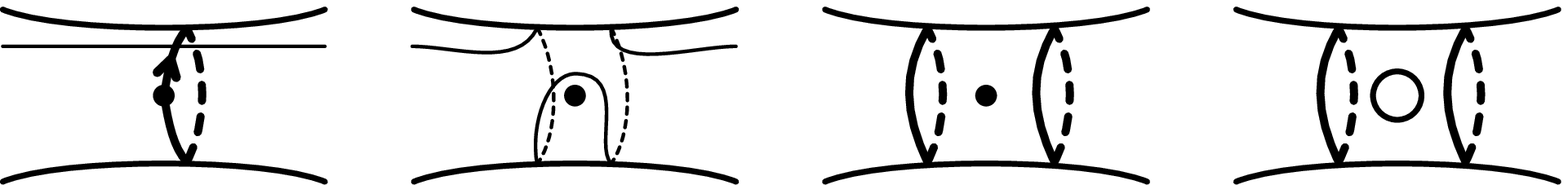}{a. A simple closed curve $\gamma \in \pi_1(\Sigma_{g,n-1})$
\CaptionSpace b. We drag $\ast$ around $\gamma$. \CaptionSpace c. $\PtPsh{\gamma} = T_{\gamma_1} T_{\gamma_2}^{-1}$
\CaptionSpace d. The lift $\LPtPsh{\gamma} = T_{\tilde{\gamma}_1} T_{\tilde{\gamma}_2}^{-1}$ of $\PtPsh{\gamma}$ to $\Mod_{g,n}$}

\noindent
The $\Z$ in the first column is the loop in the fiber, while the $\Z$ in the second column corresponds to the
Dehn twist about the
filled-in boundary component.  For $\gamma \in \pi_1(\Sigma_{g,n-1},\ast)$, 
let $\PtPsh{\gamma}$ be the element of $\Mod_{g,n-1}^{\ast}$
associated to $\gamma$ (hence $\PtPsh{\gamma}$ ``drags $\ast$ around the curve $\gamma$'').  If $\gamma$ 
is nontrivial and can be represented by a simple closed curve, then there is a nice
formula for $\PtPsh{\gamma}$ (see Figures
\ref{figure:birmanexactsequence}.a--c).  Namely, let $\gamma_1$ and $\gamma_2$ be the boundary of
a regular neighborhood of $\gamma$.  The orientation of $\gamma$ induces an orientation on $\gamma_1$ and $\gamma_2$;
assume that $\gamma$ lies to the left of $\gamma_1$ and to the right of $\gamma_2$.  Then $\PtPsh{\gamma}=T_{\gamma_1} T_{\gamma_2}^{-1}$.

Continue to assume that $\gamma \neq 1$ can be represented by a simple closed curve.  Recall that we have been
considering $\Sigma_{g,n-1}$ to be $\Sigma_{g,n}$ with a disc glued to $b$.  In the other direction, we can consider
$\Sigma_{g,n}$ to be $\Sigma_{g,n-1}$ with the point $\ast$ blown up to a boundary component (i.e.\ replaced
with its circle of unit tangent vectors).  Two such identifications of $\Sigma_{g,n}$ with a blow-up of $\Sigma_{g,n-1}$
may differ by a power of $T_b$; however, since $T_b$ fixes both $\gamma_1$ and $\gamma_2$ there are well-defined
lifts $\tilde{\gamma}_1$ and $\tilde{\gamma}_2$ of the $\gamma_i$ to $\Sigma_{g,n}$ (see Figure
\ref{figure:birmanexactsequence}.d).  It is not hard to see that 
$\LPtPsh{\gamma} := T_{\tilde{\gamma}_1} T_{\tilde{\gamma}_2}^{-1}$ is a lift of $\PtPsh{\gamma}$.

\begin{warning}
It is traditional to compose elements of $\pi_1$ from left to right (concatenation order) but
to compose mapping classes from right to left (functional order).  We will (reluctantly) adhere to
these conventions, but because of them the map $\pi_1(U\Sigma_{g,n-1}) \rightarrow \Mod_{g,n}$
and all other maps derived from it are {\em anti-homomorphisms}; i.e.\ they reverse the order
of composition.
\end{warning}

\subsection{Two group-theoretic lemmas}
\label{section:grouptheory}
In this section, we prove two easy group-theoretic lemmas that will
form the basis for many of our arguments.  The first is a tool for proving that sequences are exact.
\begin{lemma}
\label{lemma:exactnesstest}
Let $j : G_2 \rightarrow G_3$ be a surjective homomorphism between two groups
$G_2$ and $G_3$, and let $G_1$ be a normal subgroup of $G_2$ with $G_1 \subset \Ker(j)$.
Additionally, let $\langle S_3 | R_3 \rangle$ be a presentation for $G_3$ and $S_2$ be
a generating set for $G_2$ satisfying $j(S_2)=S_3$.  Assume that the following
two conditions are satisfied.
\begin{enumerate}
\item For all $s,s' \in S_2 \cup \{1\}$ with $j(s)=j(s')$, there exist $k_1,k_2 \in G_1$ so that $s = k_1 s' k_2$.
\item For any relation $r_1 \cdots r_k \in R_3$, we can find $\tilde{r}_1,\ldots,\tilde{r}_k \in S_2^{\pm 1}$ with
$\tilde{r}_1 \cdots \tilde{r}_k=1$ so that $j(\tilde{r}_i) = r_i$ for $1 \leq i \leq k$.
\end{enumerate}
Then the sequence
$$1 \longrightarrow G_1 \longrightarrow G_2 \stackrel{j}{\longrightarrow} G_3 \longrightarrow 1$$
is exact.
\end{lemma}
\begin{proof}
Let $\overline{S}_2 \subset G_2 / G_1$ be the projection of $S_2$.  By condition 1 the induced
map $\overline{j} : G_2 / G_1 \rightarrow G_3$ restricts to a bijection between $\overline{S}_2$
and $S_3$.  Condition 2 then implies that there is an inverse $\overline{j}^{-1}$; i.e.\ 
that $\overline{j}$ is an isomorphism, as desired. 
\end{proof}

\begin{remark}
In the first condition of Lemma \ref{lemma:exactnesstest}, since $G_1$ is normal it is enough to assume that there
exists some $k \in G_1$ so that $s = s' k$.  We stated it the way we did to make the logic behind some of our applications
clearer.
\end{remark}

\noindent
The following special case of Lemma \ref{lemma:exactnesstest} will be used repeatedly.

\begin{corollary}
\label{corollary:isomorphismtest}
Let $j : G_2 \rightarrow G_3$ be a surjective homomorphism between two groups.  Assume that
$G_3$ has a presentation $\langle S_3 | R_3 \rangle$ and that $G_2$ has a generating set
$S_2$ so that $j$ restricts to a bijection between $S_2$ and $S_3$.  Furthermore,
assume that every relation $r_1 \cdots r_k \in R_3$ (here $r_i \in S_3^{\pm 1}$) satisfies 
$j^{-1}(r_1) \cdots j^{-1}(r_k) = 1$, where $j^{-1}(r_i)$ is the unique
element of $S_2^{\pm 1}$ that is mapped to $r_i$.  Then $j$ is an isomorphism.
\end{corollary}

\noindent
Corollary \ref{corollary:isomorphismtest} is interesting even if $G_3$ is a free group -- 
it says that if $j : G \rightarrow F(S)$ is a homomorphism from a group $G$ to the free group $F$ 
on the free generating set $S$ and if for each $s \in S$ there is some $\tilde{s} \in j^{-1}(s)$
so that the set $\{\text{$\tilde{s}$ $|$ $s \in S$}\}$ generates $G$, then $j$ is an
isomorphism.

The second lemma is a tool for proving that a set of elements generates a group.

\begin{lemma}
\label{lemma:generationtest}
Let $G$ be a group generated by a set $S$.  Assume that a group $H$ generated
by a set $T$ acts on $G$ (as a set, not necessarily as a group) and that $S' \subset S$ satisfies the following
two conditions. 
\begin{enumerate}
\item $H(S') = S$
\item For $t \in T^{\pm 1}$ and $s \in S'$, we have $t(s) \in \Span{$S'$} \subset G$.
\end{enumerate}
Then $S'$ generates $G$.
\end{lemma}
\begin{proof}
By condition 2, the group $H$ stabilizes $\Span{$S'$} \subset G$.  Condition 1 then implies
that $S \subset \Span{$S'$}$, so $\Span{$S'$} = G$, as desired.
\end{proof}

\section{Non-formal relations in the Torelli group}
\label{section:relations}

In this section, we derive the non-formal relations in our presentation.

\begin{remark}
As will become clear, all the non-formal relations in our presentation arise in some
fashion from the Birman exact sequence.
\end{remark}

\subsection{The lantern and crossed lantern relations}
\label{section:subsurfacedrag}
\subsubsection{Preliminaries}
\label{section:subsurfacedragprelim}
We first discuss relations that arise from ``dragging subsurfaces around''.  Fix a simple closed 
separating curve $b$ on $\Sigma_g$.  Cutting $\Sigma_g$ along $b$, we
obtain subsurfaces homeomorphic to $\Sigma_{h_1,1}$ and $\Sigma_{h_2,1}$ for some integers $h_1$ and $h_2$ satisfying
$h_1 + h_2 = g$.  Assume that $h_1 > 0$ and $h_2 > 1$.  
Observe that we have an injection $\Torelli_{h_2,1} \hookrightarrow \Torelli_{g}$.  
Additionally, the formulas in \S \ref{section:birmanexactsequence} imply that the kernel 
$\pi_1(U\Sigma_{h_2})$ of the Birman exact sequence
for $\Sigma_{h_2,1}$ lies in $\Torelli_{h_2,1}$, so we have an exact sequence
\begin{equation}
1 \longrightarrow \pi_1(U\Sigma_{h_2}) \longrightarrow \Torelli_{h_2,1} \longrightarrow \Torelli_{h_2} \longrightarrow 1. \label{torseq1}
\end{equation}
Combining these two observations, we obtain an injection $\pi_1(U\Sigma_{h_2}) \hookrightarrow \Torelli_{g}$.  
The element of $\Mod_{g}$ that corresponds to $\gamma \in \pi_1(U\Sigma_{h_2})$ can be informally
described as ``dragging $\Sigma_{h_1,1}$ around $\gamma$''.  We will construct relations in $\pi_1(U\Sigma_{h_2})$
using the push-maps discussed in \S \ref{section:birmanexactsequence} and then use the aforementioned
injection to map these relations into $\Torelli_g$.

\subsubsection{The lantern relation}
\label{section:lantern}

Consider simple closed curves $x,y,z \in (\pi_1(\Sigma_{h_2}) \setminus \{1\})$ that can be arranged like the curves drawn
in Figure \ref{figure:surfacerelations}.a.  Observe that $xyz=1$ and that
$$\LPtPsh{x} = T_{\tilde{x}_1, \tilde{x}_2} \in \pi_1(U\Sigma_{h_2})$$
for the curves $\tilde{x}_1$ and $\tilde{x}_2$ depicted in Figure \ref{figure:surfacerelations}.c.  Similar statements
are true for $y$ and $z$.  We conclude that in $\pi_1(U\Sigma_{h_2}) \subset \Torelli_{g}$, we must have
$$T_{\tilde{z}_1,\tilde{z}_2} T_{\tilde{y}_1, \tilde{y}_2} T_{\tilde{x}_1,\tilde{x}_2} = T_b^k$$
for some $k$ (observe that we have switched the order of composition here from concatenation order for curves
to functional order for mapping classes).  By examining the action on a properly embedded arc exactly one
of whose endpoints lies on $b$, one can check that $k=1$.  These are
the classical {\em lantern relations} (see, e.g., \cite{JohnsonFirst}).  Summing up,
we have
\begin{equation*}
T_{\tilde{z}_1,\tilde{z}_2} T_{\tilde{y}_1,\tilde{y}_2} T_{\tilde{x}_1,\tilde{x}_2} = T_b \tag{L}\label{L}\\
\end{equation*}
for all curves $\tilde{x}_1$, $\tilde{x}_2$, $\tilde{y}_1$, $\tilde{y}_2$, $\tilde{z}_1$, and $\tilde{z}_2$ embedded in $\Sigma_g$ like the
curves in Figure \ref{figure:surfacerelations}.c.

\begin{remark}
This interpretation of the lantern relation was discovered independently by Margalit and McCammond \cite{MargalitBraids}.
\end{remark}

\subsubsection{The crossed lantern relation}
\label{section:crossedlantern}

Now consider simple closed curves $x,y,z \in (\pi_1(\Sigma_{h_2}) \setminus \{1\})$ that can be arranged like the curves
drawn in Figure \ref{figure:surfacerelations}.b.  Observe that $xy=z$ and that
$$\LPtPsh{x} = T_{\tilde{x}_1} T_{\tilde{x}_2}^{-1} \in \pi_1(U\Sigma_{h_2}) \subset \Torelli_{g}$$
for the curves $\tilde{x}_1$ and $\tilde{x}_2$ depicted in Figure \ref{figure:surfacerelations}.d.  Similar statements
are true for $y$ and $z$.  We conclude that in $\pi_1(U\Sigma_{h_2}) \subset \Torelli_{g}$, we must have
$$(T_{\tilde{y}_1} T_{\tilde{y}_2}^{-1})(T_{\tilde{x}_1} T_{\tilde{x}_2}^{-1}) = (T_{\tilde{z}_1} T_{\tilde{z}_2}^{-1}) T_b^k$$
for some $k$.  By examining the action on a properly embedded arc exactly one
of whose endpoints lies on $b$, one can check that $k=0$.  We will call these the {\em crossed lantern
relations}.  Summing up, our relation is
\begin{equation*}
T_{\tilde{y}_1,\tilde{y}_2} T_{\tilde{x}_1,\tilde{x}_2} = T_{\tilde{z}_1,\tilde{z}_2} \tag{CL}\label{CL}\\
\end{equation*}
for all curves $\tilde{x}_1$, $\tilde{x}_2$, $\tilde{y}_1$, $\tilde{y}_2$, $\tilde{z}_1$, and $\tilde{z}_2$
that can be embedded in $\Sigma_g$ like the curves depicted in Figure \ref{figure:surfacerelations}.d.

\begin{alternatederivation}
Observe that for $i=1,2$ we have $\tilde{z}_i = T_{\tilde{x}_2}(\tilde{y}_i)$.  Expanding out the
$T_{\tilde{z}_1,\tilde{z}_2}$ in \eqref{CL} as $T_{\tilde{x}_2} T_{\tilde{y}_1,\tilde{y}_2} T_{\tilde{x}_2}^{-1}$ and
rearranging terms, we see that \eqref{CL} is equivalent to $T_{\tilde{y}_1,\tilde{y}_2} T_{\tilde{x}_1} T_{\tilde{y}_1,\tilde{y}_2}^{-1} = T_{\tilde{x}_2}$.
This follows from the easily verified identity $T_{\tilde{y}_1,\tilde{y}_2}(\tilde{x}_1)=\tilde{x}_2$.
\end{alternatederivation}

\subsection{The Witt--Hall and commutator shuffle relations}
\label{section:handledrag}

\subsubsection{Preliminaries}
\label{section:handledragprelim}

We now examine the relations that arise from ``dragging the end of a handle''.  For use later in
\S \ref{section:exactsequencesweak},
we will discuss a slightly more general situation.  For $g \geq 0$ and $n \geq 2$,
let $i : \Sigma_{g,n} \hookrightarrow \Sigma_{g+n-1}$ be
the embedding of $\Sigma_{g,n}$ into the surface obtained by gluing the boundary components
of a copy of $\Sigma_{0,n}$ to the boundary
components of $\Sigma_{g,n}$.  Define
$\Torelli_{g,n} = i_{\ast}^{-1}(\Torelli_{g+n-1})$.
It is not hard to see that this is well-defined.  Observe that
$i_{\ast}(\Torelli_{g,2})$ is the subgroup of $\Torelli_{g+1}$ stabilizing
the handle corresponding to the glued-in annulus.  The groups $\Torelli_{g,n}$ were
introduced by Johnson \cite{JohnsonKg} and
investigated further by van den Berg \cite{VanDenBergThesis} and the author \cite{PutmanCutPaste}
(in the notation of \cite{PutmanCutPaste}, if the boundary components
of $\Sigma_{g,n}$ are $\{b_1,\ldots,b_n\}$, then
$\Torelli_{g,n} = \Torelli(\Sigma_{g,n},\{\{b_1,\ldots,b_n\}\})$).

We will say that a mapping class $f \in \Mod_{g,n}$ is a
separating twist, etc.,\ if $i_{\ast}(f)$ is a separating twist, etc.  It follows from
\cite[Theorem 1.3]{PutmanCutPaste} that if $g \geq 1$, then separating
twists and bounding pair maps generate $\Torelli_{g,n}$.

\begin{remark}
Not all simple closed curves that separate $\Sigma_{g,n}$ are nullhomologous.  By
our definition, separating twists in $\Mod_{g,n}$ are exactly
Dehn twists about nullhomologous simple closed curves.
\end{remark}

Let $b$ be a boundary component of $\Sigma_{g,n}$.  The kernel
$\pi_1(U\Sigma_{g,n-1})$ of the map $\Mod_{g,n} \rightarrow \Mod_{g,n-1}$
induced by gluing a disc to $b$ does not lie
in $\Torelli_{g,n}$ (for instance, $T_b \notin \Torelli_{g,n}$).
Instead, \cite[Theorem 4.1]{PutmanCutPaste} says that we have an exact sequence
\begin{equation}
1 \longrightarrow [\pi_1(\Sigma_{g,n-1}),\pi_1(\Sigma_{g,n-1})] \longrightarrow \Torelli_{g,n} \longrightarrow \Torelli_{g,n-1} \longrightarrow 1. \label{torseq2}
\end{equation}
The group $[\pi_1(\Sigma_{g,n-1}),\pi_1(\Sigma_{g,n-1})]$ is embedded
in $\pi_1(U\Sigma_{g,n-1}) \cong \pi_1(\Sigma_{g,n-1}) \otimes \Z$
as the graph of a homomorphism $\phi:[\pi_1(\Sigma_{g,n-1}),\pi_1(\Sigma_{g,n-1})] \rightarrow \Z$,
that is, as the set of all pairs $(x,\phi(x))$ for
$x \in [\pi_1(\Sigma_{g,n-1}),\pi_1(\Sigma_{g,n-1})]$.  The identification
$\pi_1(U\Sigma_{g,n-1}) \cong \pi_1(\Sigma_{g,n-1}) \otimes \Z$ (or, equivalently, the
splitting $\pi_1(\Sigma_{g,n-1}) \rightarrow \pi_1(U\Sigma_{g,n-1})$ of the natural
surjection $\pi_1(U\Sigma_{g,n-1}) \rightarrow \pi_1(\Sigma_{g,n-1})$) is not natural,
but once a splitting $\rho : \pi_1(\Sigma_{g,n-1}) \rightarrow \pi_1(U\Sigma_{g,n-1})$ is chosen
$\phi$ is uniquely defined by the requirement that the image of the homomorphism
$[\pi_1(\Sigma_{g,n-1}),\pi_1(\Sigma_{g,n-1})] \rightarrow \pi_1(U\Sigma_{g,n-1})$ defined by
$x \mapsto \rho(x) T_b^{\phi(x)}$ must be
contained in the pullback of $\Torelli_{g,n}$ under the inclusion
$\pi_1(U\Sigma_{g,n-1}) \hookrightarrow \Mod_{g,n}$.

\Figure{figure:torelliexactsequence}{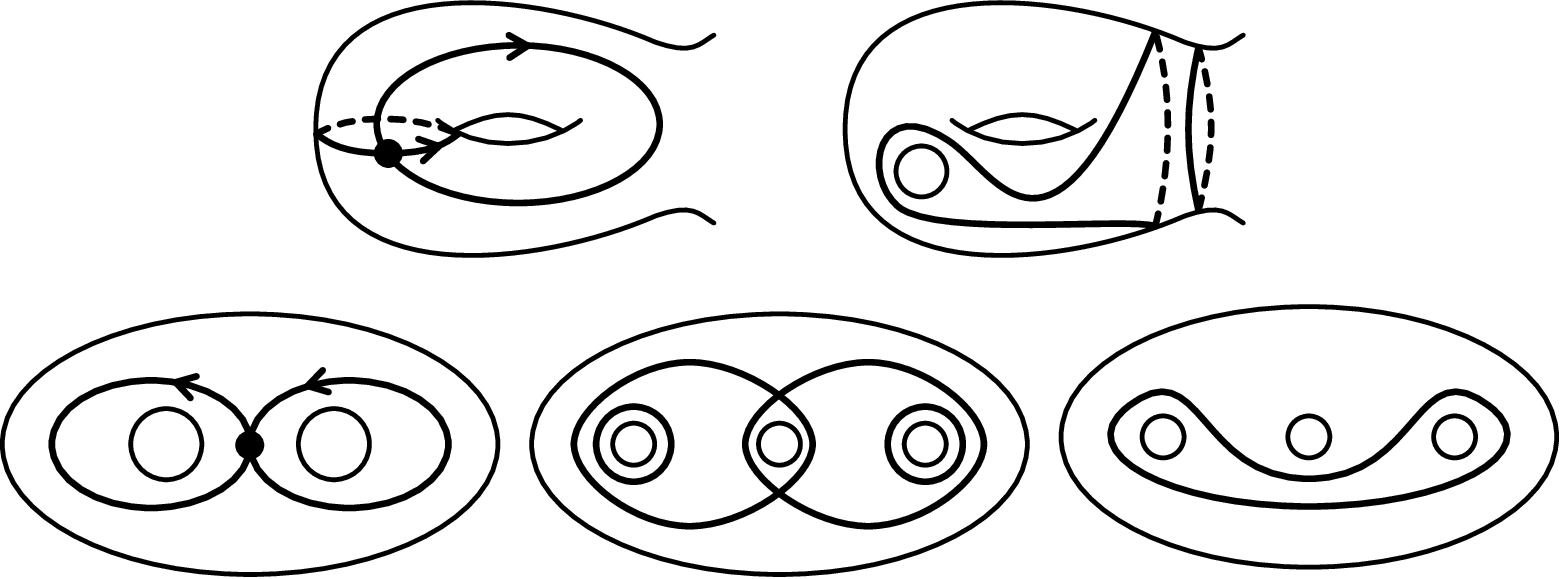}{a,c. The two
configurations of nontrivial simple closed
curves $\gamma^1,\gamma^2 \in \pi_1(\Sigma_{g,n-1},\ast)$ with $\gamma^1 \cap \gamma^2 = \{\ast\}$.
\CaptionSpace b,d. Lifts of the corresponding elements of
$[\pi_1(\Sigma_{g,n-1},\ast),\pi_1(\Sigma_{g,n-1},\ast)]$ to $\Torelli_{g,n}$
\CaptionSpace e. $\DComm{(\gamma^2)^{-1} (\gamma^1)^{-1}, (\gamma^2)^{-1}} = [T_{\tilde{\gamma}^2_1}, T_{\delta}^{-1}]$}

For curves $\gamma^1,\gamma^2 \in (\pi_1(\Sigma_{g,n-1}) \setminus \{1\})$,
define $\DComm{\gamma^1,\gamma^2}$ to be the element of $\Torelli_{g,n}$ associated to
$[\gamma^1,\gamma^2]$.  To simplify
our notation, if $\eta \in \pi_1(\Sigma_{g,n-1})$ is another simple closed curve, then we define
$$\DComm{\gamma^1,\gamma^2}^{\eta} := \DComm{(\eta^{-1})(\gamma^1)(\eta), (\eta^{-1})(\gamma^2)(\eta)} = \DComm{\PtPsh{\eta}(\gamma^1),\PtPsh{\eta}(\gamma^2)}.$$
Finally,
if $\gamma \in (\pi_1(\Sigma_{g,n-1}) \setminus \{1\})$ is already an element of
the commutator subgroup, then let $\DComm{\gamma}$
be the element of $\Torelli_{g,n}$ associated to $\gamma$.

We will need some explicit formulas for $\DComm{\cdot,\cdot}$.  Consider two completely distinct simple
closed curves $\gamma^1,\gamma^2 \in (\pi_1(\Sigma_{g,n-1}) \setminus \{1\})$
that only intersect at the base point.  From the above description,
we see that the following procedure will yield $\DComm{\gamma^1,\gamma^2}$.
\begin{enumerate}
\item Choose some $\psi \in \Mod_{g,n}$ which is associated to an element of
$\pi_1(U\Sigma_{g,n-1})$ that projects
to $[\gamma^1,\gamma^2] \in \pi_1(\Sigma_{g,n-1})$.
\item Determine $k \in \Z$ so that $\psi T_b^k \in \Torelli_{g,n}$.  We will
then have $\DComm{\gamma^1,\gamma^2} = \psi T_b^k$.
\end{enumerate}
There are two cases.  In the first (see Figure \ref{figure:torelliexactsequence}.a), a
regular neighborhood of $\gamma^1 \cup \gamma^2$ is homeomorphic to $\Sigma_{1,1}$.  Observe
that $[\gamma^1,\gamma^2]$ is homotopic
to a simple closed separating curve.  Our element $\psi$ in this case will be
$\LPtPsh{[\gamma^1,\gamma^2]}$.  Observe that
$$\LPtPsh{[\gamma^1,\gamma^2]} = T_{\widetilde{[\gamma^1,\gamma^2]}_1} T_{\widetilde{[\gamma^1,\gamma^2]}_2}^{-1}$$
for simple closed curves $\widetilde{[\gamma^1,\gamma^2]}_1$ and $\widetilde{[\gamma^1,\gamma^2]}_2$ like
the curves pictured in Figure \ref{figure:torelliexactsequence}.b.

Note that exactly one element of the pair
$\{\widetilde{[\gamma^1, \gamma^2]}_1,\widetilde{[\gamma^1, \gamma^2]}_2\}$ is a
separating curve (both curves separate $\Sigma_{g,n}$, but only one of
them maps to a separating curve on $\Sigma_{g+n-1}$).  In Figure
\ref{figure:torelliexactsequence}.b, the curve $\widetilde{[\gamma^1, \gamma^2]}_2$
is separating, but in other situations
$\widetilde{[\gamma^1,\gamma^2]}_1$ will be the separating curve (for
instance, this will happen if we flip the labels
on the curves $\gamma^1$ and $\gamma^2$ in Figure \ref{figure:torelliexactsequence}.a).
Now, the nonseparating curve and $b$ form a bounding pair on $\Sigma_{g,n}$.
We conclude that either
$$\DComm{\gamma^1,\gamma^2} = T_{\widetilde{[\gamma^1, \gamma^2]}_1, b} T_{\widetilde{[\gamma^1, \gamma^2]}_2}^{-1}$$
or
$$\DComm{\gamma^1,\gamma^2} = T_{\widetilde{[\gamma^1, \gamma^2]}_1} T_{b, \widetilde{[\gamma^1, \gamma^2]}_2},$$
depending on which curve is separating.  In a similar way, if $\gamma$ is a separating
curve then $\DComm{\gamma}$ equals
the product of a separating twist and a bounding pair map.

In the second case, a regular neighborhood of $\gamma^1 \cup \gamma^2$ is homeomorphic to $\Sigma_{0,3}$
(see Figure \ref{figure:torelliexactsequence}.c).  In this case, our element $\psi$ will be
$$\LPtPsh{\gamma^2} \LPtPsh{\gamma^1} \LPtPsh{\gamma^2}^{-1} \LPtPsh{\gamma^1}^{-1}.$$
Lifting everything to $\Sigma_{g,n}$, we see that
$$\LPtPsh{\gamma^i} = T_{\tilde{\gamma}^i_1}^{e_i} T_{\tilde{\gamma}^i_2}^{-e_i}$$
for the curves depicted in Figure \ref{figure:torelliexactsequence}.d and some $e_i=\pm 1$
(the $e_i$ depend on the orientations of $\gamma^1$ and $\gamma^2$).  Observe that
$$\LPtPsh{\gamma^2} \LPtPsh{\gamma^1} \LPtPsh{\gamma^2}^{-1} \LPtPsh{\gamma^1}^{-1} = T_{\tilde{\gamma}^2_1}^{e_2} T_{\tilde{\gamma}^1_1}^{e_1} T_{\tilde{\gamma}^2_1}^{-e_2} T_{\tilde{\gamma}^1_1}^{-e_1} = [T_{\tilde{\gamma}^2_1}^{-e_i},T_{\tilde{\gamma}^1_1}^{-e_i}] \in \Torelli_{g,n}.$$
We conclude that
$\DComm{\gamma^1,\gamma^2} = [T_{\tilde{\gamma}^2_1}^{-e_2},T_{\tilde{\gamma}^1_1}^{-e_2}]$.  Now, this
is the commutator of the simply intersecting pair $\{\tilde{\gamma}^2_1, \tilde{\gamma}^1_1\}$ if
$e_2=e_1=-1$; we will
call a pair of curves $\gamma^1$ and $\gamma^2$ with this property {\em positively aligned}.
If $\gamma^1$ and $\gamma^2$ are not positively aligned, however, then by
repeatedly applying the commutator identity $[g_1^{-1},g_2] = [g_2,g_1]^{g_1^{-1}}$
and the fact that $T_{x} T_{y} T_{x}^{-1} = T_{T_{x}(y)}$ for simple closed curves $x$ and $y$, we can
find a simply intersecting pair $C_{\rho^1,\rho^2}$ with
$\DComm{\gamma^1,\gamma^2} = C_{\rho^1,\rho^2}$.  We conclude that
$\DComm{\gamma^1,\gamma^2}$ is a commutator of {\em some} simply
intersecting pair no matter how $\gamma^1$ and $\gamma^2$ are
aligned.

\begin{example}
\label{example:nonunique}
We can now give an example of the non-uniqueness of the expression of a mapping
class as a commutator of a simply intersecting pair.  Orienting $\gamma^1$ and
$\gamma^2$ as shown in Figure \ref{figure:torelliexactsequence}.a, we have
$\DComm{(\gamma^1)^{-1},(\gamma^2)^{-1}} = C_{\tilde{\gamma}^2_1,\tilde{\gamma}^1_1}$ (verifying
this is a good exercise in understanding the above construction).  Let $\delta$
be the curve in Figure \ref{figure:torelliexactsequence}.e.  Observe that
$\DComm{(\gamma^2)^{-1} (\gamma^1)^{-1}, (\gamma^2)^{-1}} = [T_{\tilde{\gamma}^2_1}, T_{\delta}^{-1}]$.  
Since we have the commutator identity
$[(\gamma^2)^{-1} (\gamma^1)^{-1}, (\gamma^2)^{-1}] = [(\gamma^1)^{-1}, (\gamma^2)^{-1}]$, we conclude that
$C_{\tilde{\gamma}^2_1,\tilde{\gamma}^1_1} = [T_{\tilde{\gamma}^2_1}, T_{\delta}^{-1}]$.  The right
hand side of this is not a commutator of a simply intersecting pair, but the above procedure
shows that it equals $C_{\delta,T_{\delta}^{-1}(\tilde{\gamma}^2_1)}$.
\end{example}

We conclude with the following lemma.

\begin{lemma}
\label{lemma:recognizesip}
Let $s \in \Torelli_{g,n}$ be a commutator of a simply intersecting pair whose image
under the map $\Torelli_{g,n} \rightarrow \Torelli_{g,n-1}$ is $1$.  Then there are completely distinct simple
closed curves $\gamma^1,\gamma^2 \in \pi_1(\Sigma_{g,n-1})$ that only intersect at the basepoint
so that $s = \DComm{\gamma^1,\gamma^2}$.
\end{lemma}
\begin{proof}
Let $s = C_{x,y}$.  Then a regular neighborhood $N$ of $x \cup y$ satisfies $N \cong \Sigma_{0,4}$.  Moreover,
our assumptions imply that some boundary component of $N$ must be isotopic to $b$.  The lemma
then follows from Figures \ref{figure:torelliexactsequence}.c--d and the above discussion.
\end{proof}

\subsubsection{Witt--Hall relations}
\label{section:witthall}
\Figure{figure:witthall}{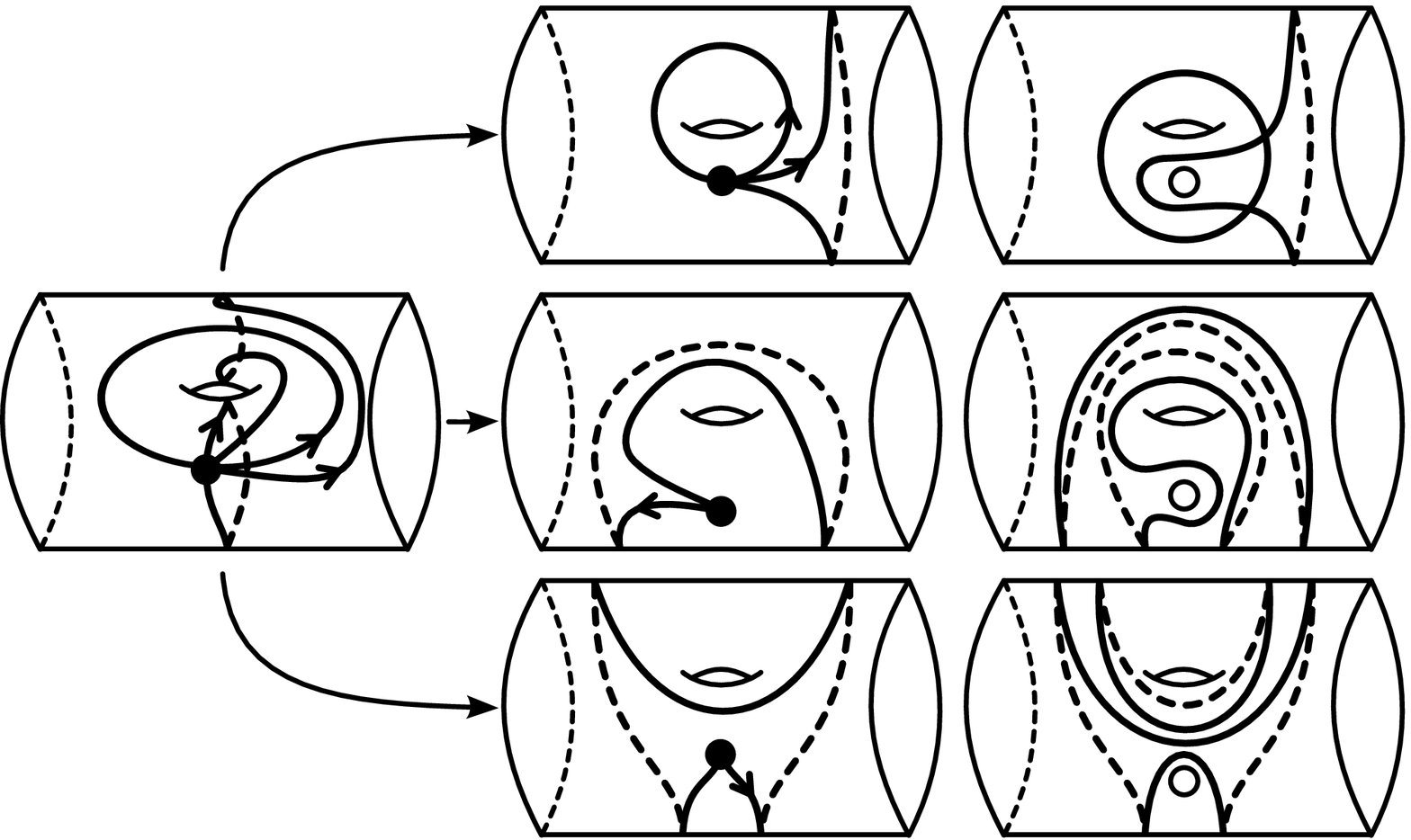}{a. One configuration of curves yielding a Witt--Hall relation \CaptionSpace b,c. The
curves needed for Example \ref{example:witthall}}

In this section and in \S \ref{section:commutatorshuffle}, we will derive relations in
the group $\Torelli_{g,2}$.  These relations give us relations in the Torelli groups of closed surfaces
in the following way.  For $g' \geq g$, let $\Sigma_{g,2} \hookrightarrow \Sigma_{g'}$ 
be {\em any} embedding (not just the embedding $\Sigma_{g,2} \hookrightarrow \Sigma_{g+1}$
discussed in \S \ref{section:handledragprelim}).  There is then an induced map
$\Torelli_{g,2} \rightarrow \Torelli_{g'}$ (``extend by the identity''; see \cite[Theorem Summary 1.1]{PutmanCutPaste}).  This
induced map takes separating twists, bounding pair maps, and simply intersecting pair maps to generators
of the same type (possibly degenerate ones, such as bounding pair maps $T_{x,y}$ with $x$ isotopic
to $y$).  If 
$$s_1^{e_1} \cdots s_k^{e_k}=1 \quad (e_i = \pm 1)$$
is a relation between separating twists, bounding pair maps, and simply intersecting pair maps
in $\Torelli_{g,2}$ and $s_i'$ is the image of $s_i$ in $\Torelli_{g'}$ via the above
map, then we obtain a relation between our generators in $\Torelli_{g'}$ by 
deleting all the degenerate generators in the relation $(s_1')^{e_1} \cdots (s_k')^{e_k}=1$.

The two families of relations that we derive from exact sequence \eqref{torseq2} come from 
commutator identities.  First, consider 
the Witt--Hall commutator identity
$$[g_1 g_2, g_3] = [g_1,g_3]^{g_2} [g_2,g_3].$$

\begin{remark}
The Witt--Hall commutator identity first appeared in \cite{HallPGroups}.  Later, it appeared in a
list of basic commutator identities dubbed the ``Witt--Hall identities'' in \cite{MagnusKarassSolitar}.
\end{remark}

Fix $x,y,z \in (\pi_1(\Sigma_{g,1}) \setminus \{1\})$ so that for each of the sets $\{x,y,z\}$ and
$\{xy,z\}$, the elements of the set can be represented by completely distinct simple closed curves
that only intersect at the basepoint.  There are several
different topological types of configurations of curves with these
properties; an example is in Figure \ref{figure:witthall}.a.
The Witt--Hall commutator identity then yields the following relation, which we will call the {\em Witt--Hall relation}.
\begin{equation*}
\DComm{xy,z} = \DComm{y,z} \DComm{x, z}^y. \label{WH}\tag{WH}\\
\end{equation*}
We now give an example.

\begin{example}
\label{example:witthall}
The curves $x$, $y$, and $z$ depicted in Figure \ref{figure:witthall}.a satisfy the conditions for the
Witt--Hall relations.  In the surface group, the relation is $[xy,z]=[x,z]^y [y,z]$.  In Figure
\ref{figure:witthall}.b, we depict the curves involved in this surface group relation.  Let $z_1$, $(xy)_2$, $c_1$, $c_2$,
$a_1$, and $a_2$ be the curves depicted in Figure \ref{figure:witthall}.c.  We then have 
$\LPtPsh{[x,z]^y} = T_{a_1,a_2}$ and $\LPtPsh{[y,z]}=T_{c_1,c_2}$.  The corresponding relation in Torelli
is $[T_{z_1}^{-1},T_{(xy)_2)}] = T_{c_1,b} T_{c_2}^{-1} T_{a_1} T_{b,a_2}$ (the counterintuitive form of the initial commutator
comes from the fact that the map from the kernel of the Birman exact sequence to Torelli is an anti-homomorphism).  However,
$[T_{z_1}^{-1},T_{(xy)_2)}]$ is {\em not} a commutator of a simply intersecting pair (i.e.\ $z$ and $xy$ are not positively
aligned).  Using the relation $[g_1^{-1},g_2]=[g_1 g_2 g_1^{-1},g_1]$,
we transform this into the Witt--Hall relation $C_{T_{z_1}((xy)_2),z_1} = T_{c_1,b} T_{c_2}^{-1} T_{a_1} T_{b,a_2}$.
\end{example}

\subsubsection{Commutator shuffle relations}
\label{section:commutatorshuffle}
We now use another, somewhat less standard commutator identity to find
relations in the Torelli group.  Our commutator identity, which is easily verified, is the following.
$$[g_1,g_2]^{g_3} = [g_3, g_1] [g_3, g_2]^{g_1} [g_1,g_2] [g_1,g_3]^{g_2} [g_2,g_3].$$
Though it may seem a bit odd, it will become apparent in \S \ref{section:commutatorsubgroup} that this
is exactly the relation we need to complete our picture.  We will apply it to completely distinct simple
closed curves $x,y,z \in (\pi_1(\Sigma_{g,1}) \setminus \{1\})$ that only intersect at the basepoint. 
Again, there are finitely many topological types of such configurations.
Our relation is then
\begin{equation*}
\DComm{x, y}^z = \DComm{y,z} \DComm{x,z}^y \DComm{x,y} \DComm{z,y}^x \DComm{z,x}. \label{CS}\tag{CS}\\
\end{equation*}
We will call these relations the {\em commutator shuffles}.  Pictures of them 
are left as an exercise for the reader.

\section{The Main Theorem}
\label{section:maintheorem}

\subsection{A stronger version of the Main Theorem}
\label{section:precisestatement}
To facilitate our induction, we will have to consider not only the
case of a closed surface but also the case of a surface with boundary.  In
this section, we state a version of our Main Theorem that applies
to these cases.  We begin with a definition.

\begin{definition}
For $g \geq 2$ and $n \geq 0$, define $\Gamma_{g,n}$ to be the group whose
generating set is the set of all separating twists, all bounding pair maps, and
all commutators of simply intersecting pairs on $\Sigma_{g,n}$ and whose relations
are the following.  For $n = 0$, they are relations \eqref{F.1}-\eqref{F.8} from \S \ref{section:introduction},
relations \eqref{L} and \eqref{CL} from \S \ref{section:subsurfacedrag}, and
relations \eqref{WH} and \eqref{CS} from \S \ref{section:handledrag} (for the relations
\eqref{WH} and \eqref{CS}, we use all ways of ``embedding them in the closed surface''
as described in the beginning of \S \ref{section:witthall}).
For $n=1$, they are the set of all words $r$ in the generators of $\Gamma_{g,n}$ so that $i_{\ast}(r)$ is one of the above relations, 
where $i : \Sigma_{g,1} \hookrightarrow \Sigma_{g+1}$ is the embedding obtained by
gluing a copy of $\Sigma_{1,1}$ to $\Sigma_{g,1}$ and $i_{\ast}$ is the obvious
map defined on the generators.  For $n>1$, they are the set of all words
$r$ in the generators of $\Gamma_{g,n}$ so that $i_{\ast}(r)$ is one of the above relations, where
$i : \Sigma_{g,n} \hookrightarrow \Sigma_{g+n-1,1}$ obtained by gluing $n$ boundary
components of a copy of $\Sigma_{0,n+1}$ to the boundary components of
$\Sigma_{g,n}$ and $i_{\ast}$ is the obvious
map defined on the generators.
\end{definition}

\begin{remark}
The generators for $\Gamma_{g,n}$ are {\em mapping classes}, not merely abstract symbols.  For
bounding pair maps and separating twists, this is unimportant, as their defining curves are
determined by their mapping classes.  For commutators of simply intersecting pairs, however, different
pairs of curves determine the same mapping class (see Example \ref{example:nonunique}), 
and we identify these in $\Gamma_{g,n}$.
\end{remark}

\noindent
Since all of the relations of $\Gamma_{g,n}$ also hold in $\Torelli_{g,n}$, there is a natural
homomorphism $\Gamma_{g,n} \rightarrow \Torelli_{g,n}$.  A stronger version of Theorem \ref{theorem:maintheoreminformal} is then the following.
\begin{theorem}[{Main Theorem, Stronger Version}]
\label{theorem:maintheoremprecise}
For $n \leq 2$ and $g \geq 2$, the natural map $\Gamma_{g,n} \rightarrow \Torelli_{g,n}$ is
an isomorphism.
\end{theorem}

\begin{remark}
In fact, this is also true for $n > 2$, but Theorem \ref{theorem:maintheoremprecise} is all we need.  We will use the groups
$\Gamma_{g,n}$ for $n > 2$ later for technical purposes.
\end{remark}

\subsection{Obtaining presentations from group actions}
In this section, we discuss a theorem of the author \cite{PutmanPresentation}
that we will use to prove Theorem \ref{theorem:maintheoremprecise}.  In order to state it, we begin by noting that
an argument of Armstrong \cite{ArmstrongGenerators}
says that if $X$ is a simply connected
simplicial complex and a group $G$ acts without rotations on $X$ (that
is, for all simplices $s$ of $G$ the stabilizer $G_s$ stabilizes $s$ pointwise;
this can be arranged by subdividing $X$), then if $X/G$ is also simply connected
we can conclude that $G$ is generated by elements
that stabilize vertices.  In other words, we have a surjective
map
$$\pi : \BigFreeProd_{v \in X^{(0)}} G_v \longrightarrow G.$$
As notation, for $v \in X^{(0)}$ denote the inclusion map
$$G_v \hookrightarrow \BigFreeProd_{v \in X^{(0)}} G_v$$
by $i_v$.  

There are then some obvious
elements $\Ker(\pi)$, which we write as relations $f=g$ rather than as elements
$f g^{-1}$.  First, we have $i_v(g) i_w(h) i_v(g^{-1}) = i_{g \cdot w}(ghg^{-1})$
for $g \in G_v$ and $h \in G_w$.  We call these relations the {\em conjugation
relations}.  Second, we have $i_v(g) = i_{v'}(g)$ if $g \in G_v \cap G_{v'}$ and
$\{v,v'\} \in X^{(1)}$ (here $\{v,v'\} \in X^{(1)}$ means that $\{v,v'\}$ forms an
edge in the 1-skeleton of $X$).  We call these the {\em edge relations}.
The following theorem of the author says that under favorable circumstances these two families of relations
yield the entire kernel of the aforementioned map.
\begin{theorem}[{\cite{PutmanPresentation}}]
\label{theorem:presentation}
Let a group $G$ act without rotations on a simply connected simplicial
complex $X$.  Assume that $X/G$ is $2$-connected.  Then
$$G = (\BigFreeProd_{v \in X^{(0)}} G_v)/R,$$
where $R$ is the normal subgroup generated by the conjugation
relations and the edge relations.
\end{theorem}

\subsection{The proof of the Main Theorem}
\label{section:proof}
In this section, we will give the outline of the proof of Theorem \ref{theorem:maintheoremprecise}.  Our main 
tool will be Theorem \ref{theorem:presentation} 
together with two other results whose proofs are postponed until later sections.

The first major ingredient in our proof will be the following proposition, which is proven in
\S \ref{section:curvestorelliconnected}.  Recall that the complex $\ModCurves_g$
was defined at the end of \S \ref{section:introduction}.

\begin{proposition}
\label{proposition:curvestorelliconnected}
The simplicial complex $\ModCurves_g$ satisfies the following two properties.
\begin{enumerate}
\item The complex $\ModCurves_g$ is $(g-2)$-connected.
\item The complex $\ModCurves_g / \Torelli_{g}$ is $(g-1)$-connected.
\end{enumerate}
\end{proposition}

\begin{remark}
In fact, using similar methods 
one can prove that $\Curves_g / \Torelli_{g}$ is $(g-1)$-connected, but Proposition 
\ref{proposition:curvestorelliconnected} suffices for our purposes, and the details of its proof are less technical.
In the end, one would get the same presentation no matter which of the two complexes one used.
\end{remark}

Theorem \ref{theorem:presentation} and Proposition \ref{proposition:curvestorelliconnected} will allow
us to give an inductive decomposition of $\Torelli_{g,n}$.  To show that the groups $\Gamma_{g,n}$
fit into this inductive picture, we will show that the groups $\Gamma_{g,n}$ fit into exact sequences
like exact sequence \eqref{torseq1} from \S \ref{section:subsurfacedragprelim} and exact sequence
\eqref{torseq2} from \S \ref{section:handledragprelim}.  More precisely, observe that there exist natural 
``disc-filling'' homomorphisms $\Gamma_{g,1} \rightarrow \Gamma_{g}$ and $\Gamma_{g,2} \rightarrow \Gamma_{g,1}$ (defined
on the generators).  In \S \ref{section:gammaexactsequences}, we will prove the following.  

\begin{proposition}
\label{proposition:exactsequences}
The aforementioned homomorphisms fit into the following exact sequences.
\begin{align}
1 \longrightarrow \pi_1(U\Sigma_{g}) \longrightarrow \Gamma_{g,1} \longrightarrow \Gamma_{g} \longrightarrow 1, \label{gammaseq2}\\
1 \longrightarrow [\pi_1(\Sigma_{g,1}),\pi_1(\Sigma_{g,1})] \longrightarrow \Gamma_{g,2} \longrightarrow \Gamma_{g,1} \longrightarrow 1. \label{gammaseq1}
\end{align}
\end{proposition}

\noindent
Next, we will need the following lemma, which forms part of Lemma \ref{lemma:commutatorexpress} below.

\begin{lemma}
\label{lemma:commutatorexpresssimple}
For $g \geq 2$ and $0 \leq n \leq 2$, using the relations in $\Gamma_{g,n}$ we can write
any commutator of a simply intersecting pair as a product of bounding pairs maps and
separating twists.
\end{lemma}

\noindent
Finally, we will need some results of Mess and Johnson about separating twists.  Recall
that by convention, all homology groups have $\Z$-coefficients.  Observe
that if $\gamma$ is a separating curve on $\Sigma_{g}$ that cuts $\Sigma_{g}$ into
two subsurfaces $S_1$ and $S_2$, then we have an splitting
$$\HH_1(\Sigma_{g}) \cong \HH_1(S_1) \oplus \HH_1(S_2),$$
where the $\HH_1(S_i)$ are symplectic $\Z$-modules which are orthogonal
with respect to the intersection form.  We will call such a splitting a
{\em symplectic splitting}.  Observe that the symplectic splitting associated
to $\gamma$ is a conjugacy invariant of $T_{\gamma} \in \Torelli_{g}$.  We
then have the following two theorems.

\begin{theorem}[{Mess, \cite{MessTorelli}}]
\label{theorem:mess}
$\Torelli_{2}$ is an infinitely generated free group.  Moreover, there exists
a free generating set of separating twists $S$ containing exactly one separating twist
associated to each symplectic splitting of $\HH_1(\Sigma_{2})$.
\end{theorem}

\begin{theorem}[{Johnson, \cite{JohnsonConjugacy}}]
\label{theorem:johnsonconjugacy}
For $g \geq 2$, two separating twists $T_{\gamma_1}$ and $T_{\gamma_2}$ in $\Torelli_{g}$
are conjugate if and only if they induce the same symplectic splitting of
$\HH_1(\Sigma_{g})$.
\end{theorem}

We now assemble these ingredients to prove Theorem \ref{theorem:maintheoremprecise}.

\begin{proof}[{Proof of Theorem \ref{theorem:maintheoremprecise}}]
The proof will be by induction on $g$ and $n$.  We begin with the base case $(g,n)=(2,0)$.

\BeginClaims
\begin{claim}
The natural map $\Gamma_{2} \rightarrow \Torelli_{2}$ is an isomorphism.
\end{claim}
\BeginClaimProof
Observe first that $\Sigma_{2}$ does not contain any bounding pairs.  Also, using Lemma
\ref{lemma:commutatorexpresssimple} we see that $\Gamma_{2}$ is generated by
separating twists.  Let $S$ be the generating set for $\Torelli_{2}$ given
by Theorem \ref{theorem:mess}.  Using the conjugation relation \eqref{F.6} together with Theorem
\ref{theorem:johnsonconjugacy}, we conclude that $\Gamma_{g,2}$ is generated by $\{\text{$T_{\gamma}$ $|$ $\gamma \in S$}\}$.  Corollary 
\ref{corollary:isomorphismtest} therefore implies that the natural map $\Gamma_{2} \rightarrow \Torelli_{2}$
is an isomorphism, as desired.
\EndClaimProof

\noindent
Now assume by induction that for some $g \geq 2$ the natural map $\Gamma_{g} \rightarrow \Torelli_{g}$ is
an isomorphism.

\begin{claim}
The natural maps $\Gamma_{g,1} \rightarrow \Torelli_{g,1}$ and $\Gamma_{g,2} \rightarrow \Torelli_{g,2}$ are isomorphisms.
\end{claim}
\BeginClaimProof
Using Proposition \ref{proposition:exactsequences}, we have the following commutative
diagram of exact sequences.
\begin{center}
\begin{tabular}{c@{\hspace{0.05 in}}c@{\hspace{0.05 in}}c@{\hspace{0.05 in}}c@{\hspace{0.05 in}}c@{\hspace{0.05 in}}c@{\hspace{0.05 in}}c@{\hspace{0.05 in}}c@{\hspace{0.05 in}}c}
$1$ & $\rightarrow$ & $\pi_1(U\Sigma_{g})$         & $\rightarrow$ & $\Gamma_{g,1}$                         & $\rightarrow$ & $\Gamma_{g}$               & $\rightarrow$ & $1$ \\
    &               & $\parallel$  &               & $\downarrow$                     &               & $\downarrow$            &               &     \\
$1$ & $\rightarrow$ & $\pi_1(U\Sigma_{g})$         & $\rightarrow$ & $\Torelli_{g,1}$ & $\rightarrow$ & $\Torelli_{g}$ & $\rightarrow$ & $1$\\
\end{tabular}
\end{center}

\noindent
The right hand map is an isomorphism by induction, so the five lemma 
implies that the center map is an isomorphism; i.e.\ that $\Gamma_{g,1} \cong \Torelli_{g,1}$.  The proof that
$\Gamma_{g,2} \cong \Torelli_{g,2}$ is similar.
\EndClaimProof

\noindent
We now prove the following.

\begin{claim}
The natural map $\Gamma_{g+1} \rightarrow \Torelli_{g+1}$ is an isomorphism.
\end{claim}
\BeginClaimProof
Since no two curves in a simplex of $\ModCurves_{g+1}$ are homologous, the group $\Torelli_{g+1}$
acts on $\ModCurves_{g+1}$ without rotations.  Since $g+1 \geq 3$, 
Proposition \ref{proposition:curvestorelliconnected} and Theorem \ref{theorem:presentation} thus imply that
$$\Torelli_{g+1} \cong (\BigFreeProd_{\text{$\gamma \in (\ModCurves_{g+1})^{(0)}$}} (\Torelli_{g+1})_{\gamma}) / R,$$
where $(\Torelli_{g+1})_{\gamma}$ denotes the stabilizer in $\Torelli_{g+1}$ of $\gamma$ and 
where $R$ is the normal subgroup generated by the edge relations and the conjugation relations coming from 
the action of $\Torelli_{g+1}$ on $\ModCurves_{g+1}$.  Now,
consider a simple closed nonseparating curve $\gamma$, and let $b$ and $b'$ be the boundary components
of the copy of $\Sigma_{g,2}$ that results from cutting $\Sigma_{g+1}$ along $\gamma$.  
By \cite[Theorem 4.1]{ParisRolfsonMod}, we have an exact sequence
$$1 \longrightarrow \langle T_{b,b'} \rangle \longrightarrow \Torelli_{g,2} \longrightarrow (\Torelli_{g+1})_{\gamma} \longrightarrow 1.$$
If we denote by $(\Gamma_{g+1})_{\gamma}$ the subgroup of $\Gamma_{g+1}$ generated by the subset of generators
that do not intersect $\gamma$, then there is a surjective homomorphism $\Gamma_{g,2} \rightarrow (\Gamma_{g+1})_{\gamma}$.
Letting $K$ denote the kernel of this surjection, we have
a commutative diagram of exact sequences
\begin{center}
\begin{tabular}{c@{\hspace{0.05 in}}c@{\hspace{0.05 in}}c@{\hspace{0.05 in}}c@{\hspace{0.05 in}}c@{\hspace{0.05 in}}c@{\hspace{0.05 in}}c@{\hspace{0.05 in}}c@{\hspace{0.05 in}}c}
$1$ & $\rightarrow$ & $K$ & $\rightarrow$ & $\Gamma_{g,2}$           & $\rightarrow$ & $(\Gamma_{g+1})_\gamma$& $\rightarrow$ & $1$ \\
    &               & $\downarrow$                         &               & $\downarrow$             &               & $\downarrow$             &               &     \\
$1$ & $\rightarrow$ & $\langle T_{b,b'}\rangle$  & $\rightarrow$ & $\Torelli_{g,2}$ & $\rightarrow$ & $(\Torelli_{g+1})_\gamma$ & $\rightarrow$ & $1$\\
\end{tabular}
\end{center}
By induction, the center map is an isomorphism.  Also, we have $T_{b,b'} \in K$, so the left hand vertical map
is surjective.  By the five lemma, we conclude that the map 
$(\Gamma_{g+1})_\gamma \rightarrow (\Torelli_{g+1})_{\gamma}$ is an isomorphism.  

Now, every generator of $\Gamma_{g+1}$ lies in $(\Gamma_{g+1})_{\gamma}$ for some simple closed nonseparating
curve $\gamma$.  Hence there is a surjection
$$\BigFreeProd_{\text{$\gamma \in (\ModCurves_{g+1})^{(0)}$}} (\Gamma_{g+1})_{\gamma} \longrightarrow \Gamma_{g+1}.$$
Since the map $(\Gamma_{g+1})_\gamma \rightarrow (\Torelli_{g+1})_{\gamma}$ is an isomorphism, we conclude
that there is a surjective map
$$\BigFreeProd_{\text{$\gamma \in (\ModCurves_{g+1})^{(0)}$}} (\Torelli_{g+1})_{\gamma} \longrightarrow \Gamma_{g+1}.$$
The edge relations in $R$ project to trivial relations in $\Gamma_{g+1}$.  Also, using relations \eqref{F.6}--\eqref{F.8}, we
see that the conjugation relations in $R$ project to relations in $\Gamma_{g+1}$.
We conclude that we have a sequence of surjections
$$(\BigFreeProd_{\text{$\gamma \in (\ModCurves_{g+1})^{(0)}$}} (\Torelli_{g+1})_{\gamma}) / R \longrightarrow \Gamma_{g+1} \longrightarrow \Torelli_{g+1}.$$
Since the composition of these two maps is an isomorphism, we conclude that the natural map $\Gamma_{g+1} \rightarrow \Torelli_{g+1}$ is an isomorphism, as desired.
\EndClaimProof

\noindent
This completes the proof of Theorem \ref{theorem:maintheoremprecise}, which we recall
is stronger than Theorem \ref{theorem:maintheoreminformal} from the introduction.
\end{proof}

It remains to prove Propositions \ref{proposition:curvestorelliconnected} and 
\ref{proposition:exactsequences} and Lemma \ref{lemma:commutatorexpresssimple}.  The
proofs of Proposition \ref{proposition:curvestorelliconnected} and Lemma \ref{lemma:commutatorexpresssimple}
are contained in \S \ref{section:gammaexactsequences}, while the proof of Proposition 
\ref{proposition:curvestorelliconnected} is contained in \S \ref{section:curvestorelliconnected}.

\section{Exact sequences for $\Gamma_{g,n}$ : The proof of Proposition \ref{proposition:exactsequences}}
\label{section:gammaexactsequences}
The proof of Proposition \ref{proposition:exactsequences} will be split into two pieces.  Before discussing
these two pieces, recall the following.

\begin{itemize}
\item In \S \ref{section:birmanexactsequence}, we defined a bounding pair map $\LPtPsh{\gamma} \in \Torelli_{g,1}$
for every nontrivial $\gamma \in \pi_1(\Sigma_g)$ that can be realized by a simple closed curve.  
Together with the twist about the boundary component,
these bounding pair maps generate the kernel of the Birman exact sequence
$$1 \longrightarrow \pi_1(U\Sigma_{g}) \longrightarrow \Torelli_{g,1} \longrightarrow \Torelli_{g} \longrightarrow 1;$$
the key observation is that $\pi_1(\Sigma_g)$ is generated by simple closed curves.
\item Consider $n \geq 2$.  In \S \ref{section:handledragprelim}, we defined an element $\DComm{x,y} \in \Torelli_{g,n}$
for every pair $x,y \in \pi_1(\Sigma_{g,n-1})$ of completely distinct nontrivial elements that can be realized by simple closed
curves that only intersect at the basepoint.  Additionally, we showed that $\DComm{x,y}$ is either
a commutator of a simply intersecting pair or a well-defined product of a bounding pair map and a separating
twist.  The group $[\pi_1(\Sigma_{g,n-1}),\pi_1(\Sigma_{g,n-1})]$ is generated by the set of all $[x,y]$ where 
$x,y \in \pi_1(\Sigma_{g,n-1})$ range over pairs satisfying the above conditions.  Thus the elements
$\DComm{x,y}$ generate the kernel of the Birman exact sequence
$$1 \longrightarrow [\pi_1(\Sigma_{g,n-1}),\pi_1(\Sigma_{g,n-1})] \longrightarrow \Torelli_{g,n} \longrightarrow
\Torelli_{g,n-1} \longrightarrow 1.$$
\end{itemize}

For most of this section, we will only consider $\Torelli_{g,n}$ for $n \leq 2$; the cases where
$n>2$ will play a small role in \S \ref{section:consequences1}.  
In the following definition, we will abuse notation and identify $\LPtPsh{\gamma}$ and $\DComm{x,y}$ with
the corresponding products of generators in $\Gamma_{g,1}$ and $\Gamma_{g,2}$.

\begin{definition}
For $g \geq 2$, let $K_{g,1}$ be the subgroup of $\Gamma_{g,1}$ generated by the set $S^K_{g,1}$ that is
defined as follows (here $b$ is the boundary component of $\Sigma_{g,1}$).
$$S^K_{g,1}:=\{T_b\} \cup \{\LPtPsh{\gamma} \text{ $|$ $\gamma \in (\pi_1(\Sigma_{g}) \setminus \{1\})$ can be realized by a simple closed curve}\}.$$
Also, let $K_{g,2}$ be the subgroup of $\Gamma_{g,2}$ generated by the set $S^K_{g,2}$ that
is defined as follows.
\begin{align*}
S^K_{g,2}:=\{\DComm{x,y} \text{ $|$ } &\text{$x,y \in (\pi_1(\Sigma_{g,1},\ast) \setminus \{1\})$ are completely distinct and can be realized by} \\
                           &\text{simple closed curves that only intersect at the basepoint} \}.
\end{align*}
\end{definition}

\begin{remark}
The set $S^K_{g,1}$ is contained in the generating set for $\Gamma_{g,1}$, but the set
$S^K_{g,2}$ is not contained in the generating set for $\Gamma_{g,2}$.
\end{remark}

\noindent
The first part of our proof of Proposition \ref{proposition:exactsequences} is the following lemma, 
which will be proven in \S \ref{section:exactsequencesweak}.

\begin{lemma}
\label{lemma:exactsequencesweak}
For $g \geq 2$ and $1 \leq n \leq 2$ we have an exact sequence
$$1 \longrightarrow K_{g,n} \longrightarrow \Gamma_{g,n} \longrightarrow \Gamma_{g,n-1} \longrightarrow 1.$$
\end{lemma}

\noindent
The second part of our proof is the following lemma, which will be proven in \S \ref{section:kernelidentification}.

\begin{lemma}
\label{lemma:kernelidentification}
For $g \geq 2$, the natural maps $K_{g,1} \rightarrow \pi_1(U\Sigma_{g})$ and 
$K_{g,2} \rightarrow [\pi_1(\Sigma_{g,1}),\pi_1(\Sigma_{g,1})]$ are isomorphisms.
\end{lemma}

Proposition \ref{proposition:exactsequences} is an immediate consequence of
Lemmas \ref{lemma:exactsequencesweak} and \ref{lemma:kernelidentification}.

\subsection{Constructing the exact sequences : Lemma \ref{lemma:exactsequencesweak}}
\label{section:exactsequencesweak}
The goal of this section is to prove Lemma \ref{lemma:exactsequencesweak}.  There
are three parts.

\begin{itemize}
\item In \S \ref{section:discfilling}, we investigate the effect of the map
$\Gamma_{g,n} \rightarrow \Gamma_{g,n-1}$ on the generators of $\Gamma_{g,n}$.
\item In \S \ref{section:consequences1} -- \S \ref{section:consequences2}, we work
out several consequences of the relations in $\Gamma_{g,n}$.
\item In \S \ref{section:proofexactsequencesweak}, we give the proof of Lemma \ref{lemma:exactsequencesweak}.
\end{itemize}

\subsubsection{The effect on generators of filling in boundary components}
\label{section:discfilling}
In this section, fix $g \geq 0$ and $1 \leq n \leq 2$.  Also, fix a boundary component $b$ of $\Sigma_{g,n}$, and let
$i:\Sigma_{g,n} \hookrightarrow \Sigma_{g,n-1}$ be the embedding induced by
gluing a disc to $b$.  This induces a map $i_{\ast} : \Mod_{g,n} \rightarrow \Mod_{g,n-1}$
(``extend by the identity'').  We begin with the following definition.

\begin{definition}
Let $x$ and $x'$ be two nontrivial simple closed
curves on $\Sigma_{g,n}$.  We say that $x$ and $x'$ {\em differ by $b$} if there is an embedding
$\Sigma_{0,3} \hookrightarrow \Sigma_{g,n}$ that takes the boundary components of $\Sigma_{0,3}$ to
$x$, $x'$, and $b$.
\end{definition}

\noindent
The following lemma is immediate.

\begin{lemma}
\label{lemma:differbybsep}
If $x$ and $x'$ are nontrivial simple closed curves that differ by $b$, then $i_{\ast}(T_{x})=i_{\ast}(T_{x'})$,
and additionally there is a simple closed curve $\gamma \in \pi_1(\Sigma_{g,n-1})$ with $T_{x,x'} = \LPtPsh{\gamma}$.
\end{lemma}

\noindent
Also, the following lemma follows from the discussion in \S \ref{section:handledragprelim} (see
especially Figure \ref{figure:torelliexactsequence}).

\begin{lemma}
\label{lemma:differbybsepbd}
Assume that $n=2$ and that $x$ and $x'$ are simple closed curves that differ by $b$.  Also, assume that
$T_{x}$ is a separating twist.  Thus $x'$ separates the two boundary components, so $T_{b,x'}$ is a bounding pair map.  
Then there is some
$\gamma \in \pi_1(\Sigma_{g,1})$ that can be realized by a simple closed separating curve so that
$T_{x} T_{b,x'} = \DComm{\gamma}$.
\end{lemma}

\noindent
Lemma \ref{lemma:differbybsep} shows that $i_{\ast}(s) = i_{\ast}(s')$ if the generators 
$s$ and $s'$ differ by the following moves.

\begin{definition}
Let $s$ and $s'$ be either separating twists, bounding pair maps, or
commutators of simply intersecting pairs.  We say that {\em $s$ differs from $s'$ by $b$} if
they satisfy one of the following conditions.
\begin{itemize}
\item $s = T_{x}$ and $s' = T_{x'}$ for separating curves $x$ and $x'$ that differ by $b$.  This
can only occur if $n=1$.
\item Either $s = T_{x,y}$ and $s' = T_{x',y}$ or $s = T_{y,x}$ and $s' = T_{y,x'}$ for
bounding pairs $\{x,y\}$ and $\{x',y\}$ so that $x$ differs from $x'$ by $b$.  This can 
only occur if $n=1$.
\item Either $s = T_{x,b}$ and $s' = T_{x'}$ or $s = T_{x'}$ and $s' = T_{x,b}$
for a bounding pair $\{x,b\}$ and a separating curve $x'$ so that either $x = x'$ or $x$ differs from $x'$ by $b$.
This can occur if $n=1$ or $n=2$; if $n=1$, then $T_x$ is also a separating twist.
\item Either $s = C_{x,y}$ and $s' = C_{x',y}$ or $s = C_{y,x}$ and $s' = C_{y,x'}$ for
simply intersecting pairs $\{x,y\}$ and $\{x',y\}$ so that $x$ differs from $x'$ by $b$.  This
can occur if $n=1$ or $n=2$.
\end{itemize}
Also, we say that $s$ and $s'$ {\em differ by a $b$-push map} if there exists some 
$\phi \in \pi_1(U\Sigma_{g,n-1}) = \Ker(i_{\ast})$
so that $s$ and $s'$ satisfy one of the following conditions.
\begin{itemize}
\item For a separating curve $x$ we have $s = T_x$ and $s' = T_{\phi(x)}$.
\item For a bounding pair $\{x,y\}$ we have $s = T_{x,y}$ and $s' = T_{\phi(x),\phi(y)}$.
\item For a simply intersecting pair $\{x,y\}$ we have $s = C_{x,y}$ and $s' = C_{\phi(x),\phi(y)}$.
\end{itemize}
We say that $s$ and $s'$ are {\em $b$-equivalent} if there is a sequence $s_1,\ldots,s_k$
of separating twists, bounding pair maps, or commutators of simply
intersecting pairs so that $s = s_1$, so that $s' = s_k$, and so that for $1 \leq j < k$ either
$s_j$ differs from $s_{j+1}$ by $b$ or $s_j$ and $s_{j+1}$ differ by a $b$-push map.
\end{definition}

We now prove the following.
\begin{lemma}
\label{lemma:generatoridentification1}
Let $s,s' \in \Torelli_{g,n}$ be separating twists, 
bounding pair maps, or commutators of simply intersecting pairs that satisfy $i_{\ast}(s) = i_{\ast} (s') \neq 1$.
Then $s$ and $s'$ are $b$-equivalent.
\end{lemma}
\begin{proof}
Assume first that $s$ and $s'$ are separating twists $T_x$ and $T_{x'}$.  Observe that the curve
$i_{\ast}(x)$ is isotopic to the curve $i_{\ast}(x')$ (here we are using the fact that if $\gamma_1$
and $\gamma_2$ are separating curves, then $T_{\gamma_1} = T_{\gamma_2}$ if and only if $\gamma_1$
is isotopic to $\gamma_2$).  Let $\phi_t : \Sigma_{g,n-1} \rightarrow \Sigma_{g,n-1}$
be an isotopy so that $\phi_0 = 1$
and $\phi_1 (i_{\ast}(x)) = i_{\ast}(x')$.  Restricting $\phi_t$ to the disc glued to $b$, we get a family of embeddings
of a disc into $\Sigma_{g,n-1}$.  If $i_{\ast}(x')$ does not separate $b$ from $\phi_1 (b)$, then we can modify $\phi_t$ so that
$\phi_1 (i_{\ast}(x)) = i_{\ast}(x')$ and $\phi_1 (b) = b$.  In this case, $\phi_t$ determines a mapping class
$\phi \in \pi_1(U\Sigma_{g,n-1}) \subset \Mod_{g,n}$ with $\phi(x)=x'$, and we are done.
If instead $i_{\ast}(x')$ separates $b$ from $\phi_1 (b)$, then we can
modify $\phi_t$ so that $\phi_1(b) = b$ but (letting
$\phi \in \pi_1(U\Sigma_{g,n-1}) \subset \Mod_{g,n}$ be the mapping class induced by $\phi_t$) so that
$\phi(x)$ and $x'$ differ by $b$ (we ``pull $b$ through $x'$'').
The desired sequence of generators is then $T_{x}, T_{\phi(x)}, T_{x'}$.

The proof is similar if $s$ and $s'$ are both bounding pair maps or both commutators of simply intersecting pairs.
Only two addenda are necessary.
\begin{itemize}
\item In both cases we may need to ``pull $b$'' through both of the curves that define $s'$.
\item While bounding pair maps are determined by their defining curves, simply intersecting
pair maps are not.  However, in the definition
of differing by $b$ and differing by a $b$-push map we only required that there be 
{\em some} simply intersecting pairs $\{x,y\}$
and $\{x',y'\}$ satisfying the conditions so that $s = C_{x,y}$ and $s' = C_{x',y'}$.  To make
the above argument work, we need to choose these pairs so that they become isotopic after gluing a
disc to $b$.
\end{itemize}

It remains to consider the case that (reordering $s$ and $s'$ if necessary) $s$ is a
bounding pair map and $s'$ is a separating twist -- it is not hard to see that the other
possibilities (for instance, that $s$ is a separating twist while $s'$ is a simply intersecting
pair map) are impossible.  In this case, we must have
$s = T_{x,b}$ (we cannot have $s = T_{b,x}$ since $s'$ is a positive twist).  An argument
similar to the argument in the previous two paragraphs then shows that $s$ and $s'$ are $b$-equivalent.
\end{proof}

\subsubsection{Consequences of our relations : commutators of simply intersecting pairs}
\label{section:consequences1}

Fix a surface $\Sigma_{g,n}$ with $g \geq 1$, with $n \geq 0$, and with $(g,n) \neq (1,1)$.  If $n \geq 1$, then 
let $b \subset \partial \Sigma_{g,n}$ be a boundary component and let
$i:\Sigma_{g,n} \hookrightarrow \Sigma_{g,n-1}$ and $i_{\ast} : \Mod_{g,n} \rightarrow \Mod_{g,n-1}$
be the maps induced by gluing a disc to $b$.  The main result of this section is the following.

\begin{lemma}
\label{lemma:commutatorexpress}
Assume that $n \leq 2$.  Let $s$ be a commutator of a simply intersecting pair on $\Sigma_{g,n}$.
\begin{enumerate}
\item Using the relations in $\Gamma_{g,n}$, we can write $s = s_1 \cdots s_k$ for some $k$, where
the $s_j$ are separating twists or bounding pair maps.
\item If $1 \leq n \leq 2$ and if $t$ is another commutator of a simply intersecting pair that differs from $s$
by $b$, then using the relations in $\Gamma_{g,n}$, we can write $s = s_1 \cdots s_k$ and $t = t_1 \cdots t_k$
for some $k$, where $s_j$ and $t_j$ are separating twists or bounding pair maps with $i_{\ast}(s_j) = i_{\ast}(t_j)$
for $1 \leq j \leq k$.
\item If $1 \leq n \leq 2$ and $i_{\ast}(s)=1$, then using the relations in $\Gamma_{g,n}$, we can write
$s = s_1 \cdots s_k$ for some $k$, where $s_j \in S^K_{g,n}$ for $1 \leq j \leq k$.
\end{enumerate}
\end{lemma}

For the proof of Lemma \ref{lemma:commutatorexpress}, we will need a lemma.  For $n \geq 2$, define
\begin{align*}
T^K_{g,n} = \{\text{$\DComm{x}$ $|$ }&\text{$x \in (\pi_1(\Sigma_{g,n-1},\ast) \setminus \{1\})$ can be
realized by a simple closed curve that}\\
&\text{cuts off a subsurface homeomorphic to $\Sigma_{1,1}$}\}.
\end{align*}
Our lemma is as follows.

\begin{lemma}
\label{lemma:commutatorexpresspower}
Consider $n \geq 2$.  Let $s$ be a commutator of a simply intersecting pair on $\Sigma_{g,n}$.
Assume that $i_{\ast}(s)=1$.  Then by using the relations in $\Gamma_{g,n}$, we can write
$s = s_1 \cdots s_k$ for some $k$, where $s_j \in T^K_{g,n}$ for $1 \leq j \leq k$.
\end{lemma}

In fact, Lemma \ref{lemma:commutatorexpresspower} follows immediately from a known result about commutator
subgroups of surface groups.  If $(\Sigma,\ast)$ is a compact surface with a basepoint $\ast \in \Interior(\Sigma)$ and
$x,y,z \in (\pi_1(\Sigma,\ast) \setminus \{1\})$ are  
such that for each of the sets $\{x,y,z\}$ and $\{xy,z\}$, all
the curves in the set can be realized by completely distinct nontrivial simple closed curves that only intersect
at the basepoint, then we will call the relation
$$[xy,z] = [x,z]^y [y,z]$$
a {\em Witt--Hall} relation in $[\pi_1(\Sigma,\ast),\pi_1(\Sigma,\ast)]$.  Observe that
since $x$ and $z$ can be realized by simple closed curves that only intersect at the basepoint, so
can $x^y = \PtPsh{y}(x)$ and $z^y = \PtPsh{y}(z)$.  Of course, the Witt--Hall relation in the Torelli
group is modeled on this commutator relation.  It is obvious that Lemma \ref{lemma:commutatorexpresspower}
follows from the following lemma combined with Lemma \ref{lemma:recognizesip}.

\begin{lemma}[{\cite[Lemma A.1]{PutmanCutPaste}}]
\label{lemma:commutatorgenerators}
Let $(\Sigma,\ast)$ be a compact surface of positive genus with a basepoint $\ast \in \Interior(\Sigma)$.
Let $\gamma^1,\gamma^2 \in \pi_1(\Sigma,\ast)$ be completely distinct simple closed curves
that only intersect at the basepoint.  Then by using a sequence of Witt--Hall relations
in $[\pi_1(\Sigma,\ast),\pi_1(\Sigma,\ast)]$, we can write
$$[\gamma^1,\gamma^2] = [\eta^1_1,\eta^2_1] \cdots [\eta^1_k,\eta^2_k],$$
where for $1 \leq j \leq k$ the curves $\eta^1_j$ and $\eta^1_j$ are completely distinct
simple closed curves so that $[\eta^1_j, \eta^2_j]$ can be realized by a simple closed
separating curve that cuts off a subsurface homeomorphic to $\Sigma_{1,1}$.
\end{lemma}

\begin{remark}
This result as stated is more precise than \cite[Lemma A.1]{PutmanCutPaste}; the proof
there actually proves the indicated result.
\end{remark}

\begin{proof}[{Proof of Lemma \ref{lemma:commutatorexpress}}]
We begin with conclusion 3.  The case $n=2$ follows from Lemma
\ref{lemma:recognizesip}, so we only need to consider the case $n=1$ (the reason
the case $n=1$ is harder is that $K_{g,1}$ does not contain any commutators
of simply intersecting pairs).
Let $i' : \Sigma_{g,2} \hookrightarrow \Sigma_{g,1}$ be an embedding so that if
the boundary components of $\Sigma_{g,2}$ are $b'$ and $b''$, then $i'(b') = b$ and $i'(b'')$ is
a simple closed curve that bounds a disc.  By \cite[Theorem Summary 1.1]{PutmanCutPaste}, there is an induced map 
$i'_{\ast} : \Torelli_{g,2} \rightarrow \Torelli_{g,1}$.  Let $\pi : \Torelli_{g,2} \rightarrow \Torelli_{g,1}$
be the map induced by gluing a disc to $b'$ (this is {\em different} from the map $i'_{\ast}$).  There
is then a simply intersecting pair $s' \in S^K_{g,2}$ so that $i'_{\ast}(s') = s$ and $\pi(s')=1$.
Lemma \ref{lemma:commutatorexpresspower} shows
that using the relations in $\Gamma_{g,2}$, we can write
$s' = \DComm{z^1} \cdots \DComm{z^k}$, where for $1 \leq j \leq k$ the element 
$z^j \in \pi_1(\Sigma_{g,1})$ can be represented by a nontrivial simple
closed separating curve.  Hence $s = i'_{\ast}(\DComm{z^1}) \cdots i'_{\ast}(\DComm{z^k})$
is a consequence of the Witt--Hall relations.  Now, for $1 \leq j \leq k$ the mapping class $i'_{\ast}(\DComm{z^j})$
is equal (up to taking inverses) to $T_{\rho_j,b} T_{\rho_j'}$, where 
$\rho_j$ and $\rho_j'$ are separating curves that differ by $b$.
This is not a generator for $K_{g,1}$, but we can use relation \eqref{F.4} twice 
together with \eqref{F.6} (which says that $T_b$ commutes with $T_{\rho_j'}$) to rewrite it as $T_{\rho_j,\rho_j'} T_{b}^{-1}$, 
which is a product of two generators for $K_{g,1}$ by Lemma \ref{lemma:differbybsep}.  This completes the proof of 
conclusion 3.  

To prove conclusion 1, we first show that for some $m \geq 2$ there exists an embedding 
$i'' : \Sigma_{1,m} \hookrightarrow \Sigma_{g,n}$ with an associated homomorphism
$i''_{\ast} : \Torelli_{1,m} \rightarrow \Torelli_{g,n}$ so that the following holds.  For
some simply intersecting pair map $s'' \in \Torelli_{1,m}$ that gets mapped to $1$ when a disc is glued to one
of the boundary components of $\Sigma_{1,m}$, we have $s = i''_{\ast}(s'')$.  Indeed, let $N \cong \Sigma_{0,4}$
be a regular neighborhood of the curves defining $s$.  We then simply choose a genus $1$ subsurface
containing $N$ and sharing a boundary component with $N$.

Now, in Lemma \ref{lemma:commutatorexpresspower} we proved that we can use the relations in $\Gamma_{1,m}$ to write $s''$ as a 
product of elements of $T^K_{1,m}$.  Since every element of $T^K_{1,m}$ is the product of a
separating twist and a bounding pair map, we obtain an expression $s'' = y_1 \cdots y_l$, where
$y_j$ is a separating twist or bounding pair map on $\Sigma_{1,m}$ for $1 \leq j \leq l$.  We conclude that
the relations in $\Gamma_{g,n}$ yield the desired expression $s = i''_{\ast}(y_1) \cdots i''_{\ast}(y_l)$.  

For conclusion 2, observe that there must exist an embedding $i''' : \Sigma_{1,m} \hookrightarrow \Sigma_{g,n}$
with an associated homomorphism $i'''_{\ast} : \Torelli_{1,m} \rightarrow \Torelli_{g,n}$ so
that $i'''_{\ast}(s'') = t$ and so that the embeddings $i \circ i'': \Sigma_{1,m} \hookrightarrow \Sigma_{g,n-1}$
and $i \circ i''' : \Sigma_{1,m} \hookrightarrow \Sigma_{g,n-1}$ are isotopic.  The desired expression
for $t$ is then $t = i'''_{\ast}(y_1) \cdots i'''_{\ast}(y_l)$.
\end{proof}

\subsubsection{Consequences of our relations : generators differing by a $b$-push map}
\label{section:consequences2}

\Figure{figure:conjugateunderb}{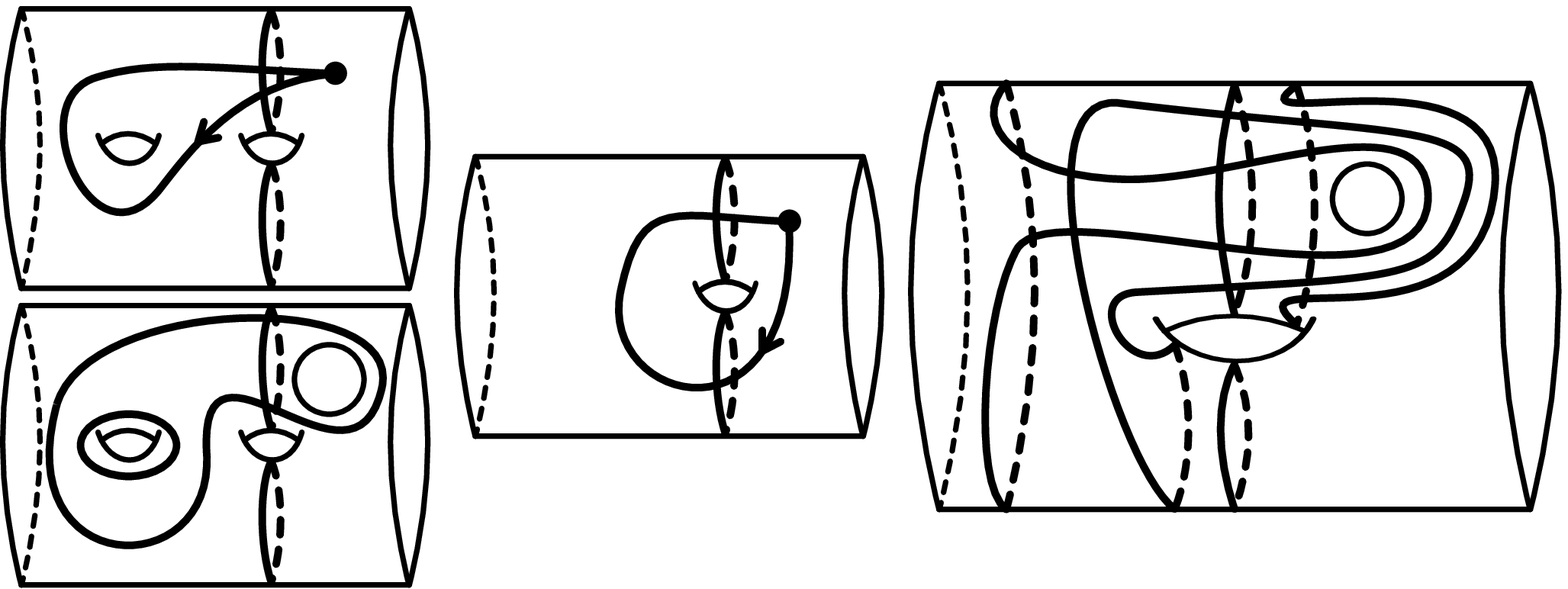}{The various configurations of curves needed for the proof of Lemma \ref{lemma:conjugateunderb}}

In this section, we prove the following.

\begin{lemma}
\label{lemma:conjugateunderb}
Fix $g \geq 2$ and $1 \leq n \leq 2$, and let $s$ and $s'$ be either separating twists,
bounding pair maps, or commutators of simply intersecting pairs.  If $s$ and $s'$ differ by 
a $b$-push map, then in $\Gamma_{g,n}$ the element $s$ is equal to
$k_1 s' k_2$ with $k_1,k_2 \in K_{g,n}$.
\end{lemma}
\begin{proof}
We begin by observing that for $n=1$, this is an immediate consequence of
the conjugation relations \eqref{F.6}--\eqref{F.8} (the point being that $\pi_1(U\Sigma_g) \subset \Torelli_{g,1}$
and $K_{g,1}$ surjects onto $\pi_1(U\Sigma_g)$).  We can therefore assume that $n=2$.  

Next, we claim that it is enough to prove the lemma for bounding pair maps $s$ and $s'$ so that
$s$ (and hence $s'$) does not equal $T_{x,y}$ with $T_{x}$ (and hence $T_y$) a separating twist.
Indeed, assume that
the lemma is true for such bounding pair maps and that $s = T_z$ and
$s' = T_{\psi(z)}$ for a separating curve $z$ and some $\psi \in \pi_1(U\Sigma_{g,1}) \subset \Mod_{g,2}$.
We can then find a simple closed curve $z'$ that differs
from $z$ by $b$.  By Lemma \ref{lemma:differbybsepbd}, we have 
$T_{z} T_{b,z'} \in K_{g,n}$ and $T_{\psi(z)} T_{b,\psi(z')} \in K_{g,n}$.  Now, neither
$T_b$ nor $T_{z'}$ is a separating twist, so by assumption there 
exists $k_1',k_2' \in K_{g,n}$ so that $T_{b,\psi(z')} = k_1' T_{z',b} k_2'$.  We
conclude that
\begin{align*}
T_{\psi(z)} &= (T_{\psi(z)} T_{b,\psi(z')}) T_{b,\psi(z')}^{-1} 
= (T_{\psi(z)} T_{b,\psi(z')}) (k_2')^{-1} T_{b,z'}^{-1} (k_1')^{-1}\\
&= (T_{\psi(z)} T_{b,\psi(z')}) (k_2')^{-1} (T_{z} T_{b,z'})^{-1} T_{z} (k_1')^{-1},
\end{align*}
so we can take $k_1 = (T_{\psi(z)} T_{b,\psi(z')}) (k_2')^{-1} (T_{z} T_{b,z'})^{-1}$ and $k_2 = (k_1')^{-1}$.

If instead $s$ is a commutator of a simply intersecting pair, then we can use 
Lemma \ref{lemma:commutatorexpress} to write $s = s_1^{\pm 1} \cdots s_k^{\pm 1}$,
where the $s_i$ are separating twists or bounding pair maps.  Since $K_{g,n}$ is normal, this 
reduces us to the previous cases.  Finally, if $s = T_{x,y}$ with $T_{x}$ (and hence $T_{y}$) a
separating twist, then we can use relation \eqref{F.4} to reduce ourselves to the case
of separating twists.

We can therefore assume that both $s$ and $s'$ are bounding pair maps of the above form.  
We claim that we can assume furthermore that either 
$s = T_{z,b}$ or $s = T_{x,y}$ with neither $x$ nor $y$ separating the surface (we remark that 
since $n = 2$, separating the surface is strictly weaker than being the curve in a separating twist).
Indeed, assume that $s = T_{x,y}$, where both $x$ and $y$ separate the surface (it is impossible for only one of them
to separate the surface) but where $T_x$ (and hence $T_y$) is not a separating twist.  Both $\{x,b\}$ and $\{y,b\}$ are bounding
pairs, and hence we can use relation \eqref{F.3} to write $s = T_{x,b} T_{b,y}$, reducing ourselves to the indicated
situation.

We will do the case that $s = T_{x,y}$ with neither $x$ nor $y$ separating the surface; the other
case is similar.  We must show that for all $\phi \in \pi_1(U\Sigma_{g,1}) \subset \Mod_{g,2}$, there
exists some $k_1,k_2 \in K_{g,2}$ so that $T_{\phi(x),\phi(y)} = k_1 T_{x,y} k_2$.  It is enough check 
this for all $\phi$ in a generating set for $\pi_1(U\Sigma_{g,1})$.  Draw $x$ and $y$ like
the curves in Figure \ref{figure:conjugateunderb}.a (we will systematically confuse the surface $\Sigma_{g,2}$ with the
surface $\Sigma_{g,1}$ that results from gluing a disc to $b$).  
Our generating set $S_{U\Sigma}$ for $\pi_1(U\Sigma_{g,1})$ will consist of $T_b$ plus the set of all $\LPtPsh{\gamma}$
for based simple closed curves $\gamma$ that are either disjoint from $x$ and $y$ or intersect $x$ and $y$ like either 
the curve depicted in the top of Figure \ref{figure:conjugateunderb}.a or the curve depicted in Figure
\ref{figure:conjugateunderb}.b.

Consider $\phi \in S_{U\Sigma}$.  Since $T_b$ fixes $x$ and $y$, the case $\phi = T_b$ is trivial.
We therefore can assume that $\phi = \LPtPsh{\gamma}$ for a based curve $\gamma$ like those described above.  If $\gamma$
is disjoint from $x$ and $y$, then the proof is trivial.  If $\gamma$
is a curve that intersects $x$ and $y$ like the curve in the top of Figure \ref{figure:conjugateunderb}.a, then
$\LPtPsh{\gamma} = T_{\gamma_1,\gamma_2}$ for the curves $\gamma_1$ and $\gamma_2$ shown in the bottom of Figure
\ref{figure:conjugateunderb}.a.  We conclude that using relation \eqref{F.5}, we have
$$T_{\LPtPsh{\gamma}(x),\LPtPsh{\gamma}(y)} = T_{x, T_{\gamma_2}^{-1}(y)} = C_{\gamma_2,y} T_{x,y}.$$
Since $C_{\gamma_2,y} \in K_{g,n}$ (see \S \ref{section:handledragprelim}), this proves the claim.

If instead $\gamma$ is a curve that intersects $x$ and $y$ like the curve in Figure \ref{figure:conjugateunderb}.b, then
observe that $T_{\LPtPsh{\gamma}(x),\LPtPsh{\gamma}(y)} = T_{x',y'}$ for the curves $x'$ and $y'$ depicted in
Figure \ref{figure:conjugateunderb}.c.  Letting $\rho$ and $\eta$ be the other curves in Figure
\ref{figure:conjugateunderb}, there is a lantern relation \eqref{L}
$$T_{\rho} = T_{b,\eta} T_{y',x'} T_{x,y}.$$
Here the four boundary components of the lantern are $\rho$, $b$, $y'$, and $x$.
Using relation \eqref{F.1}, we can rearrange this formula and get
$$T_{x,y} = T_{x',y'} (T_{\eta,b} T_{\rho}).$$
Lemma \ref{lemma:differbybsepbd} says that $T_{\eta,b} T_{\rho}$ is a generator for $K_{g,2}$, so the proof follows.
\end{proof}

\subsubsection{The proof of Lemma \ref{lemma:exactsequencesweak}}
\label{section:proofexactsequencesweak}

We now prove Lemma \ref{lemma:exactsequencesweak}.  Let the boundary component $b \subset \Sigma_{g,n}$ and the maps
$i:\Sigma_{g,n} \hookrightarrow \Sigma_{g,n-1}$ and $i_{\ast} : \Mod_{g,n} \rightarrow \Mod_{g,n-1}$
be as in \S \ref{section:discfilling}.

\begin{proof}[{Proof of Lemma \ref{lemma:exactsequencesweak}}]
Let $S_{g,n}$ be the generating set for $\Gamma_{g,n}$.  Observe that for $n=1,2$, the groups $K_{g,n}$ are 
normal subgroups of $\Gamma_{g,n}$ (this uses the conjugation relations \eqref{F.6}-\eqref{F.8}).  Additionally, they are contained in the kernels of the
disc-filling maps $\Gamma_{g,n} \rightarrow \Gamma_{g,n-1}$.  We will apply Lemma \ref{lemma:exactnesstest}.

We must verify the two conditions of Lemma \ref{lemma:exactnesstest}.  We begin with the second condition (that relations
in $\Gamma_{g,n-1}$ lift to relations in $\Gamma_{g,n}$).  
Observe that $\Sigma_{g,n} \setminus i(\Sigma_{g,n-1})$
is a disc $D$.  What we must show is that for every relation
$$s_1 \cdots s_k=1 \quad \quad (s_j \in S_{g,n-1}^{\pm 1})$$
in $\Gamma_{g,n-1}$ we can homotope the curves involved in the definitions of the $s_j$ so that $D$ is disjoint
from all these curves and so that if we let $\tilde{s}_j$ for $1 \leq j \leq k$ be the generators of $\Gamma_{g,n}$
defined by these curves, then $\tilde{s}_1 \cdots \tilde{s}_k$ is a relation of the same type (lantern, crossed lantern, etc.)
in $\Gamma_{g,n}$.  This is an easy case by case check and the details are left to the reader.

It remains to verify the first condition.  Consider $s,s' \in S_{g,n} \cup \{1\}$ that 
project to the same element of $\Gamma_{g,n-1}$.  We must find $k_1,k_2 \in K_{g,n}$ so that
$s' = k_1 s k_2$ in $\Gamma_{g,n}$.  We first assume that one of $s$ and $s'$ (say $s'$) equals $1$.
Consider the case $n=1$.  If $s$ is a bounding pair map or a separating twist,
then (using Lemma \ref{lemma:differbybsep} if $s$ is a bounding pair map) it follows that 
$s$ is a generator of $K_{g,n}$.  Hence in this case we can take $k_1=k_2=1$.  
Also, if $s$ is a commutator of a simply intersecting pair,
then by Lemma \ref{lemma:commutatorexpress} we can write $s = s_1 \cdots s_k$, where the $s_j$ are
separating twists or bounding pair maps with $i_{\ast}(s_j)=1$.  Hence by the previous case
we have $s_j \in K_{g,n}$, so $s \in K_{g,n}$.  Again we can take $k_1=k_2=1$.  
Now consider the case $n=2$.  It is easy to see that
the generator $s$ cannot be a separating twist or a bounding pair map (the key point is that both
boundary components of $\Sigma_{g,2}$ must lie in the same component of the disconnected surface
one gets when one cuts along the curves defining a separating twist or bounding pair map).  
We conclude that $s$ must be
a commutator of a simply intersecting pair, so by Lemma \ref{lemma:commutatorexpress} 
we can again take $k_1=k_2=1$.

We now assume that neither $s$ nor $s'$ equals $1$.  By Lemma \ref{lemma:generatoridentification1}, it
is enough to show that the appropriate $k_1,k_2 \in K_{g,n}$ exist if $s$ and $s'$ either differ by $b$
or differ by a $b$-push map.  The case that they differ by a $b$-push map being a consequence of
Lemma \ref{lemma:conjugateunderb}, we only need to consider the case that $s$ and $s'$ differ by $b$.  
We first assume that $n=1$.  If $s$ and $s'$ are both
bounding pair maps, then without loss of generality we can assume that $s = T_{x,y}$ and
$s = T_{x',y}$ for curves $x$ and $x'$ that differ by $b$.  By Lemma \ref{lemma:differbybsep}, $\{x,x'\}$ forms a bounding pair
and $T_{x,x'} \in K_{g,n}$, so
relation \eqref{F.2} implies that
$$s = T_{x,y} = T_{x,x'} T_{x',y} = T_{x,x'} s',$$
as desired.  The case where $s$ and $s'$ are both separating twists is dealt with in a similar way, using relation
\eqref{F.4} instead of \eqref{F.2}.  If $s$ is a bounding pair map $T_{x,b}$ and $s'$ is a separating twist 
$T_{x'}$ so that $x$ and $x'$ differ by $b$, then since $n=1$, both $T_x$ and $T_b$ are separating twists, and the proof
is similar to the case that $s$ and $s'$ are both separating twists.  Finally, if $s$ and $s'$ are both
commutators of simply intersecting pairs, then using Lemma \ref{lemma:commutatorexpress} together with
the normality of $K_{g,n}$ we can 
reduce to the previously proven cases

We conclude by considering the case $n=2$.  Observe first that $s$ and $s'$ cannot both be bounding
pair maps or separating twists.  Again, the key point is that the curves defining both $s$ and $s'$ cannot 
separate the boundary components of $\Sigma_{g,2}$.  If $s$ is a bounding pair map and $s'$ is a
separating twist, then $s (s')^{-1}$ is a generator of $K_{g,2}$ (see Lemma \ref{lemma:differbybsepbd}).  Finally,
if $s$ and $s'$ are both commutators of simply intersecting pairs, then using Lemma \ref{lemma:commutatorexpress}, we can
reduce to the previously proven cases.
\end{proof} 

\subsection{Identifying the kernels : Lemma \ref{lemma:kernelidentification}}
\label{section:kernelidentification}

The goal of this section is to prove Lemma \ref{lemma:kernelidentification}, which we recall
says that for $g \geq 2$ the natural maps $K_{g,1} \rightarrow \pi_1(U\Sigma_{g})$ and
$K_{g,2} \rightarrow [\pi_1(\Sigma_{g,1}),\pi_1(\Sigma_{g,1})]$ are isomorphisms..  There are four
parts.
\begin{itemize}
\item In \S \ref{section:modaction}, we record some formulas for the action of the $\Mod_{g,n}$
on $\pi_1(\Sigma_{g,n})$.
\item In \S \ref{section:unittangentbundle}, we construct a new presentation for $\pi_1(U\Sigma_g)$.  Along
the way, we prove Theorem \ref{theorem:surfacegroup}, giving a presentation for $\pi_1(\Sigma_g)$ whose
generating set is the set of all simple closed curves.
\item In \S \ref{section:commutatorsubgroup}, by the same method we construct a new
presentation for $[\pi_1(\Sigma_{g,1}),\pi_1(\Sigma_{g,1})]$.
\item In \S \ref{section:kernelidentificationproof}, we put these ingredients together
to prove Lemma \ref{lemma:kernelidentification}.
\end{itemize}

\subsubsection{The action of the mapping class group on $\pi_1$}
\label{section:modaction}

\Figure{figure:surfacegroup}{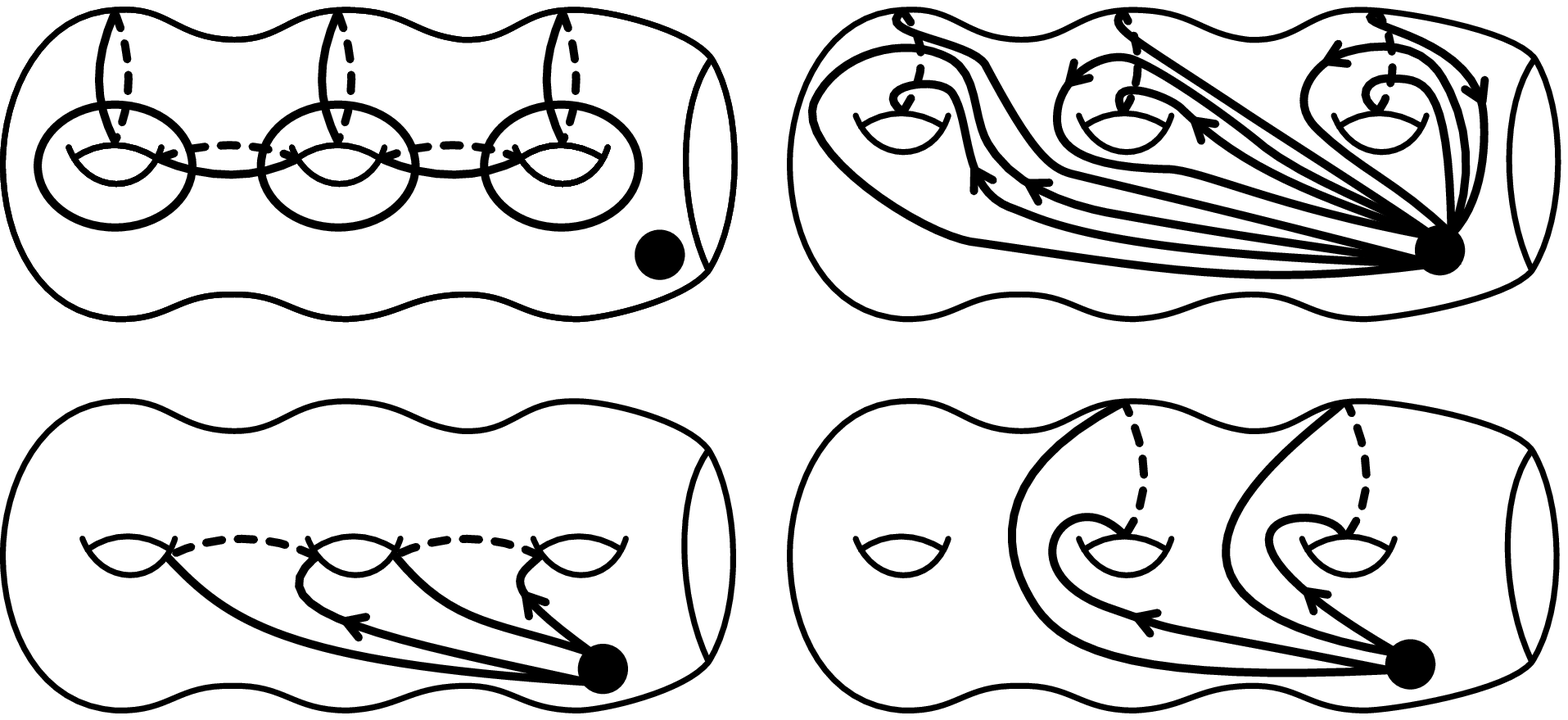}{a. Some curves on $\Sigma_{g,n}^{\ast}$ \CaptionSpace b. The generators for
$\pi_1(\Sigma_{g,n},\ast)$
\CaptionSpace c, d. Extra elements of $\pi_1(\Sigma_{g,n},\ast)$}

\begin{table}
\begin{center}
\begin{tabular}{l|l|l}
$T_{a_i}(\beta_i) = \alpha_i \beta_i$        & $T_{a_i}^{-1} (\beta_i) = \alpha_i^{-1} \beta_i$ &
\\
$T_{b_i} (\alpha_i) = \beta_i^{-1} \alpha_i$ & $T_{b_i}^{-1} (\alpha_i) = \beta_i \alpha_i$     &
\\
$T_{c_i} (\alpha_i) = \gamma_i \alpha_i$     & $T_{c_i} (\beta_i) = \gamma_i \beta_i \gamma_i^{-1}$           & $T_{c_i} (\alpha_{i+1})
= \alpha_{i+1} \gamma_i^{-1}$ \\
$T_{c_i}^{-1} (\alpha_i) = \gamma_i^{-1} \alpha_i$ & $T_{c_i}^{-1} (\beta_i) = \gamma_i^{-1} \beta_i \gamma_i$ & $T_{c_i}^{-1}
(\alpha_{i+1}) = \alpha_{i+1} \gamma_i$ \\
\end{tabular}
\end{center}
\caption{Formulas for the action of $\Mod_{g,n}^{\ast}$ on $\pi_1(\Sigma_{g,n},\ast)$.}
\label{table:modaction}
\end{table}

In this appendix, we record some formulas for the action of certain elements of $\Mod_{g,n}^{\ast}$ on
$\pi_1(\Sigma_{g,n},\ast)$ for $n \leq 1$.  The elements of $\Mod_{g,n}^{\ast}$ we will consider are the
right Dehn twists
$$\{T_{a_1},\ldots,T_{a_g},T_{b_1},\ldots,T_{b_g},T_{c_1},\ldots,T_{c_{g-1}}\},$$
where the curves $a_i$, $b_i$, and $c_i$ are as depicted in Figure \ref{figure:surfacegroup}.a, which
depicts the case $g=3$.  This figure depicts
a surface with one boundary component; our formulas will also hold on a closed surface, where we
interpret all maps as occurring on $\Sigma_{g,1}$ with a disc glued to its boundary component.  
Our generators for $\pi_1(\Sigma_{g,n},\ast)$ are the oriented loops
$$\{\alpha_1,\ldots,\alpha_g,\beta_1,\ldots,\beta_g\}$$
depicted in Figure \ref{figure:surfacegroup}.b in the case $g=3$.  To simplify our formulas, we will make use
of the additional elements
$$\{\gamma_1,\ldots,\gamma_{g-1},\eta_1,\ldots,\eta_{g-1}\} \subset \pi_1(\Sigma_{g,n},\ast)$$
depicted in Figures \ref{figure:surfacegroup}.c and \ref{figure:surfacegroup}.d in the case $g=3$.  
The following formulas express these
additional elements in terms of our generators for $\pi_1(\Sigma_{g,n},\ast)$.
\begin{align*}
\gamma_i &= \eta_i \beta_i^{-1}, \\
\eta_i &= \alpha_{i+1}^{-1} \beta_{i+1} \alpha_{i+1}.
\end{align*}
With these definitions, the formulas in Table \ref{table:modaction} hold.

\subsubsection{A presentation for $\pi_1(U\Sigma_{g})$}
\label{section:unittangentbundle}

We now prove the following.

\begin{proposition}
\label{proposition:unittangentbundle}
Let $\Gamma$ be the group whose generators are the symbols
\begin{align*}  
S=\{T_b\} \cup \{T_{x_1,x_2} \text{ $|$ } & \text{there exists some nontrivial simple closed curve 
$\gamma \in \pi_1(\Sigma_{g})$} \\
                                          & \text{so that $\LPtPsh{\gamma} = T_{x_1,x_2}$} \}
\end{align*}
subject to the relations $(L)$, $(CL)$, $T_{x_1,x_2} T_{x_2,x_1} = 1$, and $[T_b,s]=1$ for all $s \in S$.  Then the
natural map $\Gamma \rightarrow \pi_1(U\Sigma_{g})$ is an isomorphism.
\end{proposition}

\noindent
This will be a consequence of Theorem \ref{theorem:surfacegroup}, which we
now prove.
    
\begin{proof}[{Proof of Theorem \ref{theorem:surfacegroup}}]
Let $S$ be the generating set for $\Gamma$ and let
$$S_{\pi_1}=\{\alpha_1,\ldots,\alpha_g,\beta_1,\ldots,\beta_g\}$$
be the set of generators for $\pi_1(\Sigma_{g},\ast)$ depicted in Figure \ref{figure:surfacegroup} (remember
the convention we discussed in \S \ref{section:modaction} -- since we are working on a closed surface
we view $\Sigma_g$ as the surface $\Sigma_{g,1}$ in Figure \ref{figure:surfacegroup} with a disc attached
to the boundary component).  Observe
that $S_{\pi_1}$ may be naturally identified with the subset
$$S' = \{s_x \text{ $|$ $x \in S_{\pi_1}$}\}$$
of $S$.  By Corollary \ref{corollary:isomorphismtest}, to prove the theorem, it is enough to prove that
$S'$ generates $\Gamma$ and that the $s_x$ satisfies the surface relation
$$[s_{\alpha_1},s_{\beta_1}] \cdots [s_{\alpha_g},s_{\beta_g}]=1.$$
The latter claim follows from the following easy calculation, where we indicate above each $=$ sign the
relation used.
\begin{align*}
[s_{\alpha_1},s_{\beta_1}] \cdots [s_{\alpha_g},s_{\beta_g}] &\stackrel{\overline{CL}}{=} (s_{\alpha_1^{-1} \beta_1^{-1}} s_{\alpha_1 \beta_1}) \cdots (s_{\alpha_g^{-1} \beta_g^{-1}} s_{\alpha_g \beta_g}) \\
                                                             &\stackrel{\overline{L}}{=} s_{[\alpha_1, \beta_1]} \cdots s_{[\alpha_g,\beta_g]} \\
                                                             &\stackrel{\overline{L}}{=} s_{[\alpha_1,\beta_1][\alpha_2,\beta_2]} s_{[\alpha_3,\beta_3]} \cdots s_{[\alpha_g,\beta_g]}\\
                                                             &=\ldots=s_{[\alpha_1,\beta_1] \cdots [\alpha_g,\beta_g]}=1.
\end{align*}
We now prove the former claim.  Observe first that we can express $s_x$ for $x$ a separating curve
as a product of commutators of $s_y$ for nonseparating curves $y$.  Indeed, this is essentially contained
in the above calculation.  Hence $\Gamma$ is generated by
$$S_{\nosep} = \{s_x \text{ $|$ $x \in \pi_1(\Sigma_{g},\ast)$ is a nonseparating simple closed curve}\}.$$
Observe that $\Mod_{g}^{\ast}$ acts on $S_{\nosep}$ and that 
$\Mod_{g}^{\ast} \cdot S' = S_{\nosep}$.  Let
$$S_{\Mod}=\{T_{a_1},\ldots,T_{a_g},T_{b_1},\ldots,T_{b_g},T_{c_1},\ldots,T_{c_{g-1}}\}$$
be the set of generators for $\Mod_{g}^{\ast}$ defined in \S \ref{section:modaction} and let
$$\{\gamma_1,\ldots,\gamma_{g-1},\eta_1,\ldots,\eta_{g-1}\}$$
be the elements of the surface group defined in 
\S \ref{section:modaction}.
By Lemma \ref{lemma:generationtest}, to prove that $S'$ generates $\Gamma$, it is enough to prove that for
$f \in S_{\Mod}^{\pm 1}$ and $s_x \in S'$, the element $s_{f(x)}$ can be expressed
as a product of elements of $(S')^{\pm 1}$.  This is essentially immediate from
the formulas in Table \ref{table:modaction} in \S \ref{section:modaction}.  We give one of the calculations as a example.
Recall that $\gamma_i = \eta_i \beta_i^{-1}$ and $\eta_i = \alpha_{i+1}^{-1} \beta_{i+1} \alpha_{i+1}$.
\begin{align*}
s_{T_{c_i} a_i} &= s_{\gamma_i \alpha_i} \stackrel{\overline{CL}}{=} s_{\gamma_i} s_{\alpha_i} = s_{\eta_i \beta_i^{-1}} s_{\alpha_i} \stackrel{\overline{L}}{=} s_{\eta_i} s_{\beta_i}^{-1} s_{\alpha_i} = s_{\alpha_{i+1}^{-1} \beta_{i+1} \alpha_{i+1}} s_{\beta_i}^{-1} s_{\alpha_i} \\
&\stackrel{\overline{CL}}{=} s_{\alpha_{i+1}^{-1} \beta_{i+1}} s_{\alpha_{i+1}} s_{\beta_i}^{-1} s_{\alpha_i} \stackrel{\overline{CL}}{=} s_{\alpha_{i+1}}^{-1} s_{\beta_{i+1}} s_{\alpha_{i+1}} s_{\beta_i}^{-1} s_{\alpha_i}.   
\end{align*}
The others are similar.
\end{proof}

\noindent
We now prove Proposition \ref{proposition:unittangentbundle}.

\begin{proof}[{Proof of Proposition \ref{proposition:unittangentbundle}}]
Let $\Gamma'$ be the group from Theorem \ref{theorem:surfacegroup}.  Observe that
$\Gamma' \cong \Gamma / \langle T_b \rangle$.  We therefore have the following commutative diagram of exact sequences.

\begin{center}
\begin{tabular}{c@{\hspace{0.05 in}}c@{\hspace{0.05 in}}c@{\hspace{0.05 in}}c@{\hspace{0.05 in}}c@{\hspace{0.05 in}}c@{\hspace{0.05 in}}c@{\hspace{0.05 in}}c@{\hspace{0.05 in}}c}
$1$ & $\rightarrow$ & $\Z$         & $\rightarrow$ & $\Gamma$                         & $\rightarrow$ & $\Gamma'$               & $\rightarrow$ & $1$ \\
    &               & $\parallel$  &               & $\downarrow$                     &               & $\downarrow$            &               &     \\
$1$ & $\rightarrow$ & $\Z$         & $\rightarrow$ & $\pi_1(U\Sigma_{g},\tilde{\ast})$ & $\rightarrow$ & $\pi_1(\Sigma_{g},\ast)$ & $\rightarrow$ & $1$ \\
\end{tabular}
\end{center}

\noindent
By Theorem \ref{theorem:surfacegroup}, the right hand arrow is an isomorphism.  The five lemma therefore
implies that the center arrow is also an isomorphism, as desired.
\end{proof}

\subsubsection{A presentation for $[\pi_1(\Sigma_{g,1}),\pi_1(\Sigma_{g,1})]$}
\label{section:commutatorsubgroup}

Throughout this section, we will assume that $g \geq 1$.  We begin with some definitions (these definitions will not
be used outside of this section).  We define the group $\Gamma$ to be the group
whose generating set is the set of symbols
\begin{align*}
S = \{[x,y]_0 \text{ $|$ } &\text{$x,y \in (\pi_1(\Sigma_{g,1},\ast) \setminus \{1\})$ are completely distinct and can be realized by} \\
                           &\text{simple closed curves that only intersect at the basepoint} \}
\end{align*}
subject to following set of relations.  For simplicity, for $z \in \pi_1(\Sigma_{g,1},\ast)$, we define
$$[x,y]_0^z := [z^{-1} x z, z^{-1} y z]_0 = [\PtPsh{z} (x), \PtPsh{z} (y)]_0.$$
Also, call a set $X \subset \pi_1(\Sigma_{g,1})$ a {\em good} set if the elements
of $X$ are completely distinct, nontrivial, and can be represented by simple closed curves that only intersect at the basepoint.
The first set of relations are the Witt--Hall relations
\begin{equation}
\tag{$\overline{\text{WH}}$}
\label{WH0}
[g_1 g_2, g_3]_0 = [g_1,g_3]_0^{g_2} [g_2,g_3]_0
\end{equation}
for all $g_1,g_2,g_3 \in \pi_1(\Sigma_{g,1})$ so that the sets
$\{g_1,g_2,g_3\}$ and $\{g_1 g_2, g_3\}$ are good.
Next, we will need the commutator shuffle relation
\begin{equation}
\tag{$\overline{\text{CS}}$}
\label{CS0}
[g_1,g_2]_0^{g_3} = [g_3, g_1]_0 [g_3, g_2]_0^{g_1} [g_1,g_2]_0 [g_1,g_3]_0^{g_2} [g_2,g_3]_0
\end{equation}
for all $g_1,g_2,g_3 \in \pi_1(\Sigma_{g,1})$ so that $\{g_1,g_2,g_3\}$ is a good set.
Next, we will need the relation
\begin{equation}
\tag{ID}
\label{ID}
[g_1,g_2]_0 = [g_3,g_4]_0
\end{equation}
for all $g_1,g_2,g_3,g_4 \in \pi_1(\Sigma_{g,1})$ with $[g_1,g_2]=[g_3,g_4]$ (we emphasize that this is equality
in the commutator subgroup; an example of this phenomenon is $[yx,y]=[x,y]$) so 
that the sets $\{g_1,g_2\}$ and $\{g_3,g_4\}$ are good.
Finally, we will need the following relations for all $x,y,z,w \in \pi_1(\Sigma_{g,1},\ast)$ so that
each of the sets $\{x,y\}$ and $\{z,w\}$ are good.
\begin{align*}
[x,y]_0 [y,x]_0 &= 1, \tag{R.1}\label{R1}\\
[z,w]_0^{-1} [x,y]_0 [z,w]_0 &= [x,y]_0^{[z,w]}. \tag{R.2}\label{R2}\\
\end{align*}

Observe the following.

\begin{lemma}
\label{lemma:surjectivecommutator}
The map $[x,y]_0 \mapsto \DComm{x,y}$ induces a surjective homomorphism $\Gamma \rightarrow K_{g,2}$.
\end{lemma}
\begin{proof}
We must check that relations go to relations.  The only relations for which this
is not clear are the relations \eqref{ID}.  Consider such a relation $[g_1,g_2]_0 = [g_3,g_4]_0$.
There are two cases.  In the first, $[g_1,g_2]$ can be represented by a simple closed separating
curve $\gamma$.  By definition $\DComm{g_1,g_2}$ only depends on $\gamma$, so since $[g_3,g_4] = \gamma$
it follows that $\DComm{g_1,g_2} = \DComm{g_3,g_4}$.  In the other case, $\DComm{g_1,g_2} = C_{a,b}$ and
$\DComm{g_3,g_4} = C_{a',b'}$ for simply intersecting pairs $\{a,b\}$ and $\{a',b'\}$.  We might
not have $a = a'$ and $b = b'$, but we must have $C_{a,b} = C_{a',b'}$ in $\Mod_{g,2}$.  Since
the generators of $K_{g,2}$ are {\em mapping classes}, we must have $C_{a,b} = C_{a',b'}$ in $K_{g,2}$,
as desired.
\end{proof}

Let
$$\psi : \Gamma \rightarrow [\pi_1(\Sigma_{g,1},\ast),\pi_1(\Sigma_{g,1},\ast)]$$
be the homomorphism defined on the generators of $\Gamma$ by $\psi([x,y]_0) = [x,y]$.
Our main result will be the following.

\begin{proposition}
\label{proposition:commutatorgenerators}
The map $\psi$ is an isomorphism.
\end{proposition}

The proof will be modeled on the proof of Theorem \ref{theorem:surfacegroup} above.
To that end, we will need a useful free generating set for the commutator subgroup
of the free group $\pi_1(\Sigma_{g,1},\ast)$.  Let 
$$S_{\pi_1}=\{\alpha_1,\ldots,\alpha_{g},\beta_1,\ldots,\beta_{g}\}$$
be the set of generators for $\pi_1(\Sigma_{g,1},\ast)$ described in \S \ref{section:modaction}, and let $\prec$
be any total ordering on $S_{\pi_1}$.  We then have the following theorem of Tomaszewski.
\begin{theorem}[{\cite{Tomaszewski}}]
\label{theorem:commutatorbasis}
The set
\begin{align*}
\{[x,y]^{z_1^{d_1} \cdots z_k^{d_k}} \text{ $|$ } &\text{$x,y \in S_{\pi_1}$, $x \prec y$, $z_i \in S_{\pi_1}$ and $d_i \in \Z$ for all $i$,} \\
                                                              &\text{and $x \preceq z_1 \prec z_2 \prec \ldots \prec z_k$} \},
\end{align*}
is a free generating set for $[\pi_1(\Sigma_{g,1},\ast),\pi_1(\Sigma_{g,1},\ast)]$.
\end{theorem}

The proof of Proposition \ref{proposition:commutatorgenerators} will be preceeded by four lemmas.
For the first, define
\begin{align*}
S_1 = \{[x,y]_0^{z_1^{d_1} \cdots z_k^{d_k}} \text{ $|$ } &\text{$x,y \in S_{\pi_1}$, $x \prec y$, $z_i \in S_{\pi_1}$ and $d_i \in \Z$ for all $i$,} \\
                                                              &\text{and $x \preceq z_1 \prec z_2 \prec \ldots \prec z_k$} \},
\end{align*}
and let $\Gamma'$ be the subgroup of $\Gamma$ generated by $S_1$.  We then have the following.

\begin{lemma}
\label{lemma:s1iso}
The map $\psi$ maps
$\Gamma'$ isomorphically onto $[\pi_1(\Sigma_{g,1},\ast),\pi_1(\Sigma_{g,1},\ast)]$.
\end{lemma}
\begin{proof}
The set $\psi(S_1)$ is the free generating given by Theorem \ref{theorem:commutatorbasis},
so the lemma follows from Corollary \ref{corollary:isomorphismtest}.
\end{proof}

\begin{remark}
No relations were used in the proof of Lemma \ref{lemma:s1iso}!  The purpose
of the relations is to show that $S_1$ generates $\Gamma$.
\end{remark}

Our goal is thus to prove that $\Gamma' = \Gamma$.  Define
\begin{align*}
S_4 := \{[x,y]_0^f &\text{ $|$ $x,y \in S_{\pi_1}$, $x \prec y$, and $f \in \pi_1(\Sigma_{g,1},\ast)$} \}.
\end{align*}
The first step is the following lemma.

\begin{lemma}
\label{lemma:s4in}
$S_4 \subset \Gamma'$.
\end{lemma}
\begin{proof}
\BeginClaims
This will be a three step process.  We will first prove that we can
reorder the generators in the exponents of elements of $S_1$.  Define
\begin{align*}
S_2 = \{[x,y]_0^{z_1^{d_1} \cdots z_k^{d_k}} \text{ $|$ } &\text{$x,y \in S_{\pi_1}$, $x \prec y$, and for all $i$ we  }\\
                                                              &\text{ have $z_i \in S_{\pi_1}$, $d_i \in \Z$, and $x \preceq z_i$} \}.
\end{align*}
\begin{claim}
$S_2 \subset \Gamma'$.
\end{claim}
\BeginClaimProof
Consider $\mu=[x,y]_0^{z_1^{d_1} \cdots z_k^{d_k}} \in S_2$.  Observe that relation \eqref{R2} (from the
definition of $\Gamma$) says that
by conjugating $\mu$ by elements of $S_1$, we may multiply the exponent $z_1^{d_1} \cdots z_k^{d_k}$ of $\mu$
by any element of $[\pi_1(\Sigma_{g,1},\ast),\pi_1(\Sigma_{g,1},\ast)]$ in $\psi(\Gamma')$.  Lemma \ref{lemma:s1iso} 
says that $\psi(\Gamma')$ is the entire commutator subgroup, so we can multiply the exponent of $\mu$ by 
any desired commutator.  By
doing this, we can reorder the terms in it in an arbitrary way.  We conclude that by conjugating $\mu$ by elements of $S_1$, we
can transform it into an element of $S_1$; i.e.\ that $\mu \in \Gamma'$, as desired.
\EndClaimProof
Next, we will show that we can have any generators we want in the exponents (in other words, in the
exponent of $[x,y]_0$ we can have $z$ with $z \prec x$).  Define
\begin{align*}
S_3 = \{[x,y]_0^{z_1^{d_1} \cdots z_k^{d_k}} \text{ $|$ } &\text{$x,y \in S_{\pi_1}$, $x \prec y$, $z_i \in S_{\pi_1}$ and $d_i \in \Z$ for all $i$,} \\
                                                              &\text{ and $z_1 \prec z_2 \prec \cdots \prec z_k$} \}.
\end{align*}
\begin{claim}
$S_3 \subset \Gamma'$.
\end{claim}
\BeginClaimProof
Consider $\mu=[x,y]_0^{z_1^{d_1} \cdots z_k^{d_k}} \in S_3$ with $d_1 \neq 0$.  Set
$$N = \sum_{z_i \prec x} |d_i|.$$
We will prove that $\mu \in \Gamma'$ by induction on $N$.  The base case $N=0$ being a consequence of the 
fact that $S_2 \subset \Gamma'$, we assume that
$N>0$.  We consider the case $d_1>0$; the case $d_1<0$ is exactly the same.  Set $f = z_1^{d_1-1} \cdots z_k^{d_k}$.
Observe that the following is a consequence of \eqref{CS}, \eqref{R1}, and \eqref{R2} (this calculation is the purpose of the commutator
shuffle).
$$\mu = [x, y]_0^{z_1 f} = [z_1, x]_0^f [z_1, y]_0^{x f} [x,y]_0^f [x,z_1]_0^{y f} [y,z_1]_0^f$$
By the relation \eqref{R1}, the $1^{\text{st}}$, $2^{\text{nd}}$, $4^{\text{th}}$ and $5^{\text{th}}$ terms on the right hand 
side or their inverses are in $S_2$, and hence in $\Gamma'$.  Also, by induction, 
the $3^{\text{rd}}$ term is in $\Gamma'$.  We conclude that $\mu \in \Gamma'$, as desired.
\EndClaimProof

An argument identical to the proof that $S_2 \subset \Gamma'$ now establishes that $S_4 \subset \Gamma'$, as desired.
\end{proof}

Now let
$$\{\gamma_1,\ldots,\gamma_{g-1},\eta_1,\ldots,\eta_{g-1}\}$$
be the elements of the surface group defined in \S \ref{section:modaction}.

\begin{lemma}
Fix $1 \leq i \leq g-1$.  For any $x \in S_{\pi_1}$ and $f \in \pi_1(\Sigma_{g,1},\ast)$, 
the group $\Gamma'$ contains $[\gamma_i,x]_0^f$ and $[\eta_i,x]_0^f$.
\end{lemma}
\begin{proof}
The proofs for $[\gamma_i,x]_0^f$ and $[\eta_i,x]_0^f$ are similar.  We will do the case
of $[\eta_i,x]_0^f$ and leave the other case to the reader.  Assume first that
$x \neq \alpha_{i+1},\beta_{i+1}$.  Since $\eta_i = \alpha_{i+1}^{-1} \beta_{i+1} \alpha_{i+1}$, 
we can perform the following calculation.
\begin{align*}
[\eta_i,x]_0^f &= [\alpha_{i+1}^{-1} \beta_{i+1} \alpha_{i+1},x]_0^f \stackrel{\text{\ref{WH0}}}{=} [\alpha_{i+1}^{-1} \beta_{i+1},x]_0^{\alpha_{i+1} f} [\alpha_{i+1},x]_0^f \\
&\stackrel{\text{\ref{WH0}}}{=} [\alpha_{i+1}^{-1},x]_0^{\beta_{i+1} \alpha_{i+1} x} [\beta_{i+1},x]_0^{\alpha_{i+1} f} [\alpha_{i+1},x]_0^f \\
&\stackrel{\text{\ref{ID}}}{=} [x,\alpha_{i+1}]_0^{\alpha_{i+1}^{-1} \beta_{i+1} \alpha_{i+1} x} [\beta_{i+1},x]_0^{\alpha_{i+1} f} [\alpha_{i+1},x]_0^f.
\end{align*}
Each of these terms is in $S_4$, so by Lemma \ref{lemma:s4in} we conclude that $[\eta_i,x]_0^f \in \Gamma'$, as desired.
Next, if $x = \alpha_{i+1}$ we have $[\eta_i,x]_0^f = [\beta_{i+1},\alpha_{i+1}]_0^{\alpha_{i+1} f} \in S_4$, so the
lemma is trivially true.  Finally, if $x = \beta_{i+1}$, then we have the following calculation.
\begin{align*}
[\eta_i,x]_0^f &= [\alpha_{i+1}^{-1} \beta_{i+1} \alpha_{i+1},\beta_{i+1}]_0^f \stackrel{\text{\ref{WH0}}}{=} [\alpha_{i+1}^{-1} \beta_{i+1}, \beta_{i+1}]_0^{\alpha_{i+1} f} [\alpha_{i+1}, \beta_{i+1}]_0^f\\
&\stackrel{\text{\ref{ID}}}{=} [\alpha_{i+1}^{-1}, \beta_{i+1}]_0^{\beta_{i+1} \alpha_{i+1} f} [\alpha_{i+1}, \beta_{i+1}]_0^f \stackrel{\text{\ref{ID}}}{=} [\beta_{i+1},\alpha_{i+1}]_0^{\alpha_{i+1}^{-1} \beta_{i+1} \alpha_{i+1} f} [\alpha_{i+1}, \beta_{i+1}]_0^f
\end{align*}
Again, each of these terms is in $S_4$, so by Lemma \ref{lemma:s4in} we are done.
\end{proof}

\begin{lemma}
\label{lemma:sprimegenerates}
Let $\Mod_{g,1}^{\ast}$ act on $\Gamma$ in the natural way.  Then $\Mod_{g,1}^{\ast} \cdot S_1$
generates $\Gamma$.
\end{lemma}
\begin{proof}
Observe that $[\alpha_1,\beta_1]$ can be realized by a simple closed separating curve 
which cuts off a subsurface homeomorphic to $\Sigma_{1,1}$.  By the classification of
surfaces, $\Mod_{g,1}^{\ast}$ acts transitively on such curves (ignoring their
orientations).  Hence for every $\rho \in \pi_1(\Sigma_{g,1},\ast)$
that can be realized by a simple closed separating curve that cuts off a subsurface homeomorphic
to $\Sigma_{1,1}$, there is some $[\alpha,\beta]_0 \in \Mod_{g,1}^{\ast} \cdot S_1$ so
that either $[\alpha,\beta] = \rho$ or $[\alpha,\beta] = \rho^{-1}$.  Combining Lemma \ref{lemma:commutatorgenerators} with the
relations \eqref{ID}, \eqref{WH0}, and \eqref{R1}, we conclude that every generator of $\Gamma$ or
its inverse is contained in the subgroup generated by $\Mod_{g,1}^{\ast} \cdot S_1$, as desired. 
\end{proof}

\begin{proof}[{Proof of Proposition \ref{proposition:commutatorgenerators}}]
Recall that our goal is to show that $\Gamma = \Gamma'$.  We will now use Lemma \ref{lemma:generationtest}.
Consider the natural action of $\Mod_{g,1}^{\ast}$ on $\Gamma$.
By Lemmas \ref{lemma:sprimegenerates}, \ref{lemma:s1iso} and \ref{lemma:generationtest}, to prove
that $\Gamma = \Gamma'$ it is enough to find some set of generators for $\Mod_{g,1}^{\ast}$ that takes
$S_1$ into $\Gamma'$.  Recall that $\Mod_{g,1}^{\ast}$ fits into the Birman exact sequence
$$1 \longrightarrow \pi_1(\Sigma_{g,1},\ast) \longrightarrow \Mod_{g,1}^{\ast} \longrightarrow \Mod_{g,1} \longrightarrow 1.$$
Now, the kernel $\pi_1(\Sigma_{g,1},\ast)$ acts on $S_1$ by conjugation.  Since $S_4$ contains all conjugates (by elements
of the surface group) of elements of $S_1$, by Lemma \ref{lemma:s4in} it is enough to find some set of elements of
$\Mod_{g,1}^{\ast}$ which project to generators for $\Mod_{g,1}$ and that take $S_1$ into $\Gamma'$.  Let
$$S_{\Mod}=\{T_{a_1},\ldots,T_{a_g},T_{b_1},\ldots,T_{b_g},T_{c_1},\ldots,T_{c_{g-1}}\}$$
be the elements of $\Mod_{g,1}^{\ast}$ from \S \ref{section:modaction}.  Observe
that $S_{\Mod}$ projects to a set of generators for $\Mod_{g,1}$.  We conclude
by observing that the formulas in Table \ref{table:modaction} in \S \ref{section:modaction} imply that
$S_{\Mod}^{\pm 1} (S_1) \subset \Gamma'$; the calculations
are similar to the ones that showed that $[\eta_i,x]_0^f \in \Gamma'$.
\end{proof}

\subsubsection{The proof of Lemma \ref{lemma:kernelidentification}}
\label{section:kernelidentificationproof}

We now prove Lemma \ref{lemma:kernelidentification}, completing the proof of Proposition \ref{proposition:exactsequences}.

\begin{proof}[{Proof of Lemma \ref{lemma:kernelidentification}}]
Observe that Proposition \ref{proposition:unittangentbundle}
tells us that $K_{g,1}$ is a quotient of $\pi_1(U\Sigma_{g})$.  Since the map $\Gamma_{g,1} \rightarrow \Torelli_{g}$
fits into the commutative diagram
\begin{center}
\begin{tabular}{c@{\hspace{0.05 in}}c@{\hspace{0.05 in}}c@{\hspace{0.05 in}}c@{\hspace{0.05 in}}c@{\hspace{0.05 in}}c@{\hspace{0.05 in}}c@{\hspace{0.05 in}}c@{\hspace{0.05 in}}c}
$1$ & $\rightarrow$ & $K_{g,1}$         & $\rightarrow$ & $\Gamma_{g,1}$                         & $\rightarrow$ & $\Gamma_{g}$               & $\rightarrow$ & $1$ \\
    &               & $\downarrow$  &               & $\downarrow$                     &               & $\downarrow$            &               &     \\
$1$ & $\rightarrow$ & $\pi_1(U\Sigma_{g})$         & $\rightarrow$ & $\Torelli_{g,1}$ & $\rightarrow$ & $\Torelli_{g}$ & $\rightarrow$ & $1$ \\
\end{tabular}
\end{center}
we conclude that in fact $K_{g,1} \cong \pi_1(U\Sigma_{g})$.  In a similar way (using 
Lemma \ref{lemma:surjectivecommutator} and Proposition \ref{proposition:commutatorgenerators}
instead of Proposition \ref{proposition:unittangentbundle}), we prove that 
$K_{g,2} \cong [\pi_1(\Sigma_{g,1}),\pi_1(\Sigma_{g,1})]$, as desired.
\end{proof}

\section{The proof of Proposition \ref{proposition:curvestorelliconnected}}
\label{section:curvestorelliconnected}

This section is devoted to the proof of Proposition \ref{proposition:curvestorelliconnected},
which we recall has the following two conclusions for $g \geq 1$.  
\begin{enumerate}
\item The complex $\ModCurves_g$ is $(g-2)$-connected.
\item The complex $\ModCurves_g / \Torelli_{g}$ is $(g-1)$-connected.
\end{enumerate}
We begin in \S \ref{section:simplicialprelims} with some preliminary material on simplicial complexes.  Next,
in \S \ref{section:curvestorellifirst} we recall the definition of $\ModCurves_{g}$ and
prove the first conclusion of Proposition \ref{proposition:curvestorelliconnected}.  Next, in
\S \ref{section:curvestorellisecond} we give a linear-algebraic reformulation of the second
conclusion of Proposition \ref{proposition:curvestorelliconnected}.  The skeleton of the proof
of this linear-algebraic reformulation is contained in \S \ref{section:proofskeleton}.  This
proof depends on a proposition whose proof is contained in 
\S \ref{section:mainproposition12} - \S \ref{section:mainproposition4}.

\begin{remark}
The proof shares many ideas with the proof of \cite[Theorem 5.3]{PutmanCutPaste}, though
the details are more complicated.
\end{remark}

\subsection{Generalities about simplicial complexes}
\label{section:simplicialprelims}

Our basic reference for simplicial complexes is \cite[Chapter 3]{Spanier}.  Let us recall the definition
of a simplicial complex given there.

\begin{definition}
A {\em simplicial complex} $X$ is a set of nonempty finite sets (called {\em simplices}) so that 
if $\Delta \in X$ and $\emptyset \neq \Delta' \subset \Delta$, then $\Delta' \in X$.  
If $\Delta,\Delta' \in X$ and $\Delta' \subset \Delta$, 
then we will say that $\Delta'$ is a {\em face} of $\Delta$.  The {\em dimension} of a simplex
$\Delta \in X$ is $|\Delta|-1$ and is denoted $\Dim(\Delta)$.  A simplex of dimension $0$ will be called a {\em vertex}
and a simplex of dimension $1$ will be called an {\em edge}; we will abuse notation and confuse
a vertex $\{v\} \in X$ with the element $v$.
For $k \geq 0$, the subcomplex of $X$ consisting of all simplices
of dimension at most $k$ (known as the {\em $k$-skeleton of $X$}) will be denoted $X^{(k)}$.
If $X$ and $Y$ are simplicial complexes, then a {\em simplicial map} from $X$ to $Y$ is a function
$f : X^{(0)} \rightarrow Y^{(0)}$ so that if $\Delta \in X$, then $f(\Delta) \in Y$.
\end{definition}

If $X$ is a simplicial complex, then we will define 
the geometric realization $|X|$ of $X$ in the standard way (see \cite[Chapter 3]{Spanier}).  When
we say that $X$ has some topological property (e.g.\ simple-connectivity), we will mean that $|X|$ possesses
that property.
  
Next, we will need the following definitions.

\begin{definition}
Consider a simplex $\Delta$ of a simplicial complex $X$.
\begin{itemize}
\item The {\em star} of $\Delta$ (denoted $\Star_X(\Delta)$) is the subcomplex of $X$ consisting
of all $\Delta' \in X$ so that there is some $\Delta'' \in X$ with $\Delta,\Delta' \subset \Delta''$.
By convention, we will also define $\Star_X(\emptyset) = X$.
\item The {\em link} of $\Delta$ (denoted $\Link_X(\Delta)$) is the subcomplex of $\Star_X(\Delta)$
consisting of all simplices that do not intersect $\Delta$.  By convention, we will also define
$\Link_X(\emptyset) = X$.
\end{itemize}
If $X$ and $Y$ are simplicial complexes, then the {\em join} of $X$ and $Y$ (denoted
$X \ast Y$) is the simplicial complex whose simplices are all sets $\Delta \sqcup \Delta'$ satisfying the following.
\begin{itemize}
\item $\Delta$ is either $\emptyset$ or a simplex of $X$.
\item $\Delta'$ is either $\emptyset$ or a simplex of $Y$.
\item One of $\Delta$ or $\Delta'$ is nonempty.
\end{itemize}
Observe that $\Star_X(\Delta) = \Delta \ast \Link_X(\Delta)$ (this is true even if $\Delta = \emptyset$).  
\end{definition}

For $n \leq -1$, we will say that the empty set is both an $n$-sphere and a closed $n$-ball.  Also, if $X$
is a space then we will say that $\pi_{-1}(X)=0$ if $X$ is nonempty and that $\pi_{k}(X)=0$ for all $k \leq -2$.
With these conventions, it is true for all $n \in \Z$ that
a space $X$ satisfies $\pi_n(X)=0$ if and only if every map of an $n$-sphere into $X$ can be extended to a map
of a closed $(n+1)$-ball into $X$.

Finally, we will need the following definition.  A basic reference is \cite{RourkeSanderson}.  

\begin{definition}
For $n \geq 0$, a {\em combinatorial $n$-manifold} $M$ is a nonempty simplicial complex that
satisfies the following inductive property.  If $\Delta \in M$, then $\Dim(\Delta) \leq n$. 
Additionally, if $n-\Dim(\Delta)-1 \geq 0$, then $\Link_M(\Delta)$ is a combinatorial $(n-\Dim(\Delta)-1)$-manifold
homeomorphic to either an $(n-\Dim(\Delta)-1)$-sphere or a closed $(n-\Dim(\Delta)-1)$-ball.  We will denote by
$\partial M$ the subcomplex of $M$ consisting of all simplices $\Delta$ so that $\Dim(\Delta) < n$ and so that
$\Link_M(\Delta)$ is homeomorphic to a closed $(n-\Dim(\Delta)-1)$-ball.
If $\partial M = \emptyset$ then $M$  is said to be {\em closed}.  A
combinatorial $n$-manifold homeomorphic to an $n$-sphere (resp. a closed $n$-ball) will be called
a {\em combinatorial $n$-sphere} (resp. a {\em combinatorial $n$-ball}).
\end{definition}

It is well-known that if $\partial M \neq \emptyset$, then $\partial M$ is a closed combinatorial $(n-1)$-manifold
and that if $B$ is
a combinatorial $n$-ball, then $\partial B$ is a combinatorial $(n-1)$-sphere.  Also, if $M_1$
and $M_2$ are combinatorial manifolds and if $M_1 \times M_2$ is the standard triangulation
of $|M_1| \times |M_2|$, then $M_1 \times M_2$ is a combinatorial manifold.  Finally, subdivisions
of combinatorial manifolds are combinatorial manifolds.  

\begin{warning}
There exist simplicial complexes that are homeomorphic to manifolds but are {\em not} combinatorial
manifolds.
\end{warning}

The following is an immediate consequence of the Zeeman's extension \cite{ZeemanSimp} of
the simplicial approximation theorem.

\begin{lemma}
\label{lemma:simpapprox}
Let $X$ be a simplicial complex and $n \geq 0$.  The following hold.
\begin{enumerate}
\item Every element of $\pi_n(X)$ is represented by a simplicial map $S \rightarrow X$, 
where $S$ is a combinatorial $n$-sphere.  
\item If $S$ is a combinatorial $n$-sphere and $f : S \rightarrow X$ is a nullhomotopic 
simplicial map, then there is a combinatorial
$(n+1)$-ball $B$ with $\partial B = S$ and a simplicial map $g : B \rightarrow X$ so that $g|_{S} = f$.
\end{enumerate}
\end{lemma}

A consequence of the first conclusion of Lemma \ref{lemma:simpapprox} is that we can 
prove that simplicial complexes are $n$-connected by attempting to simplicially
homotope maps of combinatorial $n$-spheres to constant maps.  
The basic move by which we will do this is the following (see Figure \ref{figure:linkmove}
for examples).

\Figure{figure:linkmove}{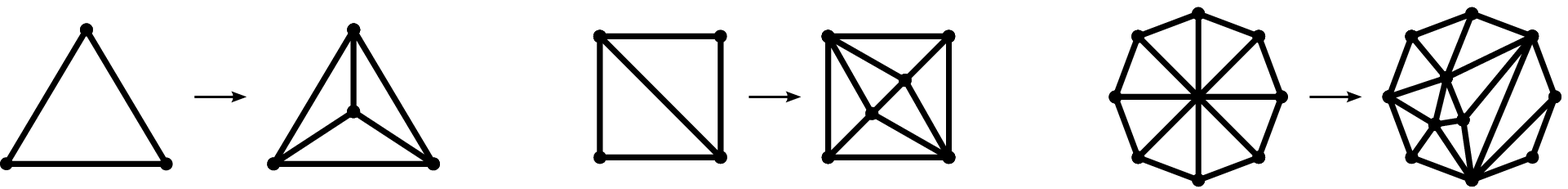}{Effect on $\Star_{S}(\Delta) \subset S$ of a link move with
a. $\Dim(\Delta) = 2$ \CaptionSpace b. $\Dim(\Delta)=1$ \CaptionSpace c. $\Dim(\Delta)=0$.}

\begin{definition}
Let $\phi : S \rightarrow X$ be a simplicial map of a combinatorial $n$-sphere into
a simplicial complex.  For some $\Delta \in S$, let $T$ be a combinatorial
$(n-\Dim(\Delta))$-ball so that $\partial T = \Link_S(\Delta)$ and let
$f : T \rightarrow \Star_X(\phi(\Delta))$ be a simplicial map so that $f|_{\partial T} = \phi|_{\Link_S(\Delta)}$.  Define
$S'$ to be $S$ with $\Star_{S}(\Delta)$ replaced with $T$ and define $\phi' : S' \rightarrow X$ in the following
way.  For $v \in (S')^{(0)} \setminus T^{(0)}$, define $\phi'(v) = \phi(v)$.  For $v \in B^{(0)}$, define
$\phi'(v) = f(v)$.  Observe that $\phi'$ extends linearly to a simplicial map.  We
will call $\phi' : S' \rightarrow X$ the result of performing a {\em link move} to $\phi : S \rightarrow X$ on $\Delta$
with $f$.
\end{definition}

Observe that if a map $S' \rightarrow X$ is the result of performing a link move
on a map $S \rightarrow X$, then $|S'|$ is naturally homeomorphic to $|S|$ and the
induced maps $|S'| \rightarrow |X|$ and $|S| \rightarrow |X|$ are homotopic.

\subsection{$\ModCurves_g$ and the proof of the first conclusion
of Proposition \ref{proposition:curvestorelliconnected}}
\label{section:curvestorellifirst}

We begin by recalling the definition of $\ModCurves_g$ and giving names
to the various types of simplices.

\Figure{figure:modifiedcomplex2}{ModifiedComplex}{a. A standard simplex \CaptionSpace b. A simplex of type $\sigma$ \CaptionSpace c. A simplex of type $\delta$}

\begin{definition}
The complex $\ModCurves_g$ is the simplicial complex whose $(k-1)$-simplices are sets
$\{\gamma_1,\ldots,\gamma_k\}$ of isotopy classes of simple closed nonseparating curves
on $\Sigma_{g}$ satisfying one of the following three conditions (for some ordering of the $\gamma_i$).
\begin{itemize}
\item The $\gamma_i$ are disjoint and $\gamma_1 \cup \cdots \cup \gamma_k$ does not separate
$\Sigma_{g}$ (see Figure \ref{figure:modifiedcomplex2}.a).  These will be called the {\em standard simplices}.
\item The $\gamma_i$ satisfy
$$\GeomI(\gamma_i,\gamma_j) = \begin{cases}
        1 & \text{if $(i,j)=(1,2)$} \\
        0 & \text{otherwise}
\end{cases}$$
and $\gamma_1 \cup \cdots \cup \gamma_k$ does not separate $\Sigma_{g}$ (see Figure \ref{figure:modifiedcomplex2}.b).
These will be called {\em simplices of type $\sigma$}.
\item The $\gamma_i$ are disjoint, $\gamma_1 \cup \gamma_2 \cup \gamma_3$ cuts off a copy of $\Sigma_{0,3}$
from $\Sigma_{g}$, and $\{\gamma_1,\ldots,\gamma_k\} \setminus \{\gamma_3\}$ is a standard simplex (see
Figure \ref{figure:modifiedcomplex2}.c).  These will be called {\em simplices of type $\delta$}.
\end{itemize}
\end{definition}

We now wish to prove the first conclusion of Proposition \ref{proposition:curvestorelliconnected},
which we recall says that $\ModCurves_g$ is $(g-2)$-connected.
We will need the following theorem of Harer.  Recall that $\CNoSep_{g,n}$ is the simplicial
complex whose $(k-1)$-simplices are sets $\{\gamma_1,\ldots,\gamma_k\}$ of isotopy classes
of simple closed curves on $\Sigma_{g,n}$ which can be realized so that 
$\Sigma_{g,n} \setminus (\gamma_1 \cup \cdots \cup \gamma_k)$ is connected.

\begin{theorem}[{\cite[Theorem 1.1]{HarerStability}}]
\label{theorem:harerconnected}
For $g \geq 1$ and $n \geq 0$, the complex $\CNoSep_{g,n}$ is $(g-2)$-connected.
\end{theorem}

\begin{proof}[{Proof of Proposition \ref{proposition:curvestorelliconnected}, first conclusion}]
For some $-1 \leq i \leq g-2$, let $S$ be a combinatorial $i$-sphere (remember our conventions
about the $(-1)$-sphere!) and let
$\phi : S \rightarrow \ModCurves_g$ be a simplicial map.
By Lemma \ref{lemma:simpapprox} and Theorem \ref{theorem:harerconnected}, it
is enough to homotope $\phi$ so that $\phi(S) \subset \CNoSep_{g}$.  If
$e \in S^{(1)}$ is such that $\phi(e)$ is a 1-simplex of type $\sigma$ (this
can only happen if $i \geq 1$), then
$\Sigma_{g}$ cut along the curves in $\phi(e)$ is homeomorphic to $\Sigma_{g-1,1}$.  This implies
that $\phi(\Link_S(e)) \subset \Link_{\ModCurves_g}(\phi(e)) \cong \CNoSep_{g-1,1}$.
Now, $\Link_S(e)$ is a combinatorial $(i-2)$-sphere, so Theorem \ref{theorem:harerconnected} and Lemma \ref{lemma:simpapprox}
imply that there is some map $f : B \rightarrow \Link_{\ModCurves_g}(\phi(e))$, where
$B$ is a combinatorial $(i-1)$-ball with $\partial B = \Link_S(e)$ and $f|_{\partial B} = \phi|_{\Link_S(e)}$.  
We can therefore perform a link move to $\phi$ on $e$ with $f$, eliminating $e$.
This allows us to remove all simplices of $S$ mapping to simplices of type $\sigma$.  A similar argument
allows us to remove all simplices of $S$ mapping to simplices of type $\delta$, and we are done.
\end{proof}

\subsection{A linear-algebraic reformulation of the second conclusion of Proposition \ref{proposition:curvestorelliconnected}}
\label{section:curvestorellisecond} 

The second conclusion of Proposition \ref{proposition:curvestorelliconnected} asserts that
$\ModCurves_g / \Torelli_g$ is $(g-1)$-connected.  In this section, we will reformulate this
by giving a concrete description of $\ModCurves_g / \Torelli_{g}$.  One
obvious thing associated to a nonseparating curve $\gamma$ on $\Sigma_{g}$ that is invariant
under $\Torelli_{g}$ is the 1-dimensional submodule $\Span{$[\gamma]$}$ of $\HH_1(\Sigma_{g})$
(the vector $[\gamma]$ is not well-defined since $\gamma$ is unoriented).  Now, $\Span{$[\gamma]$}$ is
not an arbitrary submodule of $\HH_1(\Sigma_{g})$ : since $\{[\gamma]\}$ can be completed
to a symplectic basis for $\HH_1(\Sigma_{g})$, it follows that $\Span{$[\gamma]$}$ is 
actually a 1-dimensional summand of $\HH_1(\Sigma_{g})$.
The following definition is meant to mimic the definition of $\ModCurves_g$ in terms
of summands of $\HH_1(\Sigma_{g})$.

\begin{definition}
A subspace $X$ of $\HH_1(\Sigma_g)$ is {\em isotropic} if $i(x,y)=0$ for all $x,y \in X$.  The {\em genus $g$ complex
of unimodular isotropic lines}, denoted $\Lines(g)$, is the simplicial complex whose $(k-1)$-simplices
are sets $\{L_1,\ldots,L_k\}$ of $1$-dimensional summands $L_i$ of $\HH_1(\Sigma_g)$ so that
$\Span{$L_1,\ldots,L_k$}$ is a $k$-dimensional isotropic summand of $\HH_1(\Sigma_g)$.  These
will be called the {\em standard simplices}.  Now consider a set
$\Delta=\{\Span{$v_1$},\ldots,\Span{$v_k$}\} \subset (\Lines(g))^{(0)}$.
\begin{itemize}
\item $\Delta$ forms a {\em simplex of type $\sigma$} if
$$\AlgI(v_i,v_j) = \begin{cases}
        \pm 1 & \text{if $(i,j)=(1,2)$} \\
        0     & \text{otherwise}
\end{cases}$$
and $\Span{$v_1,\ldots,v_k$}$ is a $k$-dimensional summand of $\HH_1(\Sigma_{g})$.
\item $\Delta$ forms a {\em simplex of type $\delta$} if $v_3 = \pm v_1 \pm v_2$ and $\Delta \setminus \{\Span{$v_3$}\}$ is
a standard simplex.
\end{itemize}
We will denote $\Lines(g)$ with all simplices of type $\sigma$ and $\delta$ attached by $\Lines_{\sigma,\delta}(g)$.  Similarly,
$\Lines_{\sigma}(g)$ (resp. $\Lines_{\delta}(g)$) will denote $\Lines(g)$ with all simplices of type $\sigma$ (resp. $\delta$)
attached.
\end{definition}

The map $\gamma \mapsto \Span{$[\gamma]$}$ induces a map 
$\pi : \ModCurves_g / \Torelli_g \rightarrow \Lines_{\sigma,\delta}(g)$
that is invariant under the action of $\Torelli_{g}$ and preserves the types of simplices.  We now prove the 
following (this generalizes \cite[Lemma 6.2]{PutmanCutPaste}).

\begin{lemma}
\label{lemma:modcurveslines}
For $g \geq 1$, the map $\pi$ induces an isomorphism from 
$\ModCurves_g / \Torelli_{g}$ to $\Lines_{\sigma,\delta}(g)$.
\end{lemma}

For the proof of Lemma \ref{lemma:modcurveslines}, we will need the following lemma (cf.\ \cite[{Lemma 8.3}]{PutmanCutPaste}).

\Figure{figure:basisrealize}{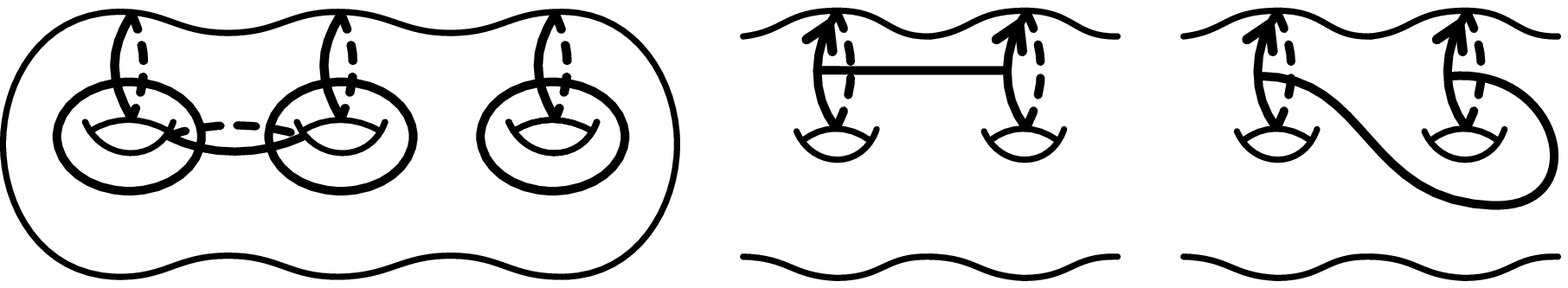}{a. Curves used in Lemma \ref{lemma:basisrealize} \CaptionSpace b,c. With an appropriate choice of orientation, 
a component of $\alpha_{h+1}' \cup \ell \cup \alpha_n$ is homologous to the following : (b) $[\alpha_{h+1}'] - [\alpha_n]$, (c) $[\alpha_{h+1}'] +
[\alpha_n]$.}

\begin{lemma}
\label{lemma:basisrealize}
Let $g \geq 1$, let $0 \leq k \leq h < g$, let $\{a_1,b_1,\ldots,a_g,b_g\}$ be a symplectic basis for $\HH_1(\Sigma_g)$, and
let $\{\alpha_1,\ldots,\alpha_h,\beta_1,\ldots,\beta_k\}$ be a set of oriented simple closed
curves on $\Sigma_g$.  If $h \geq 2$, then we are also possibly given some curve $\alpha_{1,2}$.  Assume that
our curves satisfy the following conditions for $1 \leq i,i' \leq h$ and $1 \leq j,j' \leq k$ (see Figure \ref{figure:basisrealize}.a).
\begin{enumerate}
\item $[\alpha_i] = a_i$ and $[\beta_j] = b_j$
\item $\GeomI(\alpha_i,\alpha_{i'}) = \GeomI(\beta_j,\beta_{j'}) = 0$.  
Also, $\GeomI(\alpha_i,\beta_j)$ is $1$ if $i=j$ and is $0$ otherwise.  Finally, if $\alpha_{1,2}$ is given,
then $\GeomI(\alpha_i,\alpha_{1,2}) = 0$ and $\GeomI(\beta_j,\alpha_{1,2})$ 
is $1$ if $1 \leq j \leq 2$ and is $0$ otherwise.
\item If $\alpha_{1,2}$ is given, then $\alpha_1 \cup \alpha_2 \cup \alpha_{1,2}$ 
separates $\Sigma_g$ into two components, one of which
is homeomorphic to $\Sigma_{0,3}$.
\end{enumerate}
Then there exists oriented curves $\{\alpha_{h+1},\ldots,\alpha_g,\beta_{k+1},\ldots,\beta_g\}$ 
so that that the above three conditions
are satisfied for all $1 \leq i,i' \leq g$ and $1 \leq j,j' \leq g$.
\end{lemma}
\begin{proof}
Let $S$ be $\{\alpha_1,\ldots,\alpha_h,\beta_1,\ldots,\beta_k\}$ together with $\alpha_{1,2}$ if
it is given.  Assume first that $h < g$.  We will show how to find $\alpha_{h+1}$.  Let $\Sigma'$ the component
of $\Sigma_g$ cut along the curves in $S$ whose genus is positive and let 
$i : \Sigma' \rightarrow \Sigma_g$ be the inclusion.  
If $i_{\ast} : \HH_1(\Sigma') \rightarrow \HH_1(\Sigma_g)$ is the induced map, then 
$i_{\ast}(\HH_1(\Sigma')) = [S]^{\perp}$, where
by $[S]^{\perp}$ we mean the subspace of $\HH_1(\Sigma_g)$ consisting of all vectors orthogonal with respect
to the algebraic intersection form to the homology classes of all
the curves in $S$.
Next, let $\Sigma''$ be the surface that results from gluing discs to all boundary components of $\Sigma'$,
let $i' : \Sigma' \hookrightarrow \Sigma''$ be the inclusion, and let
$i'_{\ast} : \HH_1(\Sigma') \rightarrow \HH_1(\Sigma'')$ be the induced map.  Let $\tilde{a}_{h+1} \in \HH_1(\Sigma')$
be a primitive vector so that $i_{\ast}(\tilde{a}_{h+1}) = a_{h+1}$ and let $\overline{a}_{h+1}:= i'_{\ast}(\tilde{a}_{h+1})$.  
Then $\overline{a}_{h+1} \in \HH_1(\Sigma'')$ is a primitive vector in the first
homology group of a closed surface, so there exists some
simple closed curve $\overline{\alpha}_{h+1}$ on $\Sigma''$ so that $[\overline{\alpha}_{h+1}] = \overline{a}_{h+1}$.  We then
isotope $\overline{\alpha}_{h+1}$ so that it lies in $\Sigma' \hookrightarrow \Sigma''$, define
$\tilde{\alpha}_{h+1}$ to be the preimage of $\overline{\alpha}_{h+1}$ in $\Sigma'$, and define
$\alpha_{h+1}'$ to be the image in $\Sigma_g$ of $\tilde{\alpha}_{h+1}$ under the map $i$.

Observe that $[\alpha'_{h+1}] - a_{h+1} \in i_{\ast}(\Ker(i'_{\ast}))$.  Also, since $\Ker(i'_{\ast})$ is generated
by the homology classes of the boundary components of $\Sigma'$, it follows that 
$i_{\ast}(\Ker(i'_{\ast})) = \Span{$a_{k+1},\ldots,a_h$}$.
Thus there exists some $c_{k+1},\ldots,c_h \in \Z$ so that
$[\alpha'_{h+1}] = a_{h+1} + \sum_{j=k+1}^h c_j a_j$.  Assume that $\alpha'_{h+1}$ is chosen so that
$\sum_{j=k+1}^h |c_j]$ is as small as possible.  We claim that all the $c_j$ are zero.  Indeed, assume that
$c_n \neq 0$ for some $k+1 \leq n \leq h$.  We can then (see Figures \ref{figure:basisrealize}.b--c) find
some arc $\ell$ on $\Sigma_g$ satisfying the following three properties.
\begin{itemize}
\item One of the two points of $\partial \ell$ lies on $\alpha_{h+1}'$ and the other lies on $\alpha_n$.
\item $\Interior(\ell)$ is disjoint from every curve in $S$
\item Letting $e$ equal $-1$ if $c_n > 0$ and $1$ if $c_n < 0$, a boundary component $\alpha_{h+1}''$ of
a regular neighborhood of $\alpha_{h+1} \cup \ell \cup \alpha_n$ is homologous to $[\alpha_{h+1}] + e [\alpha_n]$.
\end{itemize}
We can then replace $\alpha_{h+1}'$ with $\alpha_{h+1}''$ and reduce $\sum_{j=k+1}^h |c_j|$, a contradiction.

We can therefore assume that $h = g$.  Assuming now that $k < g$, our goal is to show how to find $\beta_{k+1}$.
Let $\beta_{k+1}$ be some curve so that the
set $\{\alpha_1,\ldots,\alpha_h,\beta_1,\ldots,\beta_k,\beta_{k+1}\}$ (plus $\alpha_{1,2}$ if it is given)
satisfies conditions 2--3 but not necessarily
condition 1.  From the conditions on the geometric intersection number, it follows that
$[\beta_{k+1}] = b_{k+1} + \sum_{n=k+1}^{g} d_n a_n$ for some $d_{k+1},\ldots,d_g \in \Z$.  Choose
$\beta_{k+1}$ so that $\sum_{n=k+1}^g |d_n|$ is as small as possible.  We claim that $d_n = 0$ for all
$k+1 \leq n \leq g$.  Indeed, assume that $d_m \neq 0$ for some $k+1 \leq m \leq g$.  If $m = k+1$, then
we can replace $\beta_{k+1}$ with $T_{\alpha_{k+1}}^{-c_{k+1}}(\beta_{k+1})$, decreasing $d_m$ to $0$ without changing
the other $d_n$.  If instead $m \geq k+2$, then by an argument like in the previous paragraph
we can modify $\beta_{k+1}'$ so as to decrease $\sum_{n=k+1}^g |d_n|$, and we are done. 
\end{proof}

\begin{proof}[{Proof of Lemma \ref{lemma:modcurveslines}}]
We have a series of projections
$$\ModCurves_g \stackrel{\tilde{\pi}}{\longrightarrow} \ModCurves_g / \Torelli_{g} \stackrel{\pi}{\longrightarrow} \Lines_{\sigma,\delta}(g).$$
We must prove that for all simplices $s$ of $\Lines_{\sigma,\delta}(g)$, there is some simplex $\tilde{s}$
of $\ModCurves_g$ so that $\pi \circ \tilde{\pi} (\tilde{s}) = s$, and in addition if 
$\tilde{s}_1$ and $\tilde{s}_2$ are
simplices of $\ModCurves_g$ so that $\pi \circ \tilde{\pi} (\tilde{s}_1) = \pi \circ \tilde{\pi} (\tilde{s}_2)$, then
there is some $f \in \Torelli_{g}$ so that $f(\tilde{s}_1) = \tilde{s}_2$.  We begin with the first assertion.  Let
$s$ be a simplex of $\Lines_{\sigma,\delta}(g)$.  
There exists some simplex $\tilde{s}_0$ in $\ModCurves_{g}$ with the same dimension and type as $s$. 
Moreover, the group $\Sp_{2g}(\Z)$ acts on $\Lines_{\sigma,\delta}(g)$, and this action is clearly transitive on simplices
of the same dimension and type.  There thus exists some $f \in \Sp_{2g}(\Z)$ so that $f(\pi \circ \tilde{\pi}(\tilde{s}_0)) = s$.
Let $\tilde{f} \in \Mod_{g}$ be a mapping class that projects to $f \in \Sp_{2g}(\Z)$.  The desired simplex of $\ModCurves_{g}$
is $\tilde{s} = \tilde{f}(\tilde{s}_0)$.

We now prove the second assertion.  Let $\tilde{s}_1$ and $\tilde{s}_2$ be two simplices of $\ModCurves_g$ with
$\pi \circ \tilde{\pi} (\tilde{s}_1) = \pi \circ \tilde{\pi} (\tilde{s}_2)$.  We will do the case
that $\tilde{s}_1$ and $\tilde{s}_2$ are simplices of type $\delta$; the other cases are similar.
Let the vertices of the $\tilde{s}_i$ be
$\{\alpha_1^i,\ldots,\alpha_h^i,\alpha_{1,2}^i\}$.  Order these and pick orientations so that 
$[\alpha_j^1] = [\alpha_j^2]$ for $1 \leq j \leq h$ and so that $\alpha_1^i \cup \alpha_2^i \cup \alpha_{1,2}^i$ separates
$\Sigma_g$ into two components, one of which is homeomorphic to $\Sigma_{0,3}$.  
Set $a_j = [\alpha_j^1]$ for $1 \leq j \leq h$, and extend this to a symplectic basis $\{a_1,b_1\ldots,a_g,b_g\}$ for 
$\HH_1(\Sigma_g)$.  
For $i=1,2$, use Lemma \ref{lemma:basisrealize} to extend 
$\{\alpha_1^i,\ldots,\alpha_h^i,\alpha_{1,2}^i\}$ to a set of oriented simple closed curves
$\{\alpha_1^i,\beta_1^i,\ldots,\alpha_g^i,\beta_g^i,\alpha_{1,2}^i\}$ satisfying the conditions
of the lemma for the given symplectic basis $\{a_1,b_1,\ldots,a_g,b_g\}$.
Using the classification of surfaces, there must exist some
$f \in \Mod_{g}$ so that $f(\alpha_j^1)=\alpha_j^2$ and $f(\beta_j^1)=\beta_j^2$ for all $j$ and
so that $f(\alpha_{1,2}^1) = \alpha_{1,2}^2$.
Since we have chosen $f$ so that it fixes a basis for homology, it follows that $f \in \Torelli_{g}$.  The
proof concludes with the observation that $f(\tilde{s}_1) = \tilde{s}_2$.
\end{proof}

We conclude that the second conclusion of Proposition \ref{proposition:curvestorelliconnected} is equivalent
to the following.

\begin{proposition}
\label{proposition:linesconnected}
For $g \geq 1$, the complex $\Lines_{\sigma,\delta}(g)$ is $(g-1)$-connected.
\end{proposition}

\subsection{Skeleton of the proof of Proposition \ref{proposition:linesconnected}}
\label{section:proofskeleton}

This section is devoted to the skeleton of the proof
of Proposition \ref{proposition:linesconnected}; most of the work will
be contained in a proposition whose proof will occupy \S \ref{section:mainproposition12} - \S \ref{section:mainproposition4}.
The bulk of the proof will consist of careful modifications of spheres in the links of simplices.
To keep our modifications from getting out of hand,
we will make use of the following subcomplexes of $\Link_{\Lines_D}(\Delta)$.

\begin{definition}
For $0 \leq k \leq g$, let $\Delta^k$ be a $(k-1)$-dimensional standard
simplex of $\Lines(g)$ (when $k=0$, we interpret $\Delta^k$ as the empty set; this
is a slight abuse of notation).  We will denote by $\Lines^{\Delta^k}(g)$ the complex
$\Link_{\Lines(g)}(\Delta^k)$.  Now consider a set $\Delta' \subset (\Lines^{\Delta^k}(g))^{(0)}$.
\begin{itemize}
\item If $\Delta'$ is a simplex of type $\sigma$ in $\Lines(g)$ and $\Delta^k \cup \Delta'$
is also a simplex of type $\sigma$ in $\Lines(g)$, then we will say that $\Delta'$ is
a {\em simplex of type $\sigma$} in $\Lines^{\Delta^k}(g)$.  We remark that the key point
of this definition is that we do not allow one of the ``intersecting'' vertices of a simplex $\Delta'$
of type $\sigma$ in $\Lines^{\Delta^k}(g)$ to lie in $\Delta^k$ and the other in $\Delta'$.
\item If $\Delta^k \cup \Delta'$ is a simplex of type $\delta$ in $\Lines(g)$, let $\Span{$v_1$}$, $\Span{$v_2$}$,
and $\Span{$v_3$}$ be the vertices of $\Delta^k \cup \Delta'$ satisfying $v_3 = \pm v_1 \pm v_2$.
	\begin{itemize}
	\item If $\Span{$v_i$} \in \Delta'$ for $1 \leq i \leq 3$, then we will say that $\Delta'$ is a
{\em simplex of type $\delta_1$} in $\Lines^{\Delta^k}(g)$. 
	\item If one of the $\Span{$v_i$}$ lies in $\Delta^k$ and the other two lie in $\Delta'$, then we
will say that $\Delta'$ is a {\em simplex of type $\delta_2$} in $\Lines^{\Delta^k}(g)$.
	\item We will say that $\Delta'$ is a {\em simplex of type $\delta$} if
it is either a simplex of type $\delta_1$ or a simplex of type $\delta_2$.
        \end{itemize}
\end{itemize}
We will then denote by $\Lines_{\sigma,\delta}^{\Delta^k}(g)$ the complex $\Lines^{\Delta^k}(g)$ with 
all simplices of types $\sigma$ and $\delta$ attached.  Similarly, we will denote
by $\Lines_{\sigma}^{\Delta^k}(g)$ (resp. $\Lines_{\delta}^{\Delta^k}(g)$) the complex $\Lines^{\Delta^k}(g)$
with all simplices of type $\sigma$ (resp. $\delta$) attached.
Next, let $W$ be a submodule of $\HH_1(\Sigma_g)$.  We define $\Lines^{\Delta^k,W}(g)$ to be
the subcomplex of $\Lines^{\Delta^k}(g)$ consisting of all simplices $\{L_1,\ldots,L_k\} \in \Lines^{\Delta^k}(g)$
so that $L_i \subset W$ for all $1 \leq i \leq k$.  We define $\Lines_{\sigma,\delta}^{\Delta^k, W}(g)$, etc.\ similarly.
\end{definition}

\noindent
We can now state the following.

\begin{proposition}
\label{proposition:main}
For $g \geq 1$, let $\{a_1,b_1,\ldots,a_g,b_g\}$ be a symplectic basis for $\HH_1(\Sigma_g)$, and fix
$0 \leq k \leq g$.  Set $\Delta^k = \{\Span{$a_1$},\ldots,\Span{$a_k$}\}$ and
$W = \{\Span{$a_1,b_1,\ldots,a_{g-1},b_{g-1},a_g$}\}$.  Then the following
hold.
\begin{enumerate}
\item For $-1 \leq n \leq g-k-2$, we have $\pi_n(\Lines^{\Delta^k,W}(g))=0$.
\item For $-1 \leq n \leq g-k-2$, we have $\pi_n(\Lines^{\Delta^k}(g))=0$.
\item For $0 \leq n \leq g-k-1$, we have $\pi_n(\Lines_{\delta}^{\Delta^k,W}(g))=0$.
\item For $0 \leq n \leq g-k-1$, the map $\Lines_{\delta}^{\Delta^k}(g) \hookrightarrow \Lines_{\sigma,\delta}^{\Delta^k}(g)$ 
induces the zero map on $\pi_n$.
\end{enumerate}
\end{proposition}

\begin{remark}
The second conclusion of Proposition \ref{proposition:main} should be compared to \cite[Theorem 2.9]{CharneyVogtmann}.  We
also remark that our proof of Proposition \ref{proposition:main} is partly inspired by the
unpublished thesis of Maazen \cite{MaazenThesis}.
\end{remark}

The proof of the first and second conclusions of Proposition \ref{proposition:main} are contained in \S \ref{section:mainproposition12}, the third in \S \ref{section:mainproposition3}, and the fourth in \S \ref{section:mainproposition4}.  We
remark that conclusions one and three are used in the proofs of conclusions two and four.  Also, conclusions three
and four make strong use of the additional simplices (of type $\delta$ for conclusion three and types $\sigma$
and $\delta$ for conclusion four) -- they are precisely the reason we introduced these simplices.  We now
show that Proposition \ref{proposition:main} implies Proposition \ref{proposition:linesconnected}.

\begin{proof}[{Proof of Proposition \ref{proposition:linesconnected}}]
Fix $g \geq 1$.  We wish to show that $\pi_n(\Lines_{\sigma,\delta}(g))=0$ for $0 \leq n \leq g-1$.  
By the fourth conclusion of Proposition \ref{proposition:main},
it is enough to show that the map
$\Lines_{\delta}(g) \hookrightarrow \Lines_{\sigma,\delta}(g)$ induces a surjection on $\pi_n$ for $0 \leq n \leq g-1$.
For some $0 \leq n \leq g-1$, let $S$ be a combinatorial $n$-sphere 
and let $\phi : S \rightarrow \Lines_{\sigma,\delta}(g)$ be a simplicial map.
We must homotope $\phi$ so that $\phi(S) \subset \Lines_{\delta}(g)$.  Assume that
$e \in S^{(1)}$ is such that $\phi(e)$ is a 1-simplex of type $\sigma$.
Observe that
$$\phi(\Link_S(e)) \subset \Link_{\Lines_{\sigma,\delta}(g)}(\phi(e)) \cong \Lines(g-1).$$
Since $\Link_S(e)$ is a combinatorial $(n-2)$-sphere, the second conclusion of Proposition
\ref{proposition:main} implies that there is a combinatorial $(n-1)$-ball $B$ with
$\partial B = \Link_S(e)$ and a simplicial map $f : B \rightarrow \Link_{\Lines_{\sigma,\delta}(g)}(\phi(e))$
so that $f|_{\partial B} = \phi|_{\Link_S(e)}$.  We can thus perform a link move to $\phi$ on $e$ with $f$,
eliminating $e$.  Iterating this process, we can ensure that no simplices of $S$ are mapped to simplices
of type $\sigma$, as desired.
\end{proof}

\subsection{The proof of the first two conclusions of Proposition \ref{proposition:main}}
\label{section:mainproposition12}

We will need the following definition.

\begin{definition}
Assume that a symplectic basis $\{a_1,b_1,\ldots,a_g,b_g\}$ for $\HH_1(\Sigma_g)$ has been fixed and that
$\rho \in \{a_1,b_1,\ldots,a_g,b_g\}$.  For a 1-dimensional summand $L$ of $\HH_1(\Sigma_g)$, pick $v \in \HH_1(\Sigma_g)$
so that $L = \Span{$v$}$ (the vector $v$ is
unique up to multiplication by $\pm 1$).  Express $v$ as $\sum (c_{a_i} a_i + c_{b_i} b_i)$ with $c_{a_i},c_{b_i} \in \Z$
for $1 \leq i \leq g$.  We define the {\em $\rho$-rank} of $L$ (denoted $\Rank^{\rho}(L)$) to equal $|c_{\rho}|$.
\end{definition}

We will also need the following obvious lemma, whose proof is omitted.

\begin{lemma}
\label{lemma:changebasis}
Fix $1 \leq k < g$ and let $\Delta^k$ be a $(k-1)$-simplex in $\Lines(g)$.  Also, let
$v_1,\ldots,v_n \in \HH_1(\Sigma_g)$ be so that $\{\Span{$v_1$}, \ldots, \Span{$v_n$}\}$
is an $(n-1)$-simplex of $\Lines_{\sigma}^{\Delta^k}$.  Then for $\Span{$v$} \in \Delta^k$ and
$q_1,\ldots,q_n \in \Z$, the set $\{\Span{$v_1 + q_1 v$},\ldots,\Span{$v_n+q_n v$}\}$ is another
simplex of $\Lines_{\sigma}^{\Delta^k}(g)$ of the same type as $\{\Span{$v_1$},\ldots,\Span{$v_n$}\}$.
\end{lemma}

\begin{proof}[{Proof of Proposition \ref{proposition:main}, first conclusion}]
We must show that
$\pi_n(\Lines^{\Delta^k,W}(g)) = 0$ for $-1 \leq n \leq g-k-2$.  The proof will be by
induction on $n$.  The base case $n=-1$ is equivalent to the observation that if $k < g$, then
$\Lines^{\Delta^k,W}(g)$ is nonempty.
Assume now that $0 \leq n \leq g-k-2$ and that 
$\pi_{n'}(\Lines^{\Delta^{k'},W}(g)) = 0$ for all
$0 \leq k' < g$ and $-1 \leq n' \leq g-k'-2$ so that $n' < n$.  Let $S$ be a combinatorial $n$-sphere
and let $\phi : S \rightarrow \Lines^{\Delta^k,W}(g)$ be a simplicial map.  By Lemma \ref{lemma:simpapprox},
it is enough to show that $S$ may be homotoped to a point.

Set
$$R = \Max \{\text{$\Rank^{a_g}(\phi(x))$ $|$ $x \in S^{(0)}$}\}.$$
If $R = 0$, then $\phi(S) \subset \Star_{\Lines^{\Delta^k,W}(g)}(\Span{$a_g$})$.
Since stars are contractible, the map $\phi$ can be homotoped to a constant map.  

Assume, therefore, that $R>0$.  Let $\Delta'$ be a simplex of $S$ with 
$\Rank^{a_g}(\phi(x))=R$ for all vertices $x$ of $\Delta'$.  Choose
$\Delta'$ so that $m := \Dim(\Delta')$ is maximal, which implies that for all vertices $x$ of $\Link_S(\Delta')$,
we have $\Rank^{a_g}(\phi(x)) < R$.  Now, $\Link_S(\Delta')$ is a combinatorial $(n-m-1)$-sphere and
$\phi(\Link_S(\Delta'))$ is contained in 
$$\Link_{\Lines^{\Delta^k,W}(g)}(\phi(\Delta')) \cong \Lines^{\Delta^{k+m'},W}(g)$$ 
for some $m' \leq m$ (it may be less than $m$ if $\phi|_{\Delta'}$ is not injective).  The inductive hypothesis together
with Lemma \ref{lemma:simpapprox} therefore tells us that there a combinatorial $(n-m)$-ball $B$ with
$\partial B = \Link_S(\Delta')$ and a simplicial map $f : B \rightarrow \Link_{\Lines^{\Delta^k,W}(g)}(\phi(\Delta'))$ so that
$f|_{\partial B} = \phi|_{\Link_S(\Delta')}$.

Our goal now is to adjust $f$ so that $\Rank^{a_g}(\phi(x)) < R$ for all $x \in B^{(0)}$.
Let $\Span{v}$ be a vertex in $\phi(\Delta')$; choose $v$ so that its $a_g$-coordinate is positive.
We define a map $f' : B \rightarrow \Link_{\Lines^{\Delta^k,W}(g)}(\phi(\Delta'))$ in the following way.  Consider
$x \in B^{(0)}$, and let $v_x \in \HH_1(\Sigma_g)$ be a vector so that $f(x) = \Span{$v_x$}$.  Choose $v_x$ so
that its $a_g$-coordinate is nonnegative.  By the division algorithm, there exists a unique $q_x \in \Z$
so that $v_x + q_x v$ has a nonnegative $a_g$-coordinate and $\Rank^{a_g}(v_x + q_x v) < \Rank^{a_g}(v) = R$;
define $f'(x) = \Span{$v_x + q_x v$}$.  By Lemma \ref{lemma:changebasis}, the map
$f'$ extends to a map $f' : B \rightarrow \Link_{\Lines^{\Delta^k,W}(g)}(\phi(\Delta'))$.  Additionally,
we have that $q_x=0$ for all $x \in (\partial B)^{(0)}$ (this is where we use the maximality of $m$),
so $f'|_{\partial B} = f|_{\partial B} = \phi|_{\Link_S(\Delta')}$.

We conclude that we can perform a link move to $\phi$ that replaces $\phi|_{\Star_S(\Delta')}$ with $f'$.  Since
$\Rank^{a_g}(f'(x)) < R$ for all $x \in B$, we have removed $\Delta'$ from $S$ without introducing any
vertices whose images have $a_g$-rank greater than or equal to $R$.  Continuing in this manner allows
us to simplify $\phi$ until $R=0$, and we are done.
\end{proof}

\begin{proof}[{Proof of Proposition \ref{proposition:main}, second conclusion}]
We must show that
$\pi_n(\Lines^{\Delta^k}(g)) = 0$ for $-1 \leq n \leq g-k-2$.  The proof is nearly identical to 
the proof of the first conclusion of Proposition \ref{proposition:main} above.  The only changes needed are the following.
\begin{itemize}
\item We use the $b_g$-rank rather than the $a_g$-rank.
\item In the case $R=0$, we now have $\phi(S) \subset \Lines^{\Delta^k,W}(g)$.  We can thus
apply the first conclusion of Proposition \ref{proposition:main} to obtain the desired conclusion. \qedhere
\end{itemize}
\end{proof}

\subsection{The proof of the third conclusion of Proposition \ref{proposition:main}}
\label{section:mainproposition3}

\begin{proof}[{Proof of Proposition \ref{proposition:main}, third conclusion}]
Our goal is to prove that $\pi_n(\Lines_{\delta}^{\Delta^k,W}(g))=0$ for $0 \leq k < g$ and
$0 \leq n \leq g-k-1$.  The proof will be by induction on $n$.
The base case is $n=0$.  If $k \leq g-2$, then the first conclusion of Proposition
\ref{proposition:main} says that $\Lines^{\Delta^k,W}(g)$ is connected, and the desired result follows.
Otherwise, $k = g-1$ and we must show that $\Lines_{\delta}^{\Delta^k,W}(g)$ is connected.
An arbitrary vertex $x$ of this complex is of the form $\Span{$c_1 a_1 + \cdots + c_{g-1} a_{g-1} + a_g$}$,
where $c_i \in \Z$ for $1 \leq i \leq g-1$.  Observe that for $e = \pm 1$ and $1 \leq j \leq g-1$ the set
\begin{align*}
&\{\Span{$c_1 a_1 + \cdots + c_{g-1} a_{g-1} + a_g$}, \\
&\quad \quad \Span{$c_1 a_1 + \cdots + c_{j-1} a_{j-1} + (c_j + e) a_j + c_{j+1} a_{j+1} + \cdots + c_{g-1} a_{g-1} + a_g$}\}
\end{align*}
is an edge of type $\delta_2$; the key point is that $\Span{$a_j$} \in \Delta^k$.  
Using a sequence of such edges, we can connected $x$ to the vertex $\Span{$a_g$}$.  We conclude 
that $\Lines_{\delta}^{\Delta^k,W}(g)$ is connected, as desired.

Assume now that $1 \leq n \leq g-k-1$ and that $\pi_{n'}(\Lines_{\delta}^{\Delta^{k'},W}(g))=0$ for all
$0 \leq k' < g$ and $0 \leq n' \leq g-k'-1$ so that $n' < n$.
Let $S$ be a combinatorial $n$-sphere
and let $\phi : S \rightarrow \Lines_{\delta}^{\Delta^k,W}(g)$ be a simplicial map.  By Lemma \ref{lemma:simpapprox},
it is enough to show that $\phi$ may be homotoped to a point.

Set
$$R = \Max \{\text{$\Rank^{a_g}(\phi(x))$ $|$ $x \in S^{(0)}$}\}.$$
If $R=0$, then $\phi(S) \subset \Star_{\Lines_{\delta}^{\Delta^k,W}(g)}(\Span{$a_g$})$ (remember,
$W = \Span{$a_1,b_1,\ldots,a_{g-1},b_{g-1},a_g$}$).  Since stars are contractible,
the map $\phi$ can be homotoped to a constant map.
Assume, therefore, that $R>0$.  Our goal is to homotope $\phi$ so that
$\Rank^{a_g}(\phi(x)) < R$ for all $x \in S^{(0)}$.  Iterating this process,
we will be able to homotope $\phi$ so that $\Rank^{a_g}(\phi(x)) = 0$ for all $x \in S^{(0)}$,
as desired.  There are three steps.  

\BeginSteps
\begin{step}
We isolate vertices whose images have $a_g$-rank $R$ from the simplices whose images are of type $\delta$.  More
precisely, we will homotope $\phi$ so that if $s \in S$ is such that $\phi(s)$
is a simplex of type $\delta$, then for all vertices $x$ of $s$ we have $\Rank^{a_g}(\phi(x)) < R$.
After this homotopy, we will still have $\Rank^{a_g}(\phi(x)) \leq R$ for all $x \in S^{(0)}$.
\end{step}

We will show how to eliminate simplices that map to simplices of type $\delta_1$ containing
vertices whose $a_g$-rank is $R$; the argument that deals with simplices of type $\delta_2$ is similar
and left to the reader.  Remember that a simplex of type $\delta_1$ in $\Lines_{\delta}^{\Delta^k,W}(g)$ 
contains a unique $2$-dimensional face of type $\delta_1$.  
Let $s \in S^{(2)}$ be so that $\phi(s)$ is a simplex of type $\delta_1$.  Assume
that there is some simplex of $S$ containing $s$ as a face whose image under $\phi$
contains a vertex whose $a_g$-rank is $R$.  Next,
let $t \in S$ be a simplex of maximal dimension so that $s \subset t$ and so that
for all vertices $x$ of $t$ that do not lie in $s$, we have $\Rank^{a_g}(\phi(x)) = R$.  By
assumption, $t$ contains {\em some} vertex whose image under $\phi$ has $a_g$-rank $R$, and moreover
for all vertices $y$ of $\Link_S(t)$ we have $\Rank^{a_g}(\phi(y)) < R$.

Let $m = \Dim(t)$, and write 
$\phi(t) = \{\Span{$v_1$}, \Span{$v_2$}, \Span{$\pm v_1 \pm v_2$}, \Span{$v_4$}, \ldots, \Span{$v_{m'}$}\}$;
we may have $m' - 1 < m$ since $\phi$ need not be injective.  Now, $\Link_S(t)$ is a combinatorial $(n-m-1)$-sphere and
$\phi(\Link_S(t))$ is contained in
$$\Link_{\Lines_{\delta}^{\Delta^k,W}(g)}(\phi(t)) = \Lines^{\Delta^k \cup \{\Span{$v_1$},\Span{$v_2$},\Span{$v_4$},\ldots,\Span{$v_{m'},W$}\}}(g) \cong \Lines^{\Delta^{k+(m'-1)},W}(g).$$
Since $m'-1 \leq m$ and $n \leq g-k-1$, we have $n-m-1 \leq g-(k+m'-1)-2$.  Hence 
the first conclusion of Proposition \ref{proposition:main} together 
with Lemma \ref{lemma:simpapprox} tells us that there a combinatorial $(n-m)$-ball $B$ with
$\partial B = \Link_S(t)$ and a simplicial map $f : B \rightarrow \Link_{\Lines_{\delta}^{\Delta^k,W}(g)}(\phi(t))$ so that
$f|_{\partial B} = \phi|_{\Link_S(t)}$ and so that $f(B)$ contains no simplices of type $\delta$.  Moreover,
since $\phi(t)$ contains {\em some} vertex whose $a_g$-rank is $R$, an argument like that given
in the proof of the first and second conclusions of Proposition \ref{proposition:main} tells us that
we can assume that $\Rank^{a_g}(\phi(y)) < R$ for all vertices $y$ of $B$.  We can thus perform
a link move to $\phi$ on $t$ with $f$, eliminating $t$ while not introducing any vertices mapping
to vertices whose $a_g$-ranks are greater than or equal to $R$.  Iterating this process, we can
achieve the desired conclusion.

\begin{step}
We isolate the vertices whose images have $a_g$-rank $R$ from each other.  More precisely, we will 
homotope $\phi$ so that if $x \in S^{(0)}$ satisfies $\Rank^{a_g}(\phi(x)) = R$ and $\{x,y\} \in S^{(1)}$ is
any edge, then $\Rank^{a_g}(\phi(y)) < R$.  After this homotopy, we will still have
$\Rank^{a_g}(\phi(x)) \leq R$ for all $x \in S^{(0)}$, and moreover we will still have that
if $x \in S^{(0)}$ satisfies $\Rank^{a_g}(\phi(x)) = R$ then $\phi(\Star_S(x))$ contains 
no simplices of type $\delta$.
\end{step}

Assume that there is some simplex $s \in S$ so that $\Dim(s) \geq 1$ and $\Rank^{a_g}(\phi(x)) = R$ for
all vertices $x$ of $s$.  Choose $s$ so that $\Dim(s)$ is maximal among such simplices.  By
Step 1, for all simplices $t$ of $\Star_S(s)$ the simplex $\phi(t)$ is a standard simplex.
We will homotope $\phi$ to a new map $\phi'$ so as to remove $s$ without introducing any vertices whose images
have $a_g$-rank greater than or equal to $R$ and so as to not introduce any simplices
of type $\delta$.  Iterating this will give the desired conclusion.
There are two cases.

\BeginCases
\begin{case}
There are two vertices $x_1$ and $x_2$ of $s$ so that $\phi(x_1) \neq \phi(x_2)$.
\end{case}

Let $v_1,v_2 \in \HH_1(\Sigma_g)$ be so that $\phi(x_i) = \Span{$v_i$}$
for $1 \leq i \leq 2$; choose the $v_i$ so that their $a_g$-coordinates are positive, and hence equal
to $R$.  Let $S'$ be the result of subdividing the edge $\{x_1,x_2\}$ of $S$.  Let $x_{1,2}$ be the new
vertex.  Define $\phi' : (S')^{(0)} \rightarrow \Lines_{\delta}^{\Delta^k,W}(g)$ by the formula
$$\phi'(x) =
\begin{cases}
\Span{$v_1-v_2$} & \text{if $x = x_{1,2}$,}\\
\phi(x) & \text{otherwise.}
\end{cases}$$
We claim that $\phi'$ extends to the higher-dimensional simplices of $S'$.  Indeed, consider
$t \in S'$.  If $x_{1,2} \notin t$, then the assertion is trivial.  Otherwise, there exists
a simplex $t' \in S$ with $\{x_1,x_2\} \subset t'$ so that $t$ is one of the two simplices
that result from subdividing the edge $\{x_1,x_2\}$ of $t'$ (see Figures \ref{figure:deltasubdivisions}.a--b).
The simplex $t$ either contains $x_1$ or $x_2$.  Assume without loss of generality that it contains $x_1$.
Let the vectors $v_3,\ldots,v_l \in \HH_1(\Sigma_g)$ be so that
$\phi(t') \cup \Delta^k = \{\Span{$v_1$}, \Span{$v_2$}, \ldots, \Span{$v_l$}\}$; by assumption $\{v_1,\ldots,v_l\}$ is the basis
of an isotropic summand of $\HH_1(\Sigma_g)$.
Observe that $\phi'(t) \cup \Delta^k = \{\Span{$v_1$}, \Span{$v_1-v_2$}, \Span{$v_3$}, \ldots, \Span{$v_l$}\}$.  Since
$\{v_1, v_1-v_2, v_3, \ldots, v_l\}$ is also the basis of an isotropic summand of $\HH_1(\Sigma_g)$, it
follows that $\phi'$ extends over $t$, as desired.

\Figure{figure:deltasubdivisions}{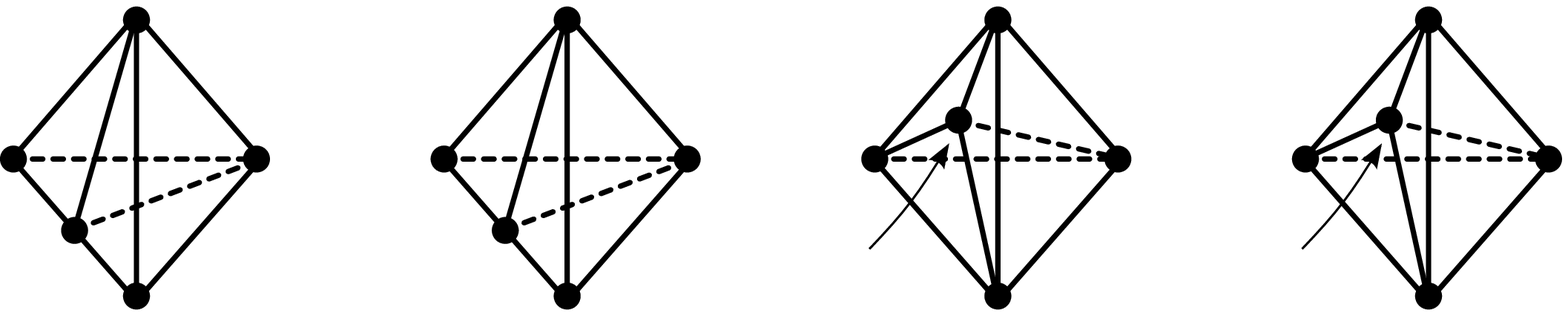}{a. Subdivided simplex $t'$ from Step 2, Case 1 \CaptionSpace b. $\phi'(t')$
\CaptionSpace c. Subdivided simplex $r'$ from Step 3 \CaptionSpace d. $g'(r')$}

Observe that $\phi$ is homotopic to $\phi'$ using simplices of type $\delta$.  Also, $x_{1,2}$ is the
only new vertex in $S'$ and $\Rank^{a_g}(\phi'(x_{1,2})) = \Rank^{a_g}(\Span{$v_1-v_2$}) = 0$; this 
calculation follows from the fact that $v_1$ and $v_2$ have the same $a_g$-coordinate.  The result
follows.

\begin{case}
For all vertices $x_1$ and $x_2$ of $s$, we have $\phi(x_1) = \phi(x_2)$.
\end{case}

Let $v \in \HH_1(\Sigma_g)$ be so that $\phi(x) = \Span{$v$}$ for all vertices $x$ of $s$.  Now,
$\Link_S(s)$ is a combinatorial $(n - \Dim(s) - 1)$-sphere and by Step 1 we have that 
$\phi(\Link_S(s))$ is contained in the following subspace of $\Link_{\Lines_{\delta}^{\Delta^k,W}(g)}(\phi(s))$ : 
$$\Lines^{\Delta^k \cup \{\Span{$v$}\},W}(g) \cong \Lines^{\Delta^{k+1},W}(g).$$
Since $\Dim(s) \geq 1$ and $n \leq g-k-1$, the dimension of $\Link_S(s)$ is at most $(g-k-1)-1-1 = g-(k+1)-2$.  The
first conclusion of Proposition \ref{proposition:main} together with
Lemma \ref{lemma:simpapprox} therefore implies that there is a combinatorial $(n-\Dim(s))$-ball $B$ with
$\partial B = \Link_S(s)$ and a simplicial map $g : B \rightarrow \Link_{\Lines_{\delta}^{\Delta^k,W}(g)}(\phi(s))$
so that $g|_{\partial B} = \phi|_{\Link_S(s)}$ and so that $g(B)$ contains no simplices of type
$\delta$.  By the maximality of the dimension of $s$, we have that
$\Rank^{a_g}(\phi(x)) < R$ for all vertices $x$ of $\Link_S(s)$, so by an argument similar to the argument in the proof
of the first and second conclusions of Proposition \ref{proposition:main}, 
we can assume that $\Rank^{a_g}(g(x)) < R$ for all vertices $x$ of
$B$.  We conclude that we can perform a link move to $\phi$ on $s$ with $g$, eliminating $s$ without
introducing any vertices whose images have $a_g$-rank greater than or equal to $R$, as desired.

\begin{step}
We eliminate all vertices whose images have $a_g$-rank $R$.  More precisely, we will homotope $\phi$ so that
for all $x \in S^{(0)}$ we have $\Rank^{a_g}(\phi(x)) < R$.  
\end{step}

Consider $x \in S^{(0)}$ so that $\Rank^{a_g}(\phi(x)) = R$.  The complex
$\Link_S(x)$ is a combinatorial $(n-1)$-sphere and by Step 2 we have $\Rank^{a_g}(\phi(y)) < R$ for
all vertices $y$ of $\Link_S(x)$.  Also, by Step 2 we have that $\phi(\Link_S(x)$ is
contained in the following subcomplex of $\Link_{\Lines_{\delta}^{\Delta^k,W}(g)}(\phi(x))$ :
$$\Lines_{\delta}^{\Delta^k \cup \{\phi(x)\},W}(g) \cong \Lines_{\delta}^{\Delta^{k+1},W}(g).$$
By induction and Lemma \ref{lemma:simpapprox}, there exists some combinatorial $n$-ball
$B$ with $\partial B = \Link_S(\{x\})$ and a simplicial map 
$g : B \rightarrow \Lines_{\delta}^{\Delta^k \cup \{\phi(x)\},W}(g)$ so that $g|_{\partial B} = \phi|_{\Link_S(\{x\})}$.
We will prove that we can modify $B$ and $g$ so that $\Rank^{a_g}(g(y)) < R$ for all $y \in B^{(0)}$.  We
will thus be able to perform a link move on $\phi$ to eliminate $x$ without introducing any vertices
whose images have $a_g$-rank greater than or equal to $R$.  Since there are no adjacent vertices in $S$ the
$a_g$-rank of whose image is equal to $R$, we can repeat this for every vertex of $S$ the $a_g$-rank of whose image
is $R$ and achieve the desired result.

For every $y \in B^{(0)}$, let $v_y \in \HH_1(\Sigma_g)$ be a vector with a nonnegative $a_g$-coordinate
so that $g(y) = \Span{$v_y$}$.  Also, let $v \in \HH_1(\Sigma_g)$ be a vector with a positive $a_g$-coordinate
so that $\phi(x) = \Span{$v$}$.  The $a_g$-coordinate of $v$ is $R$, so by the division algorithm there
exists for every $y \in B^{(0)}$ some unique $q_y \in \Z$ so that the $a_g$-coordinate of $v_y + q_y v$
is nonnegative and less than $R$.  Moreover, by assumption $q_y = 0$ for all
$y \in (\partial B)^{(0)}$.  For $y \in B^{(0)}$, define $v_y' = v_y + q_y v$ and $g'(y) = \Span{$v_y'$}$.

By Lemma \ref{lemma:changebasis}, the map $g'$ extends over all simplices
of $B$ that are mapped by $g$ to standard simplices (for later use, observe that if $g$ mapped a simplex of $B$
to a simplex of type $\sigma$, then $g'$ would extend over that simplex as well).  It will turn out that $g'$ also
extends over simplices of $B$ that are mapped by $g$ to simplices of type $\delta_2$, but
does not necessarily extend over simplices of $B$ that are mapped by $g$ to simplices of type $\delta_1$.  In
the latter case, however, we will be able to modify $B$ so as to achieve the desired extension.

We begin with the first claim, that is, that the map $g'$ extends over simplices 
$t$ of $B$ so that $g(t)$ is a simplex of type $\delta_2$ in $\Lines_{\delta}^{\Delta^k \cup \{\phi(x)\},W}(g)$.
Write $t = \{y_1,\ldots,y_l\}$, so $g(t) = \{\Span{$v_{y_1}$},\ldots,\Span{$v_{y_l}$}\}$ (since $g$ is not necessarily
injective, this latter list may have repetitions).  
Since $\Delta^k \cup \{\phi(x)\} = \{\Span{$a_1$},\ldots,\Span{$a_k$}, \Span{$v$}\}$,
after possibly reordering the $y_i$ we have the following.
\begin{itemize}
\item $v_{y_2} = v_{y_1} \pm w$ for some $w \in \{a_1,\ldots,a_k,v\}$.
\item After eliminating duplicate entries, $\{v_{y_1},v_{y_3},\ldots,v_{y_l},a_1,\ldots,a_k,v\}$ is a basis
for an isotropic summand of $\HH_1(\Sigma_g)$.
\end{itemize}
Now, clearly the set $\{v_{y_1}',v_{y_3}',\ldots,v_{y_l}',a_1,\ldots,a_k,v\}$ is also a 
basis (possibly with duplicate entries) for an isotropic summand of $\HH_1(\Sigma_g)$.
If $w \in \{a_1,\ldots,a_k\}$, then the vectors $v_{y_1}$ and $y_{y_2} = v_{y_1} \pm w$ have the same $a_g$-coordinate,
so $q_{y_1} = q_{y_2}$.  This implies that $v_{y_2}' = v_{y_1}' \pm w$, and hence $g'(t)$ is still
a simplex of type $\delta_2$.  If instead $w = v$, then $v_{y_2}' = v_{y_1}'$, so in this case
$g'(t)$ is a standard simplex.  In both cases $g'$ extends over $t$, as desired.

We conclude by showing how to modify $B$ and $g'$ so that $g'$ extends over simplices mapped by $g$ to
simplices of type $\delta_1$.  A
simplex of type $\delta_1$ has as a face a unique $2$-dimensional simplex of type $\delta_1$.  Let
$r \in B^{(2)}$ be so that $g(r) \in \Lines_{\delta}^{\Delta^k \cup \{\phi(x)\},W}(g)$ is a simplex
of type $\delta_1$.  If $r = \{z_1,z_2,z_3\}$, then by definition $v_{z_3} = \pm v_{z_1} \pm v_{z_2}$.  However,
since the $a_g$-coordinates of the $v_{z_i}$ are nonnegative, we cannot have $v_{z_3} = - v_{z_1} - v_{z_2}$.  We
conclude that after reordering the $z_i$ we can assume that $v_{z_3} = v_{z_1} + v_{z_2}$.  

Since the $a_g$-coordinates of $v_{z_1}' = v_{z_1} + q_{z_1} v$ and $v_{z_2}' = v_{z_2} + q_{z_2} v$ are nonnegative
numbers that are less than $R$, the $a_g$-coordinate of $v_{z_1} + v_{z_2} + (q_{z_1} + q_{z_2})v$ is a nonnegative
number that is less than $2R$.  Hence the $a_g$-coordinate of either $v_{z_1} + v_{z_2} + (q_{z_1} + q_{z_2})v$
or $v_{z_1} + v_{z_2} + (q_{z_1} + q_{z_2} - 1)v$ is a nonnegative number that is less than $R$.  The upshot
of all this is that either $v_{z_3}' = v_{z_1}' + v_{z_2}'$ or $v_{z_3}' = v_{z_1}' + v_{z_2}' - v$.  If 
$v_{z_3}' = v_{z_1}' + v_{z_2}'$ and if $r' \in B$ is a simplex that has $r$ as a face, then it is
clear that $g'(r')$ is a simplex of type $\delta_1$.  We can assume, therefore, that $v_{z_3}' = v_{z_1}' + v_{z_2}' - v$.

Subdivide $r$ with a new vertex $z_{z_1,z_2,z_3}$, and define $g'(z_{z_1,z_2,z_3}) = \Span{$v_{z_1}' - v$}$.  Since
the $a_g$-coordinate of $v_1' + v_2'$ is at least $R$, the $a_g$-coordinate of $v_{z_1}'$ cannot be $0$.  Hence
the $a_g$-coordinate of $v_{z_1}' - v$ is a nonpositive integer that is greater than $-R$, so
$\Rank^{a_g}(g'(z_{z_1,z_2,z_2})) < R$.

Let $r' \in B$ have $r$ as a face.  Our subdivision divides $r'$ into three simplices 
(see Figure \ref{figure:deltasubdivisions}.c), and we must check that $g'$ extends over all three of these simplices.
Write $r' = \{z_1,z_2,\ldots,z_h\}$, so $g(r') = \{\Span{$v_{z_1}$}, \ldots, \Span{$v_{z_h}$}\}$.  By
definition the set $\{v_{z_1},v_{z_2},v_{z_4},\ldots,v_{z_h},a_1,\ldots,a_k,v\}$ is a basis
for an isotropic summand of $\HH_1(\Sigma_g)$ (possibly with repetitions), so clearly 
after eliminating repetitions the set
$\{v_{z_1}',v_{z_2}',v_{z_4}',\ldots,v_{z_h}',a_1,\ldots,a_k,v\}$ is also a basis for an 
isotropic summand of $\HH_1(\Sigma_g)$.  The images under $g'$
of the three simplices that result from subdividing $r'$ are thus as follows (see Figure \ref{figure:deltasubdivisions}.d).
\begin{itemize}
\item $\{\Span{$v_{z_1}'$},\Span{$v_{z_1}'-v$},\Span{$v_{z_2}'$},\Span{$v_{z_4}'$},\ldots,\Span{$v_{z_h}'$}\}$, 
a simplex of type $\delta_2$.
\item $\{\Span{$v_{z_1}'$},\Span{$v_{z_1}'-v$},\Span{$v_{z_2}'+v_{z_1}'-v$},\Span{$v_{z_4}'$},\ldots,\Span{$v_{z_h}'$}\}$,
a simplex of type $\delta_2$.
\item $\{\Span{$v_{z_2}'+(v_{z_1}'-v)$},\Span{$v_{z_1}'-v$},\Span{$v_{z_2}'$},\Span{$v_{z_4}'$},\ldots,\Span{$v_{z_h}'$}\}$,
a simplex of type $\delta_1$.
\end{itemize}
Since $g'$ extends over all three of these, we are done
\end{proof}

\begin{remark}
For use later in the proof of the fourth conclusion
of Proposition \ref{proposition:main}, observe that
the procedure outlined in the first two steps would remain valid if we redefined $W$ to equal
$\HH_1(\Sigma_g)$; the only change needed would be to replace all references to the
first conclusion of Proposition \ref{proposition:main} with references to the
second conclusion of Proposition \ref{proposition:main}.
\end{remark}

\subsection{The proof of the fourth conclusion of Proposition \ref{proposition:main}}
\label{section:mainproposition4}

We finally come to the proof of the fourth conclusion of Proposition \ref{proposition:main}.  This
proof follows the same basic outline as the proof of \cite[Lemma 6.3]{PutmanCutPaste}, though the details
are more complicated.  To control the homotopies we construct, we will need the following definitions.

\begin{definition}
Let $\Sphere{0}$ denote the $0$-dimensional simplicial complex containing two vertices and
let $\Ball{1}$ denote the $1$-dimensional simplicial complex containing two vertices and one
edge joining those vertices.  For $n \geq 1$, the
{\em $n$-dimensional cross complex} $C_n$ (so-called because it is a subdivision of the
{\em cross polytope}; cf. \cite{CoxeterPoly}) is the join of $n-1$ copies of $\Sphere{0}$ and one copy of $\Ball{1}$.  See
Figure \ref{figure:crosscomplex}.b for pictures of $C_2$ and $C_3$.
\end{definition}

\begin{definition}
Let $0 \leq k < g$ and let $\Delta^k$ be a standard $(k-1)$-simplex of $\Lines(g)$ if $k > 0$ and $\emptyset$ if $k=0$.
\begin{itemize}
\item For $1 \leq n \leq g-k$, a {\em symplectic cross map} is a simplicial map
$\phi : C_n \rightarrow \Lines_{\sigma}^{\Delta^k}(g)$ satisfying
the following property.  Let $v_1,\ldots,v_{2n}$ be the vertices of $C_n$.  Then
there is a symplectic subspace of $\HH_1(\Sigma_g)$ with
a symplectic basis $\{a_1,b_1,\ldots,a_n,b_n\}$ so that
$\{\phi(v_1),\ldots,\phi(v_{2n})\} = \{\Span{$a_1$},\Span{$b_1$},\ldots,\Span{$a_n$},\Span{$b_n$}\}$.
\item A {\em $\sigma$-regular map} is a simplicial map $\psi : M \rightarrow \Lines_{\sigma,\delta}^{\Delta^k}(g)$,
where $M$ is a combinatorial $n$-manifold and where for all edges $e$ of $M$ so that
$\psi(e)$ is a simplex of type $\sigma$, the complex $\Star_M(e)$ is isomorphic to $C_n$ and
$\psi|_{\Star_M(e)}$ is a symplectic cross map.  Observe that this implies that $e \notin \partial M$.
\item If for $i=1,2$ we have combinatorial spheres $S_i$
and simplicial maps $f_i : S_i \rightarrow \Lines_{\delta}^{\Delta^k}(g) \subset \Lines_{\sigma,\delta}^{\Delta^k}(g)$, 
then we say
that $f_1$ and $f_2$ are {\em $\sigma$-regularly homotopic} if there is a combinatorial manifold
$A$ homeomorphic to $|S_1| \times [0,1]$ with $\partial A = S_1 \sqcup S_2$ and a $\sigma$-regular
map $\psi : A \rightarrow \Lines_{\sigma,\delta}^{\Delta^k}(g)$ so that $\psi|_{S_i} = f_i$ for $i=1,2$.
\item If $S$ is a combinatorial sphere and $f : S \rightarrow \Lines_{\delta}^{\Delta^k}(g)$ is a
simplicial map, then we say that $f$ is {\em $\sigma$-regularly nullhomotopic} if there is a combinatorial
ball $B$ with $\partial B = S$ and a $\sigma$-regular map $\psi : B \rightarrow \Lines_{\sigma,\delta}^{\Delta^k}(g)$
so that $\psi|_{S} = f$.
\end{itemize}
\end{definition}

The basic facts about $\sigma$-regularity are contained in the following lemma.

\begin{lemma}
\label{lemma:sigmaregularity}
Let $0 \leq k < g$ and let $\Delta^k$ be a standard $(k-1)$-simplex of $\Lines(g)$ if $k > 0$ and $\emptyset$ if $k=0$.
\begin{enumerate}
\item If $M$ is a combinatorial manifold and $f : M \rightarrow \Lines_{\sigma,\delta}^{\Delta^k}(g)$ is a simplicial
map so that $f(M) \subset \Lines_{\delta}^{\Delta^k}(g)$, then $f$ is $\sigma$-regular.
\item For $1 \leq i \leq 3$ let $S_i$ be a combinatorial sphere and $f_i : S_i \rightarrow \Lines_{\delta}^{\Delta^k}(g)$
be a simplicial map.  
\begin{enumerate}
\item If $f_1$ is $\sigma$-regularly homotopic to $f_2$ and $f_2$ is $\sigma$-regularly homotopic
to $f_3$, then $f_1$ is $\sigma$-regularly homotopic to $f_3$.
\item If $f_1$ is $\sigma$-regularly
homotopic to $f_2$ and $f_2$ is $\sigma$-regularly nullhomotopic, then $f_1$ is $\sigma$-regularly nullhomotopic.
\end{enumerate}
\item Let $S$ be a combinatorial $n$-sphere and let $f : S \rightarrow \Lines_{\delta}^{\Delta^k}(g)$
be a simplicial map.  Also, let $B$ be a combinatorial $(n+1)$-ball and let 
$g : B \rightarrow \Lines_{\sigma,\delta}^{\Delta^k}(g)$ be a $\sigma$-regular map.  Assume that $\partial B$
is decomposed into two combinatorial $n$-balls $D_1$ and $D_2$ so that $D_1 \cap D_2$ is a combinatorial
$(n-1)$-sphere.  Also, assume that there is a simplicial embedding $i : D_1 \hookrightarrow S$ so that
$g|_{D_1} = f \circ i$.  Define $S'$ to be $(S \setminus i(D_1 \setminus \partial D_1)) \cup_{\partial D_1} D_2$
and define $f' : S' \rightarrow \Lines_{\delta}^{\Delta^k}(g)$ to equal $f$ on $S \setminus i(D_1 \setminus \partial D_1)$
and $g$ on $D_2$.  Then $f$ is $\sigma$-regularly homotopic to $f'$.
\end{enumerate}
\end{lemma}
\begin{proof}
Conclusion 1 is trivial.  For conclusion 2.a, for $i=1,2$ let $A_i$ be a combinatorial manifold homeomorphic
to $|S_i| \times [0,1]$ with $\partial A_i = S_i \sqcup S_{i+1}$ and let 
$g_i : A_i \rightarrow \Lines_{\sigma,\delta}^{\Delta^k}(g)$ be a $\sigma$-regular map with $g_i|_{S_i} = f_i$
and $g_i|_{S_{i+1}} = f_{i+1}$.  We cannot simply glue $A_1$ to $A_2$, as the result may not be simplicial (see
Figure \ref{figure:crosscomplex}.a for an example).  Instead,
we define $A$ to be $S_2 \times \Ball{1}$ with $A_1$ and $A_2$ glued to the appropriate boundary components.  We
can then define $g : A \rightarrow \Lines_{\sigma,\delta}^{\Delta^k}(g)$ to equal $g_1$ on $A_1$, to equal the
composition of the projection $S_2 \times \Ball{1} \rightarrow S_2$ with $f_2$ on $S_2 \times \Ball{1}$, and
to equal $g_2$ on $A_2$.  It is clear that $g$ is a $\sigma$-regular map with the desired properties.  Conclusion
2.b is proven in a similar way.

For conclusion 3, define $A'$ to be $S \times \Ball{1}$ with $B$ glued to $S \times \{1\}$ along $D_1$.  It is not
hard to show that $|A'| \cong |S| \times [0,1]$.  We then define $g' : A' \rightarrow \Lines_{\sigma,\delta}^{\Delta^k}(g)$
to equal the composition of the projection $S \times \Ball{1} \rightarrow S$ with $f$ on $S \times \Ball{1}$ and $g$ on $B$.  It
is clear that $g'$ is the desired $\sigma$-regular map.
\end{proof}

\begin{proof}[{Proof of Proposition \ref{proposition:main}, fourth conclusion}]
For $0 \leq k < g$, our goal is to prove that the inclusion
map $\Lines_{\delta}^{\Delta^k}(g) \rightarrow \Lines_{\sigma,\delta}^{\Delta^k}(g)$
induces the zero map on $\pi_n$ for $0 \leq n \leq g-k-1$.  To facilitate our induction, we will
prove the stronger fact that if $S$ is a combinatorial $n$-sphere with $0 \leq n \leq g-k-1$ and
$\phi : S \rightarrow \Lines_{\sigma,\delta}^{\Delta^k}(g)$ is a simplicial map with
$\phi(S) \subset \Lines_{\delta}^{\Delta^k}(g)$, then $\phi$ is $\sigma$-regularly nullhomotopic (the
$\sigma$-regularity will be used exactly once towards the end of Step 3 of the proof below, but it is crucial --
see the comment at the end of the second paragraph of Step 3 below for a discussion of this).
The proof will be by induction on $n$.  The base case $n=0$ and the inductive cases $n \geq 1$ will
be handled simultaneously.  Thus assume that that $0 \leq n \leq g-k-1$ and that the above assertion holds 
for $n'$-spheres mapped into $\Lines_{\sigma,\delta}^{\Delta^{k'}}(g')$
for all $0 \leq k' < g'$ and $0 \leq n' \leq g'-k'-1$ so that $n' < n$.  Let $S$ be a combinatorial $n$-sphere
and let $\phi : S \rightarrow \Lines_{\sigma,\delta}^{\Delta^k}(g)$ be a simplicial map with
$\phi(S) \subset \Lines_{\delta}^{\Delta^k}(g)$.  

Set
$$R = \Max \{\text{$\Rank^{b_g}(\phi(x))$ $|$ $x \in S^{(0)}$}\}.$$
If $R=0$, then $\phi(S) \subset \Lines_{\delta}^{\Delta^k,W}(g)$ (remember,
$W = \Span{$a_1,b_1,\ldots,a_{g-1},b_{g-1},a_g$}$), and hence the
third conclusion of Proposition \ref{proposition:main} combined with Lemma
\ref{lemma:simpapprox} implies that there is a combinatorial $(n+1)$-ball $B$
with $\partial B = S$ and a simplicial map 
$$f : B \rightarrow \Lines_{\delta}^{\Delta^k,W}(g) \subset \Lines_{\sigma,\delta}^{\Delta^k}(g)$$
with $f|_{S} = \phi$.  By conclusion 1 of Lemma \ref{lemma:sigmaregularity} the map $f$ is 
$\sigma$-regular, so the conclusion follows.

Assume, therefore, that $R>0$.  Assume first that $n=0$, and let $x \in S^{(0)}$ be
so that $\Rank^{b_g}(\phi(x))=R$.  Pick $v \in \HH_1(\Sigma_g)$ so that
$\phi(x) = \Span{$v$}$.  By assumption, the set $\{a_1,\ldots,a_k,v\}$
is the basis for an isotropic summand of $\HH_1(\Sigma_g)$.  Let $v' \in \HH_1(\Sigma_g)$ satisfy $\AlgI(v,v')=1$
and $\AlgI(a_i,v')=0$ for $1 \leq i \leq k$.  Since
the $b_g$-coordinate of $v$ is $\pm R$, we can replace $v'$ with $v' + c v$ for some $c \in \Z$ if necessary
and assume that $\Rank^{b_g}(\Span{$v'$}) < R$.  Using a single simplex of type $\sigma$, we can homotope
$\phi$ so that $\phi(x) = v'$.  This homotopy is trivially $\sigma$-regular.  Iterating
this process allows us to homotope $\phi$ until the images of both vertices of $S$ have $b_g$-rank
$0$, and we are done.

Assume now that $n > 0$.  Our goal is to $\sigma$-regularly homotope $\phi$ so that
$\Rank^{b_g}(\phi(x)) < R$ for all $x \in S^{(0)}$ while retaining the property
that $\phi(S) \subset \Lines_{\delta}^{\Delta^k}(g)$ (during the intermediate steps
of this process we may introduce simplices whose images are of type $\sigma$,
but in the end we will remove them).  Using conclusion 2.a of Lemma \ref{lemma:sigmaregularity},
we can by iterating this process $\sigma$-regularly homotope $\phi$ so that $\Rank^{b_g}(\phi(x)) = 0$ for all $x \in S^{(0)}$.
An application of conclusion 2.b of Lemma \ref{lemma:sigmaregularity} then completes the proof.  
The proof will follow the same outline as the proof of the
third conclusion of Proposition \ref{proposition:main}; only the final step will
require new ideas.  Like in that proof, there are three steps.  At the end of each
of them, we will still have $\phi(S) \subset \Lines_{\delta}^{\Delta^k}(g)$.

\BeginSteps
\begin{step}
We isolate vertices whose images have $b_g$-rank $R$ from the simplices whose images are of type $\delta$.  More
precisely, we will $\sigma$-regularly homotope $\phi$ so that if $s \in S$ is such that $\phi(s)$
is a simplex of type $\delta$, then for all vertices $x$ of $S$ we have $\Rank^{b_g}(\phi(x)) < R$.
After this homotopy, we will still have $\Rank^{b_g}(\phi(x)) \leq R$ for all $x \in S^{(0)}$.
\end{step}

This is done exactly like in Step 1 of the proof of the third conclusion of Proposition
\ref{proposition:main} (see the remark following the proof of the third
conclusion of Proposition \ref{proposition:main}).
Since no simplices of type $\sigma$ are used, the first conclusion
of Lemma \ref{lemma:sigmaregularity} implies that the resulting homotopy is $\sigma$-regular.

\begin{step}
We isolate the vertices whose images have $b_g$-rank $R$ from each other.  More precisely, we will
homotope $\phi$ so that if $x \in S^{(0)}$ satisfies $\Rank^{b_g}(\phi(x)) = R$ and $\{x,y\} \in S^{(1)}$ is
any edge, then $\Rank^{b_g}(\phi(y)) < R$.  After this homotopy, we will still have
$\Rank^{b_g}(\phi(x)) \leq R$ for all $x \in S^{(0)}$, and moreover we will still have that
if $x \in S^{(0)}$ satisfies $\Rank^{b_g}(\phi(x)) = R$ then $\phi(\Star_S(x))$ contains
no simplices of type $\delta$.
\end{step}

Again, this is done exactly like in Step 2 of the proof of the third conclusion of Proposition
\ref{proposition:main}, and again no simplices of type $\sigma$ are used so the resulting
homotopy is $\sigma$-regular.

\begin{step}
We eliminate all vertices whose images have $b_g$-rank $R$.  More precisely, we will homotope $\phi$ so that
for all $x \in S^{(0)}$ we have $\Rank^{b_g}(\phi(x)) < R$.
\end{step}

\Figure{figure:crosscomplex}{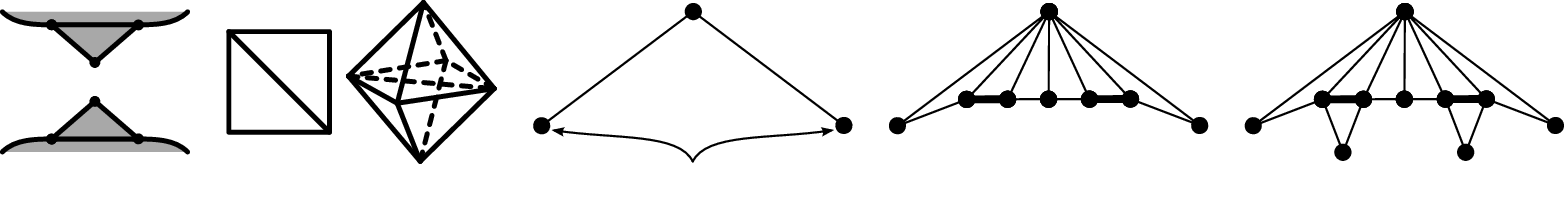}{a. Gluing the top to the bottom does not yield a simplicial
complex. \CaptionSpace b. $C_2$ and $C_3$ \CaptionSpace c. $\Link_S(x)$ is a combinatorial
$(n-1)$-sphere \CaptionSpace d. $D'$ is a combinatorial $n$-ball with $\partial D' = \Link_S(x)$.  The
edges $e_1,\ldots,e_m$ that map to edges of type $\sigma$ are in bold. \CaptionSpace e. We cone off $\Star_{D'}(e_i)$
with a new vertex $x_i$ for $1 \leq i \leq m$}

This step of the proof is illustrated in the case $n=1$ in Figures \ref{figure:crosscomplex}.c--e.
Consider $x \in S^{(0)}$ so that $\Rank^{b_g}(\phi(x)) = R$ and let $v \in \HH_1(\Sigma_g)$ be
so that $\phi(x) = \Span{$v$}$.  The complex
$\Link_S(x)$ is a combinatorial $(n-1)$-sphere and by Step 2 we have $\Rank^{b_g}(\phi(y)) < R$ for
all vertices $y$ of $\Link_S(x)$.  Our goal is to construct a combinatorial $(n+1)$-ball $B$
so that $\partial B = \Star_S(x) \cup D$ with $D$ a combinatorial $n$-ball
and $\Star_S(x) \cap D = \Link_S(x)$.  Moreover, we will also construct a $\sigma$-regular map
$g : B \rightarrow \Lines_{\sigma,\delta}^{\Delta^k}(g)$ so that $g|_{\Star_S(x)} = \phi|_{\Star_S(x)}$
and so that for all $y \in D$ we have $\Rank^{b_g}(g(y))<R$.  We can then use conclusion 3 of
Lemma \ref{lemma:sigmaregularity} to $\sigma$-regularly homotope $\phi : S \rightarrow \Link_{\sigma,\delta}^{\Delta^k}(g)$
so as to replace $\phi|_{\Star_S(x)}$ with $g|_{D}$.  This has the effect of eliminating $x$ without
introducing any vertices whose $b_g$-ranks are greater than or equal to $R$.  Iterating this
procedure will achieve the desired outcome.

As was already observed, $\Link_S(x)$ is a combinatorial $(n-1)$-sphere (see Figure \ref{figure:crosscomplex}.c).  
Also, by Step 2 we have that $\phi(\Link_S(x))$ is contained in
the following subcomplex of $\Link_{\Lines_{\sigma,\delta}^{\Delta^k}(g)}(\phi(x))$ :
$$\Lines^{\Delta^k \cup \{\phi(x)\}}(g) \cong \Lines^{\Delta^{k+1}}(g).$$
By induction, there exists some combinatorial $n$-ball
$D'$ with $\partial D' = \Link_S(x)$ and a $\sigma$-regular map
$f' : D' \rightarrow \Lines_{\sigma,\delta}^{\Delta^k \cup \{\phi(x)\}}(g)$ so that 
$f'|_{\partial D'} = \phi|_{\Link_S(\{x\})}$.
See Figure \ref{figure:crosscomplex}.d.  
Moreover, using the same argument we used in Step 3 of the proof of the third conclusion of Proposition
\ref{proposition:main} (see the parenthetical remark at the end of the first sentence of the third paragraph
of that step), we can modify $D'$ and $f'$ so that $\Rank^{b_g}(f'(y)) < R$ for all $y \in (D')^{(0)}$.
It is easy to see that these modifications do not affect the $\sigma$-regularity of $f'$.  Define $B'$
to be the join of the point $x$ with $D'$ and define $g' : B' \rightarrow \Lines_{\sigma,\delta}^{\Delta^k}(g)$
to equal $\phi$ on $x$ and $f'$ on $D'$.  It is clear that $B'$ is a combinatorial $(n+1)$-ball
and that $\partial B' = \Star_S(x) \cup D'$ with $\Star_S(x) \cap D' = \Link_S(x)$.  However, $g'$
need not be $\sigma$-regular.  In particular, $g'$ may take simplices of $D'$ to simplices of type
$\sigma$, which we wish to avoid.  The key purpose of the $\sigma$-regularity of $f'$ is
to allow us to remove these simplices of type $\sigma$.

Let $e_1,\ldots,e_m \in (D')^{(1)}$ be the edges mapping to 1-cells of type $\sigma$.  Hence for
$1 \leq i \leq m$ the complex $X_i := \Star_{D'}(e_i)$ is isomorphic to $C_n$ and $f'|_{X_i}$ is a 
symplectic cross map.  Let $V_i \subset \HH_1(\Sigma_{g})$ be the symplectic subspace of $\HH_1(\Sigma_g)$
associated to $f'|_{X_i}$ and let $\{a_1^i,b_1^i,\ldots,a_n^i,b_n^i\}$ be the associated
symplectic basis for $V_i$.  Define $W_i$ to be the orthogonal complement to $V_i$, so $W_i$ is a symplectic
subspace and we have a symplectic splitting $\HH_1(\Sigma_g) = V_i \oplus W_i$.  
Recalling that $\Delta^k = \{\Span{$a_1$},\ldots,\Span{$a_k$}\}$ and $\phi(x) = \Span{$v$}$,
we have that $\Span{$a_1,\ldots,a_k,v$}$ is an isotropic subspace of $W_i$ for each $i$.  Let $v'_i \in W_i$
be so that $\AlgI(v,v'_i) = 1$ and $\AlgI(v'_i,a_j)=0$ for $1 \leq j \leq k$.  Since $\Rank^{b_g}(\Span{$v$}) = R$, we
can replace $v'_i$ with $v'_i + c v$ for some $c \in \Z$ to ensure that $\Rank^{b_g}(\Span{$v'_i$}) < R$.  Observe
that if we set $a_{n+1}^i = v$ and $b_{n+1}^i = v_i'$, then $\{a_1^i,b_1^i,\ldots,a_{n+1}^i,b_{n+1}^i\}$
is a symplectic basis for a new symplectic subspace of $\HH_1(\Sigma_g)$.

Define $B$ to be the result of coning off the subcomplex $X_i$ of $D' \subset B'$ with a new vertex
$x_i$ for $1 \leq i \leq m$ (see Figure \ref{figure:crosscomplex}.e).  It is clear that $B$ is a
combinatorial $(n+1)$-ball and that $\partial B = \Star_S(x) \cup D$, where $D$ is the result of deleting
$X_i \setminus \partial X_i$ from $D'$ and a coning off the resulting spherical boundary component
with $x_i$ for $1 \leq i \leq m$.  Define $g : B \rightarrow \Lines_{\sigma,\delta}^{\Delta^k}(g)$
to equal $g'$ on $B'$ and to equal $\Span{$b_{n+1}^i$}$ on $x_i$.  By the previous
paragraph, $g$ is $\sigma$-regular.  Moreover, by construction we have $\Rank^{b_g}(g(y)) < R$ for all
vertices $y$ of $D$, so we are done.
\end{proof}

\noindent
Department of Mathematics; MIT, 2-306 \\
77 Massachusetts Avenue \\
Cambridge, MA 02139-4307 \\
E-mail: {\tt andyp@math.mit.edu}
\medskip

\end{document}